# Hyperseries subfields of surreal numbers

Vincent Bagayoko

IMJ-PRG

*Email:* bagayoko@imj-prg.fr

**Abstract**

We study subfields of surreal numbers, called hyperseries fields, that are suited to be equipped with derivations and composition laws. We show how to define embeddings on hyperseries fields that commute with transfinite sums and all hyperexponential and hyperlogarithmic functions.

# Introduction

Surreal numbers have many representations: Dedekind-like cuts in the set theoretic universe [14], ordinal indexed binary sequences [19], generalised power series [14, 19, 1] *à la* Hahn [20], and generalised transseries [19, 12, 13] à la Schmeling-van der Hoeven [26, 21]. Perhaps their most indicative is their recent [9, 10] presentation as hyperseries $f(\omega)$ in a variable $\omega$. Here the hyperseries $f$ is a formal series involving exponentials, logarithms, and transfinite iterates thereof [16, 17] called hyperexponentials $E_\gamma$ and hyperlogarithms $L_\gamma$ of strength $\gamma$, for arbitrary ordinals $\gamma$ (see [10, Section 7]).

A motivation for the search of representations of numbers as series in $\omega$ is the conjectured identity [23, p 6] between surreal numbers and abstract germs at infinity with similar features as germs lying in Hardy fields [25]. Those abstract germs are exactly hyperseries, provided they can be seen as infinitely differentiable monotonous functions. That is, it class for the definition of a derivation $\partial : \mathbf{No} \longrightarrow \mathbf{No}$ on the class $\mathbf{No}$ of surreal numbers and a composition law $\mathbf{No} \times \mathbf{No}^{>\mathbb{R}} \longrightarrow \mathbf{No}$ on the class of tuples $(a, b)$ with $b$ positive infinite, such that any number $a \in \mathbf{No}$, seen as the function $\tilde{a} : b \mapsto a \circ b$, should behave like a germ in a Hardy field, with derivative $\tilde{a}' = \partial(\tilde{a})$.

From the point of view of asymptotic differential algebra [3], this identity has been made precise in a series of works [2, 4, 5], based on the definition by Berarducci and Mantova [12] of a well-behaved derivation $\partial_{\mathrm{BM}} : \mathbf{No} \longrightarrow \mathbf{No}$, and culminating in the presentation of $(\mathbf{No}, +, \cdot, <, \partial_{\mathrm{BM}})$ as an elementary extension of all maximal Hardy fields.

The matters of functional composition however, both on the side of germs [8] and that of numbers, are more intricate. There is no definition of a surreal composition law $\circ : \mathbf{No} \times \mathbf{No}^{>\mathbb{R}} \longrightarrow \mathbf{No}$ which is compatible with the logarithm and exponential functions defined by Gonshor [19], let alone with their transfinite iterators. The surreal derivation $\partial_{\mathrm{BM}}$ cannot be compatible with a surreal composition law [13, Theorem 8.4]. This negative result can be interpreted as a consequence of the fact that $\partial_{\mathrm{BM}}$ is in fact incompatible with hyperexponentials and hyperlogarithms. In order to prove van der Hoeven's conjecture, it is necessary to find another definition of a derivation on $\mathbf{No}$ that is compatible with those functions, and then a composition law that enjoys similar compatibility properties. We worked on a definition of these operations (see [7, Conclusion] for a precise statement of the expected results), however our method to get there is highly technical and tedious. This article serves as a compound of tools and results required for the definition of the right derivation and composition on $\mathbf{No}$, that are designed to make their future construction easier. Our main result Theorem 2 can be seen as a partial definition of the composition.

What makes the definition of a derivation and composition law on **No** difficult is the existence of nesting in hyperserial representations of numbers, precluding simple definitions by induction on the height of expansions. A number $a$ is nested if its representation as a hyperseries is a vertically infinite right-trailing expansion

$$a = \varphi_0 + \varepsilon_0 \, \mathrm{e}^{\psi_0} (L_{\beta_0} E_{\alpha_0}(\varphi_1 + \varepsilon_1 \, \mathrm{e}^{\psi_1} (L_{\beta_1} E_{\alpha_1}(\cdots))^{\iota_1}))^{\iota_0} \qquad (1)$$

for a sequence $\Sigma = (\varphi_i, \psi_i, \varepsilon_i, \iota_i, \alpha_i, \beta_i)_{i \in \mathbb{N}}$ of parameters satisfying some technical conditions (see Section 4.2). We then say that $\Sigma$ is a nested sequence and that $a$ is $\Sigma$-nested. In [10, Section 6], we showed that each nested sequence $\Sigma$ gives rise to a proper class of $\Sigma$-nested numbers $a$, which require a nested index $z(a) \in \mathbf{No}$ in order to be distinguished. This profusion of nested numbers is not solely a complication on our way, it is also a consequence of order saturation for the surreal line, and a necessary if one hopes that many equations in first-order language of ordered differential rings with a composition law can be solved over **No**. For instance, the existence of solutions to the functional equation

$$f = \sqrt{\omega} + \mathrm{e}^{f \circ \log \omega}$$

around the approximate solution $\sqrt{\omega}$ naturally entails the existence of nested numbers with expansion

$$\sqrt{\omega} + \mathrm{e}^{\sqrt{\log \omega} + \mathrm{e}^{\sqrt{\log \log \omega} + \mathrm{e}^{\cdot^{\cdot^{\cdot}}}}}. \qquad (2)$$

In the absence of nesting, i.e. for the field of finitely nested hyperseries of [6], we already have a compatible derivation and composition law. However, the tools used to define them become insufficient for **No**, and we must combine them with a careful and fine-grained study of the hyperserial structure of the number, in particular in the nested case.

This is what we do here, by investigating the hyperserial structure of subfields $\mathbb{T}$ of **No**, called hyperseries subfields (see Definition 3.1), that are closed under hyperexponentials and hyperlogarithms of uniformly bounded strength. An embedding of a hyperseries subfield $\mathbb{T}$ into **No** is an $\mathbb{R}$-linear map $\mathbb{T} \longrightarrow \mathbf{No}$ that commutes with infinite sums, products, and hyperexponentials and hyperlogarithms, and sends monomials to monomials. In Section 4, we give relevant properties of nested numbers. In particular, we recall properties of nested sequences. Our first main result here is the definition in Section 4.3 of a fundamental tool for proving results on hyperseries subfields by induction. This is the hyperserial complexity function $\varsigma : \mathbf{No} \longrightarrow \mathbf{On}$, an ordinal measure of the complexity of numbers seen as hyperseries:

**Theorem 1.** [consequence of Theorem 4.13] *There is a unique function $\varsigma : \mathbf{No} \longrightarrow \mathbf{On}$ such that for all $a \in \mathbf{No}$, we have*

a) $\varsigma(\omega) = \omega$ *and $\varsigma(r)$ is the length of $r$ as a sign sequence [19] if $r \in \mathbb{R}$.*

b) $\varsigma(a)$ *is the ordinal sum of ordinals $\varsigma(\tau)$ where $a$ is the sum of its terms $\tau$.*

c) *If a monomial $\mathfrak{m}$ has hyperserial expansion $\mathfrak{m} = \mathrm{e}^{\psi} (L_\beta E_\alpha^u)^{\iota}$ and there is no nested sequence $\Sigma$ such that $\mathfrak{m}$ is $\Sigma$-nested, then $\varsigma(\mathfrak{m}) = \varsigma(\psi) + \varsigma(u) + \beta + 1$.*

d) *If $\Sigma = (\varphi_i, \psi_i, \ldots)_{i \in \mathbb{N}}$ is a nested sequence and $a$ is $\Sigma$-nested, then $\varsigma(a)$ is the ordinal sum of ordinals $\varsigma(\varphi_i) + \varsigma(\psi_i)$ for $i \in \mathbb{N}$.*

Thus $\varsigma$ allows one to treat $\omega$, real constants and nested numbers as elementary blocks upon which all other numbers are constructed.

In Section 5, we define nested extensions of hyperseries subfields $\mathbb{T}$ of $\mathbf{No}$, which are hyperseries subfields $\mathbb{T}_{\mathbf{P}} \supseteq \mathbb{T}$ obtained by adjoining to $\mathbb{T}$ a class $\mathbf{P}$ of $\Sigma$-nested numbers for nested sequences $\Sigma$ with parameters in $\mathbb{T}$. We prove (Proposition 5.9) that embeddings $\mathbb{T} \longrightarrow \mathbf{No}$ can be extended into embeddings $\mathbb{T}_{\mathbf{P}} \longrightarrow \mathbf{No}$. We also show (Proposition 6.2) that each hyperseries subfield is an increasing union of hyperseries subfields indexed by ordinals, where successor stages are nested and hyperexponential extensions.

Section 6 is where we prove our main general theorem (Theorem 6.3) that shows how to extend partial functions defined on hyperseries subfields $\mathbb{T}$ into embeddings $\mathbb{T} \longrightarrow \mathbf{No}$. We then apply Theorem 6.3 to the case of fields $\mathbf{No}_\lambda$ of $\lambda$-bounded numbers. Given an infinite additively indecomposable ordinal $\lambda$, the hyperseries subfield $\mathbf{No}_\lambda$ is the class of numbers whose hyperserial representation only involves ordinals below $\lambda$. For instance, $\mathbf{No}_\omega$ should be thought as the ultimate field of generalised transseries, i.e. the field of all cogent formal expressions involving exponentials and logarithms of a variable (including for instance numbers expanding as (2)). Our main concrete application is the following:

**Theorem 2.** [Section 7.5] *Let $\mathfrak{a}$ be a number such that $L_\gamma(\mathfrak{a})$ is a series of length 1 for all $\gamma < \lambda$. Then there is a unique embedding $\mathbf{No}_\lambda \longrightarrow \mathbf{No}$ that sends $\omega$ to $\mathfrak{a}$ and fixes nested indexes of nested numbers in $\mathbf{No}_\lambda$.*

Defining the composition law on $\mathbf{No}$ entails proving a similar result for any number $a > \mathbb{R}$, without the assumption on its hyperlogarithms. This is much more complicated. In a forthcoming work, we will do so in the simpler case of generalised transseries, relying the following Corollary of Theorem 2:

**Corollary 3.** *There are well-defined right compositions with $\log(\omega)$ and $\exp(\omega)$ on and onto the field $\mathbf{No}_\omega$ of generalised transseries.*

Section 7 also contains results characterising hyperexponential and nested extensions of hyperseries subfields by properties of branches in hyperserial representations of their elements. In particular, we extend to surreal numbers (Theorem 7.4) a key part of our method [6, Section 5] for constructing derivations and composition laws on finitely nested hyperseries. Indeed this is a crucial step in the definition of the composition law on $\mathbf{No}$.

# 1 Ordinals and well-based series

## 1.1 Ordinal notations

We will frequently use ordinal notations to simplify our proofs. We write $\mathbf{On}$ for the class of all ordinals and $\mathbf{On}^>$ for the class of non-zero ordinals. We write $\beta + \gamma$ and $\beta \gamma$ respectively for the ordinal sum and ordinal product of $\beta, \gamma \in \mathbf{On}$. We write $\omega^\gamma$ for the ordinal exponentiation of $\omega$ with exponent $\gamma$.

Recall that an ordinal $\alpha$ is said additively indecomposable if we have $\beta + \gamma < \alpha$ whenever $\beta, \gamma < \alpha$. Equivalently, there is $\gamma \leqslant \alpha$ with $\alpha = \omega^\gamma$.

It is sometimes practical to consider $\mathbf{On}$ itself as a generalised ordinal. Accordingly we write $\boldsymbol{\nu} \leqslant \mathbf{On}$ to say that $\boldsymbol{\nu}$ is either an ordinal or the class of all ordinals. The bold font indicates to the reader that we may allow $\boldsymbol{\nu} = \mathbf{On}$. We use the convention that $\omega^{\mathbf{On}} := \mathbf{On}$.

For $\mu \in \mathbf{On}$, we write $\mu_-$ for the unique ordinal with $\mu = \mu_- + 1$ if $\mu$ is a successor. If $\mu$ is a limit, then we set $\mu_- := \mu$. Similarly, for $\alpha := \omega^\mu$, we write $\alpha_{/\omega} := \omega^{\mu_-}$. So $\alpha = \alpha_{/\omega} \omega$ if $\mu$ is a successor, and $\alpha = \alpha_{/\omega}$ if $\mu$ is a limit.

For each ordinal $\gamma$, there is a unique family $(\gamma_\eta)_{\eta \in \mathbf{On}} \in \mathbb{N}^{\mathbf{On}}$ with such that $\gamma_\eta = 0$ for all but finitely many ordinals $\eta$, and that $\gamma = \sum_{\eta \in \mathbf{On}} \omega^\eta \gamma_\eta$. The family $(\gamma_\eta)_{\eta \in \mathbf{On}}$ is called the Cantor normal form of $\gamma$. The numbers $\eta \in \mathbf{On}$ with $\gamma_\eta \neq 0$ are called the *exponents* of $\gamma$. For $\gamma, \mu \in \mathbf{On}$, we write $\gamma \geqgg \omega^\mu$ (resp. $\gamma \ggg \omega^\mu$) if $\gamma_\eta = 0$ for all $\eta > \mu$ (resp. $\eta \geq \mu$), that is, if each exponent $\eta$ of $\gamma$ is $> \mu$ (resp. $\geq \mu$). For instance $\omega^2 + \omega 2 \geqgg \omega$ and $\omega^2 + \omega 2 \ggg 1$ but $\omega^2 + \omega 2 \not\ggg \omega$. For $\gamma, \mu \in \mathbf{On}$, there is a unique ordered pair $(\gamma', \gamma'')$ with $\gamma = \gamma' + \gamma''$ and $\gamma' \ggg \omega^\mu$ and $\gamma'' < \omega^{\mu+1}$. We have

$$\gamma' = \sum_{\eta > \mu} \omega^\eta \gamma_\eta \qquad \text{and} \qquad \gamma'' = \sum_{\eta \leq \mu} \omega^\eta \gamma_\eta.$$

## 1.2 Fields of well-based series over the reals

Let $(\mathfrak{M}, \cdot, 1, \prec)$ be a linearly ordered Abelian group. We let $\$ := \mathbb{R}[[\mathfrak{M}]]$ denote the class of functions $f : \mathfrak{M} \longrightarrow \mathbb{R}$ whose support

$$\operatorname{supp} f := \{\mathfrak{m} \in \mathfrak{M} : f(\mathfrak{m}) \neq 0\}$$

is a *well-based* set, i.e. a set which is well-ordered in the reverse order $(\mathfrak{M}, \succ)$. The elements of $\mathfrak{M}$ are called monomials, whereas those in $\mathbb{R}^\times \mathfrak{M}$ are called *terms*. We also write

$$\operatorname{term} f := \{f_\mathfrak{m} \mathfrak{m} : \mathfrak{m} \in \operatorname{supp} f\},$$

and we say that a term $\tau$ is a term in $f$ if $\tau \in \operatorname{term} f$.

We see elements $f$ of $\$$ as formal *well-based series* $f = \sum_\mathfrak{m} f_\mathfrak{m} \mathfrak{m}$ where for $\mathfrak{m} \in \mathfrak{M}$, the term $f_\mathfrak{m}$ denotes $f(\mathfrak{m}) \in \mathbb{R}$. If $\operatorname{supp} f \neq \varnothing$, then we define $\mathfrak{d}_f := \max \operatorname{supp} f \in \mathfrak{M}$ as the *dominant monomial* of $f$. For $\mathfrak{m} \in \mathfrak{M}$, we let $f_{\succ \mathfrak{m}} := \sum_{\mathfrak{n} \succ \mathfrak{m}} f_\mathfrak{n} \mathfrak{n}$ and we write $f_\succ := f_{\succ 1}$. For $f, g \in \$$, we sometimes write $f + g = f \operatorname{\dot+} g$ if $\operatorname{supp} f \succ g$. We say that a series $g \in \$$ is a *truncation* of $f$ and we write $g \trianglelefteq f$ if $\operatorname{supp}(f - g) \succ g$. The relation $\trianglelefteq$ is a well-founded large partial ordering on $\$$ with minimum $0$.

By [20], the class $\$$ is an ordered field under the pointwise sum

$$(f + g) := \sum_\mathfrak{m} (f_\mathfrak{m} + g_\mathfrak{m}) \mathfrak{m},$$

and the Cauchy product

$$fg := \sum_\mathfrak{m} \left( \sum_{\mathfrak{u}\mathfrak{v} = \mathfrak{m}} f_\mathfrak{u} g_\mathfrak{v} \right) \mathfrak{m},$$

(where each sum $\sum_{\mathfrak{u}\mathfrak{v} = \mathfrak{m}} f_\mathfrak{u} g_\mathfrak{v}$ has finite support). The positive cone $\$^> = \{f \in \$ : f > 0\}$ is

$$\$^> := \{f \in \$ : f \neq 0 \wedge f_{\mathfrak{d}_f} > 0\}.$$

The identification $f \equiv \sum_\mathfrak{m} f_\mathfrak{m} \mathfrak{m}$ induces an embedding $(\mathfrak{M}, \cdot, 1, \prec) \longrightarrow (\$^>, \cdot, 1, <)$.

The ordering $\prec$ on $\mathfrak{M}$ extends into an ordering $\prec$ on $\$$ defined by $f \prec g$ if and only if $\mathbb{R}^> |f| < |g|$. We also write $f \preccurlyeq g$ if $g \prec f$ is false, i.e. if there is $r \in \mathbb{R}^>$ with $|f| \leq r |g|$, and we write $f \asymp g$ if $f \preccurlyeq g$ and $g \preccurlyeq f$, i.e. if there is $r \in \mathbb{R}^>$ with $r |f| \geq |g|$ and $r |g| \geq f$. When $f, g$ are non-zero, we have $f \prec g$ (resp. $f \preccurlyeq g$, resp $f \asymp g$) if and only if $\mathfrak{d}_f \prec \mathfrak{d}_g$ (resp. $\mathfrak{d}_f \preccurlyeq \mathfrak{d}_g$, resp. $\mathfrak{d}_f = \mathfrak{d}_g$). We write

$$\begin{aligned}
\$_\succ &:= \{f \in \$ : \operatorname{supp} f \subseteq \mathfrak{M}^\succ\} \\
\$^\prec &:= \{f \in \$ : \operatorname{supp} f \subseteq \mathfrak{M}^\prec\} = \{f \in \$ : f \prec 1\}, \quad \text{and} \\
\$^{>,\succ} &:= \{f \in \$ : f > \mathbb{R}\} = \{f \in \$ : f \geq 0 \wedge f \succ 1\}.
\end{aligned}$$

Series in $\$_\succ$, $\$^\prec$ and $\$^{>,\succ}$ are called *purely large*, *infinitesimal*, and *positive infinite* respectively.

If $(f_i)_{i \in I}$ is a family in $\$$, then we say that $(f_i)_{i \in I}$ is *summable* if

i. $\bigcup_{i \in I} \operatorname{supp} f_i$ is well-based, and



ii. $\{i \in I : \mathfrak{m} \in \operatorname{supp} f_i\}$ is finite for all $\mathfrak{m} \in \mathfrak{M}$.

Then we may define the sum $\sum_{i \in I} f_i$ of $(f_i)_{i \in I}$ as the series

$$\sum_{i \in I} f_i := \sum_{\mathfrak{m}} \left( \sum_{i \in I} (f_i)_{\mathfrak{m}} \right) \mathfrak{m}.$$

If $\mathbb{U} = \mathbb{R}[[\mathfrak{N}]]$ is another field of well-based series and $\Psi : \mathbb{S} \longrightarrow \mathbb{U}$ is $\mathbb{R}$-linear, then we say that $\Psi$ is *strongly linear* if for every summable family $(f_i)_{i \in I}$ in $\mathbb{S}$, the family $(\Psi(f_i))_{i \in I}$ in $\mathbb{U}$ is summable, with

$$\Psi\left( \sum_{i \in I} f_i \right) = \sum_{i \in I} \Psi(f_i).$$

A map $\Phi : \mathfrak{M} \longrightarrow \mathbb{U}$ is said *well-based* if it extends into a strongly linear map $\mathbb{S} \longrightarrow \mathbb{U}$. Equivalently [22, Proposition 3.5], the family $(\Phi(\mathfrak{m}))_{\mathfrak{m} \in \mathfrak{S}}$ is summable whenever $\mathfrak{S} \subseteq \mathfrak{M}$ is a well-based subset.

### 1.3 Logarithmic hyperseries

The field $\mathbb{L}$ of logarithmic hyperseries is a field of well-based series over $\mathbb{R}$, introduced in [15], which has a structure of differential field with composition. More precisely, we have a strongly linear derivation $\partial : \mathbb{L} \longrightarrow \mathbb{L}$; $f \mapsto f'$ and a composition law $\circ : \mathbb{L} \times \mathbb{L}^{>,\succ} \longrightarrow \mathbb{L}$ with the properties listed in [15, Theorem 1.3]. We will recall here and in Remark 2.2 those w will use in this paper.

The group $\mathfrak{L}$ of monomials for $\mathbb{L}$ is the class of functions $\mathfrak{l} : \alpha \longrightarrow \mathbb{R}$ where $\alpha$ ranges in **On**, extended with zeroes on $\mathbf{On} \setminus [0, \alpha)$. The group operation is the pointwise sum and $\mathfrak{L}$ is ordered lexicographically (see [15, p 2]). We represent monomials $\mathfrak{l} \in \mathfrak{L}$ as formal products $\mathfrak{l} = \prod_{\gamma < \alpha} \ell_\gamma^{\mathfrak{l}(\gamma)}$, where each symbol $\ell_\gamma = \ell_\gamma^1, \gamma \in \mathbf{On}$ accordingly denotes the indicator function of $\{\gamma\} \subseteq \mathbf{On}$. Important features of $\partial$ and $\circ$ are that for all $\gamma, \eta \in \mathbf{On}$, we have

$$\begin{aligned} \ell'_\gamma &= \frac{1}{\prod_{\beta < \gamma} \ell_\beta}, \\ \ell_{\omega^\eta + \gamma} &= \ell_\gamma \circ \ell_{\omega^\eta} \quad \text{if } \gamma < \omega^{\eta+1}, \text{ and} \\ \ell_{\omega^{\eta+1}} \circ \ell_{\omega^\eta} &= \ell_{\omega^{\eta+1}} - 1. \end{aligned} \quad (1.1)$$

## 2 The hyperserial field of surreal numbers

We first introduce elementary properties of the class of surreal numbers. For all intents and purposes, we will consider **No** as a hyperserial field with certain properties. We purposefully omit the notions of simplicity, or birth day, which were central in previous work on surreal numbers. Thus the notions of Conway brackets, genetic definitions, equations, surreal substructures, sign sequences and so on will not appear explicitely.

### 2.1 Cuts

The linearly ordered class $(\mathbf{No}, <)$ is a class-sized saturated structure [14, 19, 18]. This means that for sets $L$ and $R$ of numbers with $L < R$, (i.e. with $\forall l \forall r ((l, r) \in L \times R \Longrightarrow l < r)$), there is $a \in \mathbf{No}$ with $L < a < R$ (i.e. with $\forall l \in L, l < a$ and $\forall r \in R, r > a$). For all subsets $L, R \subseteq \mathbf{No}$, we write $(L \mid R)$ for the class of numbers $a$ with $L < a < R$, and we call $(L \mid R)$ a *cut* in **No**. So $(L \mid R)$ is non-empty if and only if $L < R$.



**Lemma 2.1.** *Let $A, B, C, D$ be sets of numbers with $A < B$, $C < D$, and $(A \mid B) \cap (C \mid D) = \varnothing$. We have $(A \mid B) < (C \mid D)$ if and only if $B \not\succ C$.*

**Proof.** Suppose that $B \not\succ C$. So there are $b \in B$ and $c \in C$ with $b \leqslant c$. Then for $x \in (A \mid B)$, we have $x < b$ so $x < c$ so $c < (C \mid D)$. Conversely, suppose that $(A \mid B) < (C \mid D)$. So $A < D$. If $B > C$, then the number $a := \{A \cup C \mid B \cup D\}$ is well-defined and we have $a \in (A \mid B) \cap (C \mid D)$: a contradiction. So $B \not\succ C$. $\square$

See [14, Chapter 1] to understand the fundamental role which cuts play in the original definition of surreal numbers.

## 2.2 Numbers as well-based series

We recall that $(\mathbf{No}, +, \cdot, <)$ is a real-closed ordered field which canonically contains the field $\mathbb{R}$ of real numbers [14, Chapters 1 and 4]. By [14, Theorem 21], there is a subgroup $\mathbf{Mo}$ of $(\mathbf{No}^>, \times)$ such that $\mathbf{No}$ is naturally isomorphic to the field of well-based series $\mathbb{R}[[\mathbf{Mo}]]$, with which it is identified. Thus the content of Subsection 1.2 applies to $\mathbf{No}$.

## 2.3 Hyperserial structure on No

By [9, Theorem 1], there is a function $\circ : \mathbb{L} \times \mathbf{No}^{>, \succ} \longrightarrow \mathbf{No}$ which satisfies:

**C1.** For $f \in \mathbb{L}$, $g \in \mathbb{L}^{>, \succ}$ and $a \in \mathbf{No}^{>, \succ}$, we have $g \circ a \in \mathbf{No}^{>, \succ}$ and
$$f \circ (g \circ a) = (f \circ g) \circ a.$$

**C2.** For $a \in \mathbf{No}^{>, \succ}$, the function $\mathbb{L} \longrightarrow \mathbf{No}$; $f \mapsto f \circ a$ is a strongly linear morphism of ordered rings.

**C3.** For $f \in \mathbb{L}$, $a \in \mathbf{No}^{>, \succ}$ and $\delta \in \mathbf{No}$ with $\delta \prec a$, we have
$$f \circ (a + \delta) = \sum_{k \in \mathbb{N}} \frac{f^{(k)} \circ a}{k!} \delta^k.$$

**C4.** For $\gamma \in \mathbf{On}$ and $a, b \in \mathbf{No}^{>, \succ}$ with $a < b$, we have $\ell_\gamma \circ a < \ell_\gamma \circ b$.

**C5.** For $\gamma \in \mathbf{On}$ and $a \in \mathbf{No}^{>, \succ}$, there is $b \in \mathbf{No}^{>, \succ}$ with $a = \ell_\gamma \circ b$.

**Remark 2.2.** By [11, Theorems 3.1 and 4.16], the function $\circ : \mathbb{L} \times \mathbb{L}^{>, \succ} \longrightarrow \mathbb{L}$ satisfies the same properties **C1**–**C5** relativised to $\mathbb{L}$.

For $\gamma \in \mathbf{On}$, we write $L_\gamma$ for the function $\mathbf{No}^{>, \succ} \longrightarrow \mathbf{No}^{>, \succ}$; $a \mapsto \ell_\gamma \circ a$. By **C4** and **C5**, this is a strictly increasing bijection. We sometimes write $L_\gamma a := L_\gamma(a)$ for $a \in \mathbf{No}^{>, \succ}$. We write $E_\gamma$ for the functional inverse of $L_\gamma$.

For $\gamma, \rho$ with $\rho \lll \gamma$, the relation $\ell_{\gamma + \rho} = \ell_\rho \circ \ell_\gamma$ in $\mathbb{L}$, combined with **C3**, yields
$$\forall a \in \mathbf{No}^{>, \succ}, L_{\gamma + \rho} a = L_\gamma L_\rho a, \tag{2.1}$$

For $\mu \in \mathbf{On}$, the relation (1.1) in $\mathbb{L}$, combined with **C3**, yields
$$\forall a \in \mathbf{No}^{>, \succ}, L_{\omega^{\mu+1}}(L_{\omega^\mu}(a)) = L_{\omega^{\mu+1}}(a) - 1. \tag{2.2}$$

For $a \in \mathbf{No}^>$, writing $r_a$ for the coefficient of $\mathfrak{d}_a$ in $a$, we have $r_a > 0$, and there is a unique infinitesimal number $\varepsilon_a$ with $a = r_a \mathfrak{d}_a (1 + \varepsilon_a)$. Writing $\log_\mathbb{R}$ for the natural logarithm on $\mathbb{R}^> \subset \mathbf{No}$. The function defined by
$$\log a := L_1(\mathfrak{d}_a) + \log_\mathbb{R} r_a + \sum_{k \in \mathbb{N}} \frac{(-1)^k}{k+1} \varepsilon_a^{k+1}, \tag{2.3}$$



is called the *logarithm* on $\mathbf{No}^>$. This is an isomorphism $(\mathbf{No}^>, \cdot, 1, <) \longrightarrow (\mathbf{No}, +, 0, <)$ which extends $L_1$.

## 2.4 Atomicity

Given $\mu \leqslant \mathbf{On}$, we write $\mathbf{Mo}_{\omega^\mu}$ for the class of numbers $a \in \mathbf{No}^{>,\succ}$ with $L_\gamma a \in \mathbf{Mo}^\succ$ for all $\gamma < \omega^\mu$. Those numbers are called $L_{<\omega^\mu}$-atomic and they will play an important role in this paper. Indeed many properties of hyperseries can be deduced from their specialisation to $L_{<\omega^\eta}$-atomic numbers for various $\eta \in \mathbf{On}$. If $\eta = 1$, then we say that those numbers are log-atomic. We have $\mathbf{Mo}_1 = \mathbf{Mo}^\succ$ By [9, Proposition 3.21], there is a unique $L_{<\mathbf{On}}$-atomic number which is denoted $\omega$.

Let $\eta \in \mathbf{On}$. The class $\mathbf{Mo}_{\omega^\eta}$ is a proper subclass of $\mathbf{No}$ [9, Proposition 3.18]. Note that by (2.1), we have $L_{\omega^\eta} \mathbf{Mo}_{\omega^{\eta+1}} = \mathbf{Mo}_{\omega^{\eta+1}}$. Suppose that $\eta > 0$. For $a \in \mathbf{No}^{>,\succ}$, there are $\gamma < \omega^\eta$ and $\mathfrak{a} \in \mathbf{Mo}_{\omega^\eta}$ with $\delta := L_\gamma(a) - L_\gamma(\mathfrak{a}) \prec L_\gamma(a)$. Moreover, the family

$$(((\ell_{\omega^\eta}^{\uparrow \gamma})^{(k)} \circ L_\gamma(\mathfrak{a})) \delta^k)_{k \in \mathbb{N}^>}$$

is summable, and the hyperlogarithm $L_{\omega^\eta}(a)$ is given by

$$L_{\omega^\eta}(a) = L_{\omega^\eta}(\mathfrak{a}) + \sum_{k \in \mathbb{N}^>} \frac{(\ell_{\omega^\eta}^{\uparrow \gamma})^{(k)} \circ L_\gamma(\mathfrak{a})}{k!} \delta^k. \tag{2.4}$$

For $\eta \in \mathbf{On}$, we write $\mathbf{No}_{\succ,\omega^\eta} := L_{\omega^\eta}(\mathbf{Mo}_{\omega^\eta})$. We have $\mathbf{No}_{\succ,1} = \mathbf{No}_\succ \cap \mathbf{No}^{>,\succ}$. If $\varphi \in \mathbf{No}_{\succ,\omega^\eta}$ and $\delta$ is infinitesimal, with $\delta \prec (E_{\omega^\eta}^\varphi)^{-1}$, then we have [11, Lemma 6.8]

$$E_{\omega^\eta}(\varphi) \trianglelefteq E_{\omega^\eta}(\varphi + \delta). \tag{2.5}$$

Let $\eta \in \mathbf{On}$, write $\alpha := \omega^\eta$, and let $a \in \mathbf{No}^{>,\succ}$. We write

$$\mathcal{E}_\alpha[a] := \{b \in \mathbf{No} : \exists \gamma < \alpha, L_\gamma(b) \asymp L_\gamma(a)\}.$$

The class $\mathcal{E}_\alpha[a]$ is convex and the collection of classes $\mathcal{E}_\alpha[b], b \in \mathbf{No}^{>,\succ}$ forms a partition of $\mathbf{No}^{>,\succ}$. Each class $\mathcal{E}_\alpha[a]$ contains a unique $L_{<\omega^\eta}$-atomic number which is denoted $\mathfrak{d}_\alpha(a)$. The function $\mathfrak{d}_\alpha : \mathbf{No}^{>,\succ} \longrightarrow \mathbf{Mo}_\alpha$ is a non-decreasing surjection with $\mathfrak{d}_\alpha \circ \mathfrak{d}_\alpha = \mathfrak{d}_\alpha$.

We also write $\mathcal{E}_\alpha[\mathbf{S}] := \bigcup_{b \in \mathbf{S}} \mathcal{E}_\alpha[b]$ for any subclass $\mathbf{S}$ of $\mathbf{No}$.

**Lemma 2.3.** *Let $\lambda \in \omega^{\mathbf{On}^>}$ and $\mathfrak{a} \in \mathbf{Mo}_\lambda$. For $\beta < \lambda_{/\omega}$, we have $\mathfrak{d}_\lambda(L_\beta \mathfrak{a}) = \mathfrak{d}_\lambda(E_\beta \mathfrak{a}) = \mathfrak{a}$.*

**Proof.** We have $L_\beta \mathfrak{a}, E_\beta \mathfrak{a} \in \mathcal{E}_\lambda[\mathfrak{a}]$ by [9, Lemma 3.10], hence the result. □

# 3 Hyperseries subfields

We fix $\nu \leqslant \mathbf{On}$ with $\nu > 0$ and we write $\boldsymbol{\lambda} := \omega^\nu$.

## 3.1 Hyperseries subfields

**Definition 3.1.** *Let $\mathfrak{M} \subseteq \mathbf{Mo}$ be a subgroup and let $\mathbb{T} := \mathbb{R}[[\mathfrak{M}]]$. For $\nu \leqslant \mathbf{On}$, we say that $\mathbb{T}$ is a* (**hyperseries**) **subfield** *(of $\mathbf{No}$) of force $\nu$ if we have*

$$\forall \mu \leqslant \nu, \mathfrak{d}_{\omega^\mu}(\mathfrak{M}^\succ) \subseteq \mathfrak{M} \quad and \quad \mathbb{L}_{<\boldsymbol{\lambda}} \circ \mathbb{T}^{>,\succ} \subseteq \mathbb{T}.$$



We see that in particular, $\mathbb{T}$ is a confluent hyperserial field of force $\nu$ in the sense of [11, Section 6.1]. In this case, for $\mu \leqslant \nu$, we write $\mathfrak{M}_{\omega^\mu} := \mathfrak{M} \cap \mathbf{Mo}_{\omega^\mu}$. Let $\mathbf{S}$ be a subclass of $\mathbf{No}$. If $\mathbb{T}$ is smallest for the inclusion among hyperseries subfields $\mathbb{U}$ of force $\nu$ with $\mathbb{U} \supseteq \mathbf{S}$, then we write $\mathbb{T} = \mathcal{H}_\nu(\mathbf{S})$.

**Remark 3.2.** The quality of hyperseries subfield of $\mathbf{No}$ of force $\nu$ is preserved under arbitrary intersections. Consider an $\boldsymbol{\alpha} \leqslant \mathbf{On}$ with $0 < \boldsymbol{\alpha}$ and $\boldsymbol{\alpha} \notin \mathbf{On}+1$. Let $(\mathfrak{M}_\gamma)_{\gamma < \boldsymbol{\alpha}}$ be an increasing family of subgroups of $\mathbf{Mo}$ such that each $\mathbb{R}[[\mathfrak{M}_\gamma]]$ for $\gamma < \boldsymbol{\alpha}$ is a hyperseries subfield of force $\nu$. We see that $\mathbb{R}[[\bigcup_{\gamma < \boldsymbol{\alpha}} \mathfrak{M}_\gamma]]$ is a hyperseries subfield of force $\nu$.

**Proposition 3.3.** *Let* $\mathfrak{M} \subseteq \mathbf{Mo}$ *be a subgroup and assume that*

$$\begin{aligned}\mathfrak{d}_{\omega^\mu}(\mathfrak{M}^\succ) &\subseteq \mathfrak{M} \quad \textit{for all } \mu \leqslant \nu, \textit{ and} \\ \mathbb{L}_{<\omega^\mu} \circ (\mathfrak{d}_{\omega^\mu}(\mathfrak{M}^\succ)) &\subseteq \mathfrak{M} \quad \textit{for all } \mu \leqslant \nu, \textit{ and} \\ \operatorname{supp} L_{\omega^\eta}(\mathfrak{d}_{\omega^\eta}(\mathfrak{M}^\succ)) &\subseteq \mathfrak{M} \quad \textit{for all } \eta < \nu.\end{aligned}$$

*Then $\mathbb{R}[[\mathfrak{M}]]$ is a hyperseries subfield of force $\nu$.*

**Proof.** By [11, Lemma 5.9], for $a \in \mathbf{No}^{>,\succ}$ and $f \in \mathbb{L}_{<\boldsymbol{\lambda}}$, the support of $f \circ a$ is contained in the class of finite products of monomials in the class

$$\bigcup_{\mu < \nu} \bigcup_{\gamma \geqslant \omega^\mu, \gamma < \boldsymbol{\lambda}} \mathfrak{L}_{<\omega^\mu} \circ \mathfrak{d}_{\omega^\mu}(\mathfrak{d}_{L_\gamma(a)}).$$

For $\mu, \eta < \nu$ with $\eta \geqslant \mu$, we have $\mathfrak{d}_{L_{\omega^\eta}(a)} = \mathfrak{d}_{L_{\omega^\eta}(\mathfrak{d}_{\omega^\eta}(a))}$. So for $\gamma < \boldsymbol{\lambda}$ with $\gamma \geqslant \omega^\mu$, we have

$$\begin{aligned}\operatorname{supp} f \circ a &\subseteq \bigcup_{\mu < \nu} \bigcup_{\gamma \geqslant \omega^\mu, \gamma < \boldsymbol{\lambda}} \mathfrak{L}_{<\omega^\mu} \circ \mathfrak{d}_{\omega^\mu}(\mathfrak{d}_{L_\gamma(a)}) \\ &\subseteq \bigcup_{\mu < \nu} \bigcup_{\gamma \geqslant \omega^\mu, \gamma < \boldsymbol{\lambda}} \mathfrak{L}_{<\omega^\mu} \circ \mathfrak{d}_{\omega^\mu}\left(\bigcup_{\omega^\eta \not\ll \gamma} \operatorname{supp} L_{\omega^\eta}(\mathfrak{d}_{\omega^\eta}(\mathfrak{M}^\succ))\right) \\ &\subseteq \bigcup_{\mu < \nu} \mathfrak{L}_{<\omega^\mu} \circ \mathfrak{d}_{\omega^\mu}(\mathfrak{M}^\succ) \\ &\subseteq \mathfrak{M}.\end{aligned}$$

So the second condition in Definition 3.1 follows from the hypothesis, hence the result. $\square$

Write $\Lambda := \{\ell_{\boldsymbol{\lambda}_{/\omega} n} : n \in \mathbb{N}\}$ if $\nu$ is a successor and $\Lambda := \{\ell_0\}$ otherwise. Consider a non-empty subclass $\mathfrak{A} \subseteq \mathbf{Mo}_{\boldsymbol{\lambda}}$ with $\Lambda \circ \mathfrak{A} \subseteq \mathfrak{A}$. Let $\mathcal{F}_\mathfrak{A}$ denote the set of families $\mathfrak{f} := (\mathfrak{f}_\mathfrak{a})_{\mathfrak{a} \in \mathfrak{A}}$ in $\mathfrak{L}_{<\boldsymbol{\lambda}_{/\omega}}$ for which there are $n \in \mathbb{N}$ and $\mathfrak{a}_0, \ldots, \mathfrak{a}_{n-1} \in \mathfrak{A}$ with $\mathfrak{a}_0 \succ \cdots \succ \mathfrak{a}_{n-1}$ such that

$$S(\mathfrak{f}) := \{\mathfrak{a} \in \mathfrak{A} : \mathfrak{f}_\mathfrak{a} \neq 1\} \subseteq \bigcup_{i=0}^{n-1} \Lambda \circ \mathfrak{a}_i.$$

For such elements $i \in \{0, \ldots, n-1\}$, the family $(\log(\mathfrak{f}_\mathfrak{a}) \circ \mathfrak{a})_{\mathfrak{a} \in \mathfrak{A}}$ is summable. Indeed we have

$$\begin{aligned}\mathfrak{S} &:= \bigcup_{0 < i < n-1} \bigcup_{\mathfrak{a} \in \Lambda \circ \mathfrak{a}_i} \operatorname{supp}(\log(\mathfrak{f}_\mathfrak{a}) \circ \mathfrak{a}) \\ &\subseteq \bigcup_{0 < i < n-1} \{\ell_{\gamma+1} : \gamma < \boldsymbol{\lambda}\} \circ \mathfrak{a}_i,\end{aligned}$$

which is well-based. Moreover, for $\mathfrak{m} = \ell_{\boldsymbol{\lambda}_{/\omega} n + \rho} \circ \mathfrak{a}_i \in \mathfrak{S}$ where $n \in \mathbb{N}$ and $\rho < \boldsymbol{\lambda}_{/\omega}$, we have

$$\{\mathfrak{a} \in \Lambda \circ \mathfrak{a}_i : \mathfrak{m} \in \operatorname{supp} \log(\mathfrak{f}_\mathfrak{a}) \circ \mathfrak{a}\} = \{\ell_{\boldsymbol{\lambda}_{/\omega} n} \circ \mathfrak{a}_i\}.$$



So $\sum_{\mathfrak{a}\in\mathfrak{A}}(\log(\mathfrak{f}_{\mathfrak{a}})\circ\mathfrak{a})$ is defined and lies in $\mathbf{No}_{\succ}$. We set $\tilde{\mathfrak{f}}:=e^{\sum_{\mathfrak{a}\in\mathfrak{A}}\log(\mathfrak{f}_{\mathfrak{a}})\circ\mathfrak{a}}\in\mathbf{Mo}$. Write $\mathfrak{L}_{\mathfrak{A}}$ for the set of monomials $\tilde{\mathfrak{f}}$ where $\mathfrak{f}\in\mathcal{F}_{\mathfrak{A}}$. Note that this is a subgroup of $\mathbf{Mo}$.

**Lemma 3.4.** *We have $\mathcal{H}_{\boldsymbol{\nu}}(\mathfrak{A})=\mathbb{R}[[\mathfrak{L}_{\mathfrak{A}}]]$.*

**Proof.** Let $\mathbb{U}$ be a hyperseries subfield of force $\boldsymbol{\nu}$ with $\mathfrak{A}\subseteq\mathbb{U}$. We have $\mathbb{L}_{<\boldsymbol{\lambda}}\circ\mathfrak{A}\subseteq\mathbb{U}$, whence $\mathfrak{L}_{\mathfrak{A}}\subseteq\mathbb{U}$, whence $\mathbb{R}[[\mathfrak{L}_{\mathfrak{A}}]]\subseteq\mathbb{U}$. So it is enough to prove that $\mathbb{R}[[\mathfrak{L}_{\mathfrak{A}}]]$ is a hyperseries subfield of force $\boldsymbol{\nu}$. We first prove that for $\boldsymbol{\mu}\leqslant\boldsymbol{\nu}$ with $\boldsymbol{\mu}>0$, we have

$$\mathfrak{d}_{\omega^{\mu}}(\mathfrak{L}_{\mathfrak{A}})=\{\ell_{\gamma}\circ\mathfrak{a}:\mathfrak{a}\in\mathfrak{A}\wedge\gamma<\boldsymbol{\lambda}_{/\omega}\wedge\gamma\geqslant\omega^{\mu-}\}. \tag{3.1}$$

Indeed, consider $\mathfrak{f}\in\mathcal{F}_{\mathfrak{A}}$, so $\tilde{\mathfrak{f}}\in\mathfrak{L}_{\mathfrak{A}}$. Let $n\in\mathbb{N}$ and $\mathfrak{a}_0,\cdots,\mathfrak{a}_{n-1}\in\mathfrak{A}$ with $\mathfrak{a}_0\succ\cdots\succ\mathfrak{a}_{n-1}$ such that $S(\mathfrak{f})\subseteq\bigcup_{i=0}^{n-1}\Lambda\circ\mathfrak{a}_i$.

By Lemma 2.3, for $0<i<j<n-1$, we have $\mathfrak{L}_{<\boldsymbol{\lambda}_{/\omega}}^{\succ}\circ\mathfrak{a}_i\succ\mathfrak{L}_{<\boldsymbol{\lambda}_{/\omega}}^{\succ}\circ\mathfrak{a}_j$, so $\mathfrak{f}$ is uniquely determined by $\tilde{\mathfrak{f}}$. We have $\mathfrak{m}\succ 1\iff\mathfrak{f}_{\mathfrak{a}}\succ 1$ where $\mathfrak{a}:=\min S(\mathfrak{f})$. Assume now that $\mathfrak{m}\succ 1$. By the previous argument, we have $\log\mathfrak{m}\asymp\log(\mathfrak{f}_{\mathfrak{a}}\circ\mathfrak{a})\asymp\log(\ell_{\gamma}\circ\mathfrak{a})$ where $\gamma<\boldsymbol{\lambda}_{/\omega}$ is minimal such that the $\gamma$-th coefficient of $\mathfrak{f}_{\mathfrak{a}}$ is non-zero. So $\ell_{\gamma}\circ\mathfrak{a}=\mathfrak{d}_{\omega}(\mathfrak{m})$. We deduce that (3.1) holds for $\boldsymbol{\mu}=1$. Now let $\boldsymbol{\mu}\leqslant\boldsymbol{\nu}$ with $\boldsymbol{\mu}>1$ and let $\mathfrak{m}\in\mathfrak{M}^{\succ}$. We have $\mathfrak{d}_{\omega^{\mu}}(\mathfrak{m})=\mathfrak{d}_{\omega^{\mu}}(\mathfrak{d}_{\omega}(\mathfrak{m}))$ where $\mathfrak{d}_{\omega}(\mathfrak{m})\in\mathbb{L}_{<\boldsymbol{\lambda}_{/\omega}}\circ\mathfrak{a}$ for a certain $\mathfrak{a}\in\mathfrak{A}$. So the same arguments as in [11, Lemma 3.12] apply and yield (3.1) for $\boldsymbol{\mu}$.

Let $\boldsymbol{\mu}\leqslant\boldsymbol{\nu}$ and let $\mathfrak{b}\in\mathfrak{d}_{\omega^{\mu}}(\mathfrak{L}_{\mathfrak{A}})$. We will prove that $\mathfrak{L}_{<\omega^{\mu}}\circ\mathfrak{b}\subseteq\mathfrak{L}_{\mathfrak{A}}$ and that $L_{\omega^{\mu}}(\mathfrak{b})\in\mathbb{R}[[\mathfrak{L}_{\mathfrak{A}}]]$ if $\boldsymbol{\mu}<\boldsymbol{\nu}$. We write $\chi_{\mathfrak{b}}$ for the family in $\mathcal{F}_{\mathfrak{A}}$ with $\chi_{\mathfrak{b}}(\mathfrak{b})=\ell_0$ and $\chi_{\mathfrak{b}}(\mathfrak{c})=1$ if for all $\mathfrak{c}\in\mathfrak{A}$ with $\mathfrak{c}\neq\mathfrak{b}$. If $\boldsymbol{\mu}=\boldsymbol{\nu}$, then we have $\mathfrak{b}\in\mathfrak{A}$. Assume that $\boldsymbol{\nu}$ is a successor. Given $\mathfrak{l}\in\mathfrak{L}_{<\omega^{\mu}}$, write $\mathfrak{l}_n:=\prod_{\gamma<\omega^{\nu-}}\ell_{\gamma}^{\mathfrak{l}_{\omega^{\nu-}-n+\gamma}}\in\mathfrak{L}_{<\omega^{\nu-}}$. So $\mathfrak{l}=e^{\sum_{n\in\mathbb{N}}\log(\mathfrak{l}_n)\circ\ell_{\omega^{\nu-}-n}}$. Consider the family $\mathfrak{f}\in\mathcal{F}_{\mathfrak{A}}$ where for $\mathfrak{c}\in\mathfrak{A}$, we set $\mathfrak{f}_{\mathfrak{c}}:=\mathfrak{l}_n$ if $\mathfrak{c}=\ell_{\omega^{\mu-}-n}\circ\mathfrak{b}$ and $f_{\mathfrak{c}}:=1$ otherwise. We have $\tilde{\mathfrak{f}}=\mathfrak{l}\circ\mathfrak{b}$ so $\mathfrak{l}\circ\mathfrak{b}\in\mathfrak{L}_{\mathfrak{A}}$. Assume now that $\boldsymbol{\nu}$ is not a successor. So $\mathfrak{l}\in\mathfrak{L}_{<\omega^{\nu-}}$. We have $\tilde{\mathfrak{f}}=\mathfrak{l}\circ\mathfrak{b}$ where $\mathfrak{f}=(\chi_{\mathfrak{b}}(\mathfrak{c})\circ\mathfrak{l})_{\mathfrak{c}\in\mathfrak{A}}$, so $\mathfrak{l}\circ\mathfrak{b}\in\mathfrak{L}_{\mathfrak{A}}$.

If $\boldsymbol{\mu}<\boldsymbol{\nu}$, then there are $\mathfrak{a}\in\mathfrak{A}$ and $\gamma<\boldsymbol{\lambda}_{/\omega}$ with $\gamma\geqslant\omega^{\mu-}$ and $\mathfrak{b}=\ell_{\gamma}\circ\mathfrak{a}$. Write $\gamma=\gamma'+\omega^{\mu-}n$ for a certain $n\in\mathbb{N}$ and $\gamma'\gg\omega^{\mu-}$. For $\mathfrak{l}\in\mathfrak{L}_{<\omega^{\mu}}$, we have and $\mathfrak{l}\circ\mathfrak{b}=(\mathfrak{l}\circ\ell_{\gamma})\circ\mathfrak{b}$ where $\mathfrak{l}\circ\ell_{\gamma}\in\mathfrak{L}_{<\omega^{\mu}}\subseteq\mathfrak{L}_{<\omega^{\nu-}}$. So $\mathfrak{l}\circ\mathfrak{b}=\tilde{\mathfrak{f}}$ where $\mathfrak{f}=(\chi_{\mathfrak{b}}(\mathfrak{c})\circ(\mathfrak{l}\circ\ell_{\gamma}))_{\mathfrak{c}\in\mathfrak{A}}$. Moreover we have $L_{\omega^{\mu}}(\mathfrak{a})=\ell_{\gamma'+\omega^{\mu}}\circ\mathfrak{a}-n\in\mathbb{R}[[\mathfrak{L}_{\mathfrak{A}}]]$.

Finally, for $\eta<\boldsymbol{\nu}$ and $\mathfrak{m}=\ell_{\gamma}\circ\mathfrak{a}\in\mathfrak{d}_{\omega^{\eta}}(\mathfrak{L}_{\mathfrak{A}})$ where $\gamma<\boldsymbol{\lambda}_{/\omega}\wedge\gamma\geqslant\omega^{\eta-}$, We deduce with Proposition 3.3 that $\mathbb{R}[[\mathfrak{L}_{\mathfrak{A}}]]$ is a hyperseries subfield of force $\boldsymbol{\nu}$, hence $\mathbb{R}[[\mathfrak{L}_{\mathfrak{A}}]]=\mathcal{H}_{\boldsymbol{\nu}}(\mathfrak{A})$. □

**Lemma 3.5.** *Any $\boldsymbol{\lambda}$-bounded path $P$ in $\mathbb{R}[[\mathfrak{L}_{\mathfrak{A}}]]$ has length $|P|\leqslant 2$.*

**Proof.** Let $P$ be a $\boldsymbol{\lambda}$-bounded path in $\mathbb{R}[[\mathfrak{L}_{\mathfrak{A}}]]$ with $|P|>0$, so $\mathfrak{m}_{P,0}\neq 1$ and we may assume that $\mathfrak{m}_{P,0}\succ 1$. Since $P$ is $\boldsymbol{\lambda}$-bounded, we have $\mathfrak{m}_{P,i}\notin\mathfrak{A}$ for all $i\leqslant|P|$. Write $\mathfrak{m}_{P,0}=e^{\psi}(L_{\beta}E_{\alpha}^{u})^{\iota}$ in hyperserial expansion. There is $\mathfrak{f}\in\mathcal{F}_{\mathfrak{A}}$ with $\tilde{\mathfrak{f}}=\mathfrak{m}_{P,0}$. If $\alpha=0$ or $\alpha>0$ and $\sigma_{P,1}=-1$, then we have $\mathfrak{m}_{P,1}\in\mathrm{term}\log\tilde{\mathfrak{f}}$ where $\log\tilde{\mathfrak{f}}=\sum_{\mathfrak{a}\in\mathfrak{A}}\log(\mathfrak{f}_{\mathfrak{a}})\circ\mathfrak{a}$. We deduce that $\mathfrak{m}_{P,1}=L_{\gamma+1}(\mathfrak{a})$ for a certain $\gamma<\boldsymbol{\lambda}_{/\omega}$ and $\mathfrak{a}\in\mathfrak{A}$. Since $\mathfrak{a}\in\mathfrak{M}_{\boldsymbol{\lambda}}$ and $P$ is $\boldsymbol{\lambda}$-bounded, we must have $\mathfrak{a}=\omega$. So $|P|\leqslant 2$.

If $\alpha>0$ and $\sigma_{P,1}=1$, then $\mathrm{supp}\log\tilde{\mathfrak{f}}$ has a minimum $\mathfrak{n}$ with $L_{\beta}E_{\alpha}^{u}=e^{\mathfrak{n}}$. There are $\gamma<\boldsymbol{\lambda}_{/\omega}$ and $\mathfrak{a}\in\mathfrak{A}$ with $\mathfrak{n}=L_{\gamma+1}(\mathfrak{a})$, so $L_{\beta}E_{\alpha}^{u}=L_{\gamma}(\mathfrak{a})$. We deduce by the unicity of hyperserial expansions and the fact that $P$ is $\boldsymbol{\lambda}$-bounded that $u=\mathfrak{a}=\omega$. So $|P|\leqslant 2$. □

## 3.2 Hyperserial embeddings

**Definition 3.6.** *Let $\mathbb{T}=\mathbb{R}[[\mathfrak{M}]]$ be a hyperseries subfield of force $\boldsymbol{\nu}$. An **embedding of force $\boldsymbol{\nu}$** is a strongly linear, morphism of ordered rings $\Phi:\mathbb{T}\longrightarrow\mathbf{No}$ with*

   *a)* $\Phi(\mathfrak{M})\subseteq\mathbf{Mo}$.



b) $\Phi(f \circ s) = f \circ \Phi(s)$ *for all* $f \in \mathbb{L}_{<\boldsymbol{\lambda}}$ *and* $s \in \mathbb{T}^{>,\succ}$.

So $\Phi$ is a hyperserial embedding of force $\boldsymbol{\nu}$ in the sense of [11, Definition 3.4]. By [11, Proposition 3.5], we have $\Phi(\mathfrak{M}_{\omega^\mu}) = \mathfrak{M} \cap \mathbf{Mo}_{\omega^\mu}$ for all $\boldsymbol{\mu} \leqslant \boldsymbol{\lambda}$.

**Lemma 3.7.** *Let* $\mathbb{T} = \mathbb{R}[[\mathfrak{M}]]$ *be a subfield of force* $\boldsymbol{\nu}$. *Let* $\Phi : \mathbb{T} \longrightarrow \mathbf{No}$ *be a strongly linear morphism of rings with* $\Phi(L_{\omega^\eta}(\mathfrak{a})) = L_{\omega^\eta}(\Phi(\mathfrak{a}))$ *for all* $\mu < \boldsymbol{\nu}$ *and* $\mathfrak{a} \in \mathfrak{M}_{\omega^\eta}$. *We have* $\Phi(f \circ s) = f \circ \Phi(s)$ *for all* $f \in \mathbb{L}_{<\boldsymbol{\lambda}}$ *and* $s \in \mathbb{T}^{>,\succ}$.

**Proof.** Note that $\mathfrak{M} = \mathfrak{M}^{\mathbb{R}}$ is divisible, so $\mathbb{T}$ is real-closed [24]. In particular whence $\Phi$ is strictly increasing. Let $\mathbf{C}$ denote the class of series $f \in \mathbb{L}_{<\boldsymbol{\lambda}}$ with $\Phi(f \circ s) = f \circ \Phi(s)$ for all $s \in \mathbb{T}^{>,\succ}$. We prove that we have $\mathbb{L}_{<\omega^\mu} \subseteq \mathbf{C}$ by induction on $\boldsymbol{\mu} \leqslant \boldsymbol{\nu}$, starting with $\boldsymbol{\mu} = 1$.

Consider $s \in \mathbb{T}^>$ and write $s = r_s \mathfrak{d}_s (1 + \varepsilon_s)$ where $r_s \in \mathbb{R}^>$ and $\varepsilon_s \prec 1$ as in (2.3). We have $\hat{\Phi}(s) = r_s \Phi(\mathfrak{d}_s)(1 + \Phi(\varepsilon_s))$ where $\Phi(\varepsilon_s) \prec 1$, so

$$\log s = \log \mathfrak{d}_s + \log r_s + \sum_{k \in \mathbb{N}} \frac{(-1)^k}{k+1} \varepsilon_s^{k+1}, \text{ and}$$

$$\log \Phi(s) = \log \Phi(\mathfrak{d}_s) + \log r_s + \sum_{k \in \mathbb{N}} \frac{(-1)^k}{k+1} \Phi(\varepsilon_s)^{k+1}$$

$$= \hat{\Phi}(\log \mathfrak{d}_s) + \log r_s + \hat{\Phi}\left(\sum_{k \in \mathbb{N}} \frac{(-1)^k}{k+1} \varepsilon_s^{k+1}\right)$$

$$= \hat{\Phi}(\log s).$$

We deduce that $\mathbf{C}$ contains $\mathfrak{l} \in \mathfrak{L}_{<\boldsymbol{\lambda}}$ if and only if it contains $\log \mathfrak{l}$. By strong linearity of $\hat{\Phi}$, the class $\mathbf{C}$ is closed under sums of summable families. Moreover, for $f, g \in \mathbf{C}$ with $g > \mathbb{R}$, we have $f \circ g \in \mathbf{C}$. So we need only prove that we have $\ell_{\omega^\eta} \in \mathbf{C}$ for all $\eta < \boldsymbol{\nu}$. Let $\eta > 0$ such that this holds for all $\iota < \eta$. So $\mathbb{L}_{<\omega^\eta} \subseteq \mathbf{C}$ by the previous arguments. Let $s \in \mathbb{T}^{>,\succ}$ and write $\mathfrak{a} := \mathfrak{d}_{\omega^\eta(s)}$. By (2.1), there is $\gamma < \omega^\eta$ such that the number $\varepsilon := \ell_\gamma \circ s - \ell_\gamma \circ \mathfrak{a}$ is infinitesimal, with

$$\ell_{\omega^\eta} \circ s = \ell_{\omega^\eta} \circ \mathfrak{a} + \sum_{k > 0} \frac{(\ell_{\omega^\eta}^{\uparrow \gamma})^{(k)} \circ \ell_\gamma \circ \mathfrak{a}}{k!} \varepsilon^k.$$

Note that for $k \in \mathbb{N}^>$, we have $(\ell_{\omega^\eta}^{\uparrow \gamma})^{(k)} \in \mathbb{L}_{<\omega^\eta} \subseteq \mathbf{C}$. Moreover, we have

$$\ell_\gamma \circ \Phi(s) - \ell_\gamma \circ \Phi(\mathfrak{a}) = \Phi(\ell_\gamma \circ s - \ell_\gamma \circ \mathfrak{a}) = \Phi(\varepsilon) \prec 1.$$

We deduce that

$$\ell_{\omega^\eta} \circ \Phi(s) = \ell_{\omega^\eta} \circ \Phi(\mathfrak{a}) + \sum_{k > 0} \frac{(\ell_{\omega^\eta}^{\uparrow \gamma})^{(k)} \circ \ell_\gamma \circ \Phi(\mathfrak{a})}{k!} \Phi(\varepsilon)^k.$$

$$= \Phi(\ell_{\omega^\eta} \circ \mathfrak{a}) + \Phi\left(\sum_{k > 0} \frac{(\ell_{\omega^\eta}^{\uparrow \gamma})^{(k)} \circ \ell_\gamma \circ \mathfrak{a}}{k!} \varepsilon^k\right)$$

$$= \Phi(\ell_{\omega^\eta} \circ s).$$

We conclude by induction that $\mathbf{C} = \mathbb{L}_{<\boldsymbol{\lambda}}$. □

## 3.3 Hyperexponential extensions

In this subsection we consider a hyperseries subfield $\mathbb{T} = \mathbb{R}[[\mathfrak{M}]]$ of force $\boldsymbol{\nu}$ and a $\boldsymbol{\mu} \leqslant \boldsymbol{\nu}$ with $\boldsymbol{\mu} > 0$. We say that $\mathbb{T}$ has *force* $(\boldsymbol{\nu}, \boldsymbol{\mu})$ if $L_{\omega^\eta}(\mathbb{T}^{>,\succ}) = \mathbb{T}^{>,\succ}$ for all $\eta < \boldsymbol{\mu}$. For instance, the field of surreal numbers $\mathbf{No}$ has force $(\boldsymbol{\nu}, \boldsymbol{\nu})$.



In [11, Sections 6 and 7], conditions under which hyperserial fields have surjective hyperlogarithms are discussed. We recall some of the results here.

**Proposition 3.8.** [11, Corollary 7.24] *The field $\mathbb{T}$ has force $(\boldsymbol{\nu}, \boldsymbol{\mu})$ if and only if*
$$E_{\omega^\eta}(\mathbb{T} \cap \mathbf{No}_{\succ, \omega^\eta}) \subseteq \mathbb{T} \quad \text{for all } \eta < \boldsymbol{\mu}.$$

**Proposition 3.9.** [11, Theorem 7.4] *There is a smallest subfield of hyperseries $\mathbb{T}_{(<\boldsymbol{\mu})}$ of $\mathbf{No}$ of force $(\boldsymbol{\nu}, \boldsymbol{\mu})$ which contains $\mathbb{T}$.*

**Proposition 3.10.** [11, Theorem 7.4] *If $\Phi : \mathbb{T} \longrightarrow \mathbf{No}$ is an embedding of force $\boldsymbol{\nu}$, then there is a unique extension $\Phi_{(<\boldsymbol{\mu})}$ of $\Phi$ into an embedding $\mathbb{T}_{(<\boldsymbol{\mu})} \longrightarrow \mathbf{No}$ of force $\boldsymbol{\nu}$.*

Let $\eta < \mu$ and write $\mathbf{T}_\eta$ for the class of $\omega^\eta$-truncated series $\varphi$ in $\mathbb{T}$ with $E_{\omega^\eta}(\varphi) \notin \mathbb{T}$. So
$$\mathbf{T}_\eta = \mathbb{T}_{\succ, \omega^\eta} \setminus L_{\omega^\eta}(\mathfrak{M}_{\omega^\eta}).$$

For $\varphi \in \mathbf{T}_\eta$, write $\langle \varphi \rangle = \{\varphi - n : n \in \mathbb{N}\}$ if $\eta$ is a successor, and $\langle \varphi \rangle = \{\varphi\}$ otherwise. Let $\mathcal{F}_\eta$ denote the class of families $\mathfrak{f} = (\mathfrak{f}_\varphi)_{\varphi \in \mathbf{T}_\eta}$ in $(\mathfrak{L}_{<\omega^{\eta-}})^{\mathbf{T}_\eta}$ such that the set $\{\varphi \in \mathbf{T}_\eta : \mathfrak{f}_\varphi \neq 1\}$ is contained in $\bigcup_{\varphi \in S} \langle \varphi \rangle$ for a finite subset $S \subseteq \mathbf{T}_\eta$. Write $\mathfrak{L}_{<\omega^\eta}\big[e_{\omega^\eta}^{\mathbf{T}_\eta}\big]$ for the subgroup of monomials
$$\tilde{\mathfrak{f}} := e^{\sum_{\varphi \in \mathbf{T}_\eta} \log(\mathfrak{f}_\varphi) \circ E_{\omega^\eta}^\varphi},$$
where $\mathfrak{f}$ ranges in $\mathcal{F}_\eta$. We have $\mathfrak{M} \cap \mathfrak{L}_{<\omega^\eta}\big[e_{\omega^\eta}^{\mathbf{T}_\eta}\big] = \{1\}$. Write $\mathfrak{M}_{(\eta)}$ for the internal product group $\mathfrak{M} \mathfrak{L}_{<\omega^\eta}\big[e_{\omega^\eta}^{\mathbf{T}_\eta}\big]$ and set $\mathbb{T}_{(\eta)} := \mathbb{R}[[\mathfrak{M}_{(\eta)}]]$. Then $\mathbb{T}_{(\eta)}$ is a hyperseries subfield of force $\boldsymbol{\nu}$ which contains $\mathbb{T} \cup E_{\omega^\eta}^{\mathbf{T}_\eta}$.

The field $\mathbb{T}_{(<\boldsymbol{\mu})}$ can be obtained, by induction on $\boldsymbol{\mu}$, as the union of an increasing family $(\mathbb{T}_{\gamma,\boldsymbol{\mu}})_{\gamma \in \mathbf{On}}$ of subfields of hyperseries $\mathbb{T}_{\gamma,\boldsymbol{\mu}} = \mathbb{R}[[\mathfrak{M}_{\gamma,\boldsymbol{\mu}}]]$ of $\mathbf{No}$ of force $\boldsymbol{\nu}$. Indeed, for $\gamma \in \mathbf{On}$ and $\boldsymbol{\eta} \leqslant \boldsymbol{\mu}$, we set

- $\mathfrak{M}_{(0,0)} := \mathfrak{M}$.
- $\mathfrak{M}_{(\gamma,\boldsymbol{\eta})} := (\mathfrak{M}_{(\gamma,\boldsymbol{\eta}_-)})_{(\boldsymbol{\eta}_-)}$ if $\boldsymbol{\eta}$ is a successor.
- $\mathfrak{M}_{(\gamma,\boldsymbol{\eta})} := \bigcup_{\sigma < \boldsymbol{\eta}} \mathfrak{M}_{(\gamma,\sigma)}$ if $\boldsymbol{\eta}$ is a limit or $\boldsymbol{\eta} = \mathbf{On}$.
- $\mathfrak{M}_{(\gamma,0)} := \bigcup_{\rho < \gamma} \mathfrak{M}_{(\rho,\boldsymbol{\mu})}$ if $\gamma > 0$.

We set $\mathbb{T}_{(\gamma,\boldsymbol{\eta})} := \mathbb{R}[[\mathfrak{M}_{(\gamma,\boldsymbol{\eta})}]]$, so $\mathbb{T}_{(0,0)} = \mathbb{T}$ and we have the force $\boldsymbol{\nu}$ inclusion $\mathbb{T}_{(\lambda,\sigma)} \subseteq \mathbb{T}_{(\gamma,\boldsymbol{\eta})}$ whenever $\rho < \gamma$ or $\rho = \gamma$ and $\sigma \leqslant \boldsymbol{\eta}$. We set
$$\mathfrak{M}_{(<\boldsymbol{\mu})} := \bigcup_{\gamma \in \mathbf{On}} \mathfrak{M}_{(\gamma,0)}, \qquad \mathbb{T}_{(<\boldsymbol{\mu})} := \bigcup_{\gamma \in \mathbf{On}} \mathbb{T}_{(\gamma,0)}.$$

We have $\mathbb{T}_{(<\boldsymbol{\mu})} = \mathbb{R}[[\mathfrak{M}_{(<\boldsymbol{\mu})}]]$ by [11, Lemma 2.1].

**Proposition 3.11.** *Let $a \in \mathbb{T}_{(<\boldsymbol{\mu})}$ and let $P$ be an infinite path in $a$. There is $i \in \mathbb{N}$ with $\tau_{P,i} \in \mathbb{T}$.*

**Proof.** Proving the result by induction on $\boldsymbol{\mu}$, we may assume that $P$ is a path in $\mathbb{T}_{(\eta)}$ for some $\eta < \boldsymbol{\mu}$. So $\mathfrak{m}_{P,0} = \tilde{\mathfrak{f}} \mathfrak{m}$ for $\mathfrak{f} \in \mathcal{F}_\eta$ and $\mathfrak{m} \in \mathfrak{M}$.

Assume that $\alpha_{P,0} = 1$ or $\alpha_{P,0} = 0$. In the latter case we have $\sigma_{P,1} = -1$ because $P$ is infinite. So in any case we have then $\mathfrak{m}_{P,1} \in \operatorname{supp} \log \tilde{\mathfrak{f}}$ or $\mathfrak{m}_{P,1} \in \operatorname{supp} \log \mathfrak{m}$. In the latter case we are done. In the first case, we have $\mathfrak{m}_{P,1} = L_{\gamma+1} E_{\omega^\eta}^\varphi$ for a certain $\gamma < \omega^{\eta-}$ and $\varphi \in \mathbf{T}_\eta$. We have $E_{\omega^\eta}^\varphi \neq \omega$ because $P$ is infinite. So $\mathfrak{m}_{P,1} = L_{\gamma+1} E_{\omega^\eta}^\varphi$ is a hyperserial expansion. We deduce that $\tau_{P,2} \in \operatorname{term} \varphi \subseteq \mathbb{T}$.



Assume that $\alpha_{P,0} \geqslant \omega$. So $\tilde{\mathfrak{f}}\,\mathfrak{m}$ is log-atomic. We deduce since $\mathfrak{M} \cap \mathfrak{L}_{<\omega^\eta}\!\left[\mathrm{e}_{\omega^\eta}^{\mathbf{T}_\eta}\right] = \{1\}$ that $\tilde{\mathfrak{f}} = 1$ or $\mathfrak{m} = 1$. In the first case, we have $\tau_{P,0} \in \mathbb{T}$. In the second case, we have $\mathfrak{m}_{P,0} = L_\gamma E_{\omega^\eta}^\varphi$ for certain $\gamma < \omega^{\eta-}$ and $\varphi \in \mathbf{T}_\eta$, and we conclude as previously that $\tau_{P,1} \in \mathbb{T}$. □

# 4　Nested numbers

We now introduce the notions of hyperserial expansions, paths, well-nestedness, nested numbers of [10] in the more general case of hyperseries subfields. Throughout this section, we fix a $\boldsymbol{\nu} \leqslant \mathbf{On}$ with $\boldsymbol{\nu} > 0$ and a hyperseries subfield $\mathbb{T}$ of force $\boldsymbol{\nu}$, and we write $\boldsymbol{\lambda} := \omega^{\boldsymbol{\nu}}$.

## 4.1　Hyperserial expansions, paths and well-nestedness

**Definition 4.1.** *We say that a purely infinite number $\varphi \in \mathbb{T}_\succ$ is* **tail-atomic** *if $\varphi = \psi + \iota\,\mathfrak{a}$, for certain $\psi \in \mathbb{T}_\succ$, $\iota \in \{-1, 1\}$, and $\mathfrak{a} \in \mathfrak{M}_\omega$.*

**Definition 4.2.** *Let $\mathfrak{m} \in \mathfrak{M}^{\neq 1}$. Assume that there are $\psi \in \mathbb{T}_\succ$, $\iota \in \{-1,1\}$, $\alpha \in \{0\} \cup \omega^{\mathbf{On}}$, $\beta \in \mathbf{On}$ and $u \in \mathbf{No}$ such that*
$$\mathfrak{m} = \mathrm{e}^\psi \, (L_\beta E_\alpha^u)^\iota, \tag{4.1}$$
*with $\operatorname{supp} \psi \succ L_{\beta+1} E_\alpha^u$. Then we say that (4.1) is a* **hyperserial expansion of type I** *if*

- $\beta\omega < \alpha$;
- $E_\alpha^u \in \mathbf{Mo}_\alpha \setminus L_{<\alpha} \mathbf{Mo}_{\alpha\omega}$;
- $\alpha = 1 \Longrightarrow (\psi = 0$ *and $u$ is not tail-atomic*$)$.

*We say that (4.1) is a* **hyperserial expansion of type II** *if $\alpha = 0$ and $u \in \mathfrak{M}_{\boldsymbol{\lambda}}$, so that $E_\alpha^u = u$ and*
$$\mathfrak{m} = \mathrm{e}^\psi \, (L_\beta u)^\iota. \tag{4.2}$$

**Remark 4.3.** *If $\mathfrak{m} = \mathrm{e}^\psi \, (L_\beta E_\alpha^u)^\iota$ is a hyperserial expansion, then we may not have $\mathrm{e}^\psi \in \mathbb{T}$, or, equivalently, $(L_\beta E_\alpha^u)^\iota \in \mathbb{T}$. But both monomials lie in the exponential closure of $\mathbb{T}$, hence in particular in $\tilde{\mathbb{T}}$.*

**Proposition 4.4.** *Each $\mathfrak{m} \in \mathfrak{M} \setminus \{1\}$ has a unique hyperserial expansion.*

Using this, we define, as in [10, Section 5.2], the notion of path in elements of $\mathbb{T}$. Let $\rho$ be an ordinal with $0 < \rho \leqslant \omega$ and note that $i < 1 + \rho \Longleftrightarrow (i \leqslant \rho < \omega \lor i < \omega = \rho)$ for all $i \in \mathbb{N}$. Consider a sequence
$$P = (P(i))_{i<\rho} = (\tau_{P,i})_{i<\rho} = (r_{P,i}\,\mathfrak{m}_{P,i})_{i<\rho} \quad \text{in} \quad (\mathbb{R}^\times \mathfrak{M})^\rho.$$
We say that $P$ is a *path* if there exist sequences $(u_{P,i})_{i<1+\rho}$, $(\psi_{P,i})_{i<1+\rho}$, $(\iota_{P,i})_{i<\rho}$, $(\alpha_{P,i})_{i<\rho}$, and $(\beta_{P,i})_{i<1+\rho}$ with

- $u_{P,0} = \tau_{P,0}$ and $\psi_{P,0} = 0$;
- $\tau_{P,i} \in \operatorname{term} \psi_{P,i}$ or $\tau_{P,i} \in \operatorname{term} u_{P,i}$ for all $i < \rho$;
- $\tau_{P,i} \in \mathbb{R}^\times \cup \mathfrak{M}_{\boldsymbol{\lambda}} \Longrightarrow \rho = i + 1$ for all $i < \rho$;



- For $i < \rho$, the hyperserial expansion of $\mathfrak{m}_{P,i}$ is

$$\mathfrak{m}_{P,i} \;=\; \mathrm{e}^{\psi_{P,i+1}} (L_{\beta_{P,i}} E^{u_{P,i+1}}_{\alpha_{P,i}})^{\iota_{P,i}}.$$

We call $\rho$ the *length* of $P$ and we write $|P| := \rho$. We say that $P$ is *infinite* if $|P| = \omega$ and *finite* otherwise. We say that $P$ is a path in $s \in \mathbb{T}$ if $P(0)$ is a term of $s$. In that case, we set $s_{P,0} := s$. For $0 < i < |P|$, we define

$$(\sigma_{P,i}, s_{P,i}) \;:=\; \begin{cases} (-1, \psi_{P,i}) & \text{if } \mathfrak{m}_{P,i} \in \operatorname{supp} \psi_{P,i} \\ (1, u_{P,i}) & \text{if } \mathfrak{m}_{P,i} \in \operatorname{supp} u_{P,i}. \end{cases}$$

Given $\mu < \boldsymbol{\nu}$ and $\delta := \omega^\mu$, we say that a path $P$ is $\delta$-*bounded* if $\alpha_{P,i} < \delta$ for all $i < |P|$.

**Definition 4.5.** *Let $s \in \mathbb{T}$ and let $P$ be a path in $s$. We say that an index $i < |P|$ is **bad** for $(P, s)$ if one of the following conditions is satisfied*

1. *$\mathfrak{m}_{P,i}$ is not the $\preccurlyeq$-minimum of $\operatorname{supp} u_{P,i}$;*

2. *$\mathfrak{m}_{P,i} = \min \operatorname{supp} u_{P,i}$ and $\beta_{P,i} \neq 0$;*

3. *$\mathfrak{m}_{P,i} = \min \operatorname{supp} u_{P,i}$ and $\beta_{P,i} = 0$ and $r_{P,i} \notin \{-1, 1\}$;*

4. *$\mathfrak{m}_{P,i} = \min \operatorname{supp} u_{P,i}$ and $\beta_{P,i} = 0$ and $r_{P,i} \in \{-1, 1\}$ and $\mathfrak{m}_{P,i} \in \operatorname{supp} \psi_{P,i}$.*

*The index $i$ is **good** for $(P, s)$ if it is not bad for $(P, s)$.*

*If $P$ is infinite, then we say that it is **good** if $(P, \tau_{P,0})$ is good for all but a finite number of indices. In the opposite case, we say that $P$ is a **bad** path. A series $s \in \mathbb{T}$ is said to be **well-nested** every path in $s$ is good.*

## 4.2 Coding sequences

We say that a coding sequence $\Sigma$ is $\boldsymbol{\lambda}$-*bounded* if we have $\alpha_{\Sigma,i} < \boldsymbol{\lambda}$ for all $i \in \mathbb{N}$. Let $\Sigma$ be a coding sequence and let $i \in \mathbb{N}$. We write $\Sigma_{\nearrow i}$ for the coding sequence with

$$\Sigma_{\nearrow i}(j) := \Sigma(i+j) \quad \text{for all } j \in \mathbb{N}.$$

If $\alpha_{\Sigma,0} \in \omega^{\mathbf{On}+1}$, then for $n \in \mathbb{Z}$ we write $\Sigma + n$ for the coding sequence with

$$\begin{aligned} \varphi_{\Sigma+n,1} &:= \varphi_{\Sigma,1} + n, \\ \varphi_{\Sigma+n,i} &= \varphi_{\Sigma,i} \quad \text{for all } i \in \mathbb{N} \setminus \{1\}, \text{ and} \\ (\varepsilon_{\Sigma+n,i}, \psi_{\Sigma+n,i}, \iota_{\Sigma+n,i}, \alpha_{\Sigma+n,i}) &:= (\varepsilon_{\Sigma,i}, \psi_{\Sigma,i}, \iota_{\Sigma,i}, \alpha_{\Sigma,i}) \quad \text{for all } i \in \mathbb{N}. \end{aligned}$$

We write $\Sigma^+$ for the coding sequence with

$$\begin{aligned} \Sigma(0) &:= (0, 1, 0, 1, \alpha_{\Sigma,0}), \text{ and} \\ \Sigma^+(i) &:= \Sigma(i) \quad \text{for all } i \in \mathbb{N}^>. \end{aligned}$$

We say that $\Sigma$ is *positive* if $\Sigma = \Sigma^+$.

Let $\Sigma$ be a coding sequence. We say that a number $a$ is $\Sigma$-*admissible* if there is a sequence $(a_{\Sigma,i})_{i \in \mathbb{N}}$ such that for all $i \in \mathbb{N}$, we have

$$\begin{aligned} a_{\Sigma,i} &= \varphi_{\Sigma,i} \# \varepsilon_{\Sigma,i} \mathrm{e}^{\psi_{\Sigma,i}} (E_{\alpha_{\Sigma,i}} a_{\Sigma,i+1})^{\iota_i}, \\ \operatorname{supp} \psi_{\Sigma,i} &\succ \log E_{\alpha_{\Sigma,i}} a_{\Sigma,i+1}, \text{ and} \\ \varphi_{\Sigma,i+1} &\lhd \sharp_{\alpha_i}(a_{\Sigma,i+1}) \quad \text{if } \varphi_{\Sigma,i+1} \neq 0. \end{aligned}$$

We say that $\Sigma$ is *admissible* if there exists a $\Sigma$-admissible number.



We say that $a$ is $\Sigma$-*nested* if it is $\Sigma$-admissible, and if for $i \in \mathbb{N}$, the number

$$\mathfrak{m}_{a_{\Sigma,i}} := \left(\frac{a_{\Sigma,i} - \varphi_{\Sigma,i}}{r_{\Sigma,i}\,e^{\psi_{\Sigma,i}}}\right)^{\iota_{\Sigma,i}}$$

is a monomial and $\mathfrak{m}_{a_{\Sigma,i}} = E_{\alpha_{\Sigma,i}}^{a_{\Sigma,i+1}}$ is a hyperserial expansion.

**Lemma 4.6.** *If $\Sigma, \Sigma'$ are distinct coding sequences and $a$ is $\Sigma$-nested, then $a$ is not $\Sigma'$-nested.*

**Proof.** Assume for contradiction that $a$ is both $\Sigma$-nested and $\Sigma'$-nested. Consider $i \in \mathbb{N}$ minimal with $\Sigma(i) \neq \Sigma'(i)$. Taking $\Sigma_{\nearrow i}$ and $\Sigma'_{\nearrow i}$, we may assume that $i = 0$. If $\varphi_{\Sigma,0} \neq \varphi_{\Sigma',0}$, then we must have $\varphi_{\Sigma,0} \triangleleft \varphi_{\Sigma',0}$ or $\varphi_{\Sigma',0} \triangleleft \varphi_{\Sigma,0}$. But

$$\mathrm{ot}(\mathrm{supp}\,\varphi_{\Sigma,0}, \succ) + 1 = \mathrm{ot}(\mathrm{supp}\,a, \succ) = \mathrm{ot}(\mathrm{supp}\,\varphi_{\Sigma',0}, \succ) + 1,$$

so $\mathrm{ot}(\mathrm{supp}\,\varphi_{\Sigma,0}, \succ) = \mathrm{ot}(\mathrm{supp}\,\varphi_{\Sigma',0}, \succ)$: a contradiction. So $\varphi_{\Sigma,0} = \varphi_{\Sigma',0}$. We have

$$\varepsilon_{\Sigma,0} = 1 \iff a > \varphi_{\Sigma,0} \iff a > \varphi_{\Sigma',0} \iff \varepsilon_{\Sigma',0} = 1,$$

so $\varepsilon_{\Sigma,0} = \varepsilon_{\Sigma',0}$. Since $e^{\psi_{\Sigma,0}}(E_{\alpha_{\Sigma,0}}^{a_{\Sigma,1}})^{\iota_{\Sigma,0}}$ and $e^{\psi_{\Sigma',0}}(E_{\alpha_{\Sigma',0}}^{a_{\Sigma',1}})^{\iota_{\Sigma',0}}$ are both hyperserial expansions of $\varepsilon_{\Sigma,0}(a - \varphi_{\Sigma,0})$, we have $(\psi_{\Sigma,0}, \iota_{\Sigma,0}, \alpha_{\Sigma,0}) = (\psi_{\Sigma',0}, \iota_{\Sigma',0}, \alpha_{\Sigma',0})$. This contradicts the assumption that $\Sigma(0) \neq \Sigma'(0)$. $\square$

As a consequence, we may write $\Sigma_a$ for the unique coding sequence for which $a$ is $\Sigma_a$-nested. We will also simply write $a_{;i}$ instead of $a_{\Sigma,i}$ for all $i \in \mathbb{N}$, in accordance with [10, Section 6.1]. We write $\mathbf{Ad}_\Sigma$ for the class of $\Sigma$-admissible numbers, and $\mathbf{Ne}_\Sigma$ for the class of $\Sigma$-nested numbers. By [10, Proposition 6.5], there are sets $L_\Sigma$ and $R_\Sigma$ of numbers with $L_\Sigma < R_\Sigma$ and $\mathbf{Ad}_\Sigma = (L_\Sigma \mid R_\Sigma)$.

Note that we have $\mathbf{Ad}_{\Sigma_{\nearrow i}} \subseteq \varphi_{\Sigma,i} + \varepsilon_{\Sigma,i}\,e^{\psi_{\Sigma,i}}(E_{\alpha_{\Sigma,i}}\mathbf{Ad}_{\Sigma_{\nearrow i+1}})^{\iota_{\Sigma,i}}$ for all $i \in \mathbb{N}$ whenever $\Sigma$ is admissible. We say that an admissible sequence $\Sigma$ is *nested* if we have

$$\forall i \in \mathbb{N}, \mathbf{Ad}_{\Sigma_{\nearrow i}} = \varphi_{\Sigma,i} + \varepsilon_{\Sigma,i}\,e^{\psi_{\Sigma,i}}(E_{\alpha_{\Sigma,i}}\mathbf{Ad}_{\Sigma_{\nearrow i+1}})^{\iota_{\Sigma,i}}.$$

We say that a coding sequence $\Sigma$ is *good* if there is $n \in \mathbb{N}$ such that $\Sigma + n$ is nested. Note that $\Sigma_{\nearrow 1}, \Sigma - 1$ and $\Sigma^+$ are good if $\Sigma$ is good. If $\Sigma$ is good, then we generalize some of our notations with

$$L_\Sigma := \varphi_{\Sigma,0} + \varepsilon_{\Sigma,0}\,e^{\psi_{\Sigma,0}}\left(L_{(\alpha_{\Sigma,0})/\omega}n\left(\frac{L_{(\Sigma+n)} - \varphi_{\Sigma,0}}{\varepsilon_{\Sigma,0}\,e^{\psi_{\Sigma,0}}}\right)^{\iota_{\Sigma,0}}\right)^{\iota_{\Sigma,0}} \quad \text{if } \varepsilon_{\Sigma,0}\,\iota_{\Sigma,0} = 1,$$

$$L_\Sigma := \varphi_{\Sigma,0} + \varepsilon_{\Sigma,0}\,e^{\psi_{\Sigma,0}}\left(L_{(\alpha_{\Sigma,0})/\omega}n\left(\frac{R_{(\Sigma+n)} - \varphi_{\Sigma,0}}{\varepsilon_{\Sigma,0}\,e^{\psi_{\Sigma,0}}}\right)^{\iota_{\Sigma,0}}\right)^{\iota_{\Sigma,0}} \quad \text{if } \varepsilon_{\Sigma,0}\,\iota_{\Sigma,0} = -1,$$

$$R_\Sigma := \varphi_{\Sigma,0} + \varepsilon_{\Sigma,0}\,e^{\psi_{\Sigma,0}}\left(L_{(\alpha_{\Sigma,0})/\omega}n\left(\frac{R_{(\Sigma+n)} - \varphi_{\Sigma,0}}{\varepsilon_{\Sigma,0}\,e^{\psi_{\Sigma,0}}}\right)^{\iota_{\Sigma,0}}\right)^{\iota_{\Sigma,0}} \quad \text{if } \varepsilon_{\Sigma,0}\,\iota_{\Sigma,0} = 1,$$

$$R_\Sigma := \varphi_{\Sigma,0} + \varepsilon_{\Sigma,0}\,e^{\psi_{\Sigma,0}}\left(L_{(\alpha_{\Sigma,0})/\omega}n\left(\frac{L_{(\Sigma+n)} - \varphi_{\Sigma,0}}{\varepsilon_{\Sigma,0}\,e^{\psi_{\Sigma,0}}}\right)^{\iota_{\Sigma,0}}\right)^{\iota_{\Sigma,0}} \quad \text{if } \varepsilon_{\Sigma,0}\,\iota_{\Sigma,0} = -1,$$

$$\mathbf{Ad}_\Sigma := \varphi_{\Sigma,0} + \varepsilon_{\Sigma,0}\,e^{\psi_{\Sigma,0}}\left(L_{(\alpha_{\Sigma,0})/\omega}n\left(\frac{\mathbf{Ad}_{(\Sigma+n)} - \varphi_{\Sigma,0}}{\varepsilon_{\Sigma,0}\,e^{\psi_{\Sigma,0}}}\right)^{\iota_{\Sigma,0}}\right)^{\iota_{\Sigma,0}} = (L_\Sigma \mid R_\Sigma), \quad \text{and}$$

$$\mathbf{Ne}_\Sigma := \varphi_{\Sigma,0} + \varepsilon_{\Sigma,0}\,e^{\psi_{\Sigma,0}}\left(L_{(\alpha_{\Sigma,0})/\omega}n\left(\frac{\mathbf{Ne}_{(\Sigma+n)} - \varphi_{\Sigma,0}}{\varepsilon_{\Sigma,0}\,e^{\psi_{\Sigma,0}}}\right)^{\iota_{\Sigma,0}}\right)^{\iota_{\Sigma,0}}.$$



Good sequences will play an important role in the sequel. The following useful properties are easily deduced from the corresponding properties of nested sequences.

**Lemma 4.7.** [10, Lemma 6.14] *Let $\Sigma$ be a good sequence. We have*

$$\mathbf{Ad}_\Sigma = \varphi_{\Sigma,0} + \varepsilon_{\Sigma,0}\, \mathrm{e}^{\psi_{\Sigma,0}} (\mathcal{E}_{\alpha_{\Sigma,0}}[E_{\alpha_{\Sigma,0}} \mathbf{Ad}_{\Sigma_{\nearrow 1}}])^{\iota_{\Sigma,0}}.$$

**Lemma 4.8.** [10, Proposition 6.19] *Let $\Sigma$ be a good sequence. We have*

$$\mathbf{Ne}_\Sigma = \varphi_{\Sigma,0} + \varepsilon_{\Sigma,0}\, \mathrm{e}^{\psi_{\Sigma,0}} (E_{\alpha_{\Sigma,0}} \mathbf{Ne}_{\Sigma_{\nearrow 1}})^{\iota_{\Sigma,0}}.$$

Let $\Sigma$ be a good sequence and let $n \in \mathbb{N}$ such that $\Sigma + n$ is nested. For $z \in \mathbf{No}$, we define

$$\Xi_\Sigma z := \varphi_{\Sigma,0} + \varepsilon_{\Sigma,0}\, \mathrm{e}^{\psi_{\Sigma,0}} \left( L_{(\alpha_{\Sigma,0})/\omega} n \left( \frac{\Xi_{\Sigma+n}(\varepsilon_{\Sigma,0}\, \iota_{\Sigma,0}\, z) - \varphi_{\Sigma,0}}{\varepsilon_{\Sigma,0}\, \mathrm{e}^{\psi_{\Sigma,0}}} \right)^{\iota_{\Sigma,0}} \right)^{\iota_{\Sigma,0}},$$

so $\Xi_\Sigma$ is a strictly increasing bijection $\mathbf{No} \longrightarrow \mathbf{Ne}_\Sigma$. By [10, Proposition 6.21], this function does not depend on the choice of $n$.

**Lemma 4.9.** [10, Corollary 6.20] *Let $\Sigma$ be a good sequence. We have*

$$\Xi_\Sigma z = \varphi_{\Sigma,0} + \varepsilon_{\Sigma,0}\, \mathrm{e}^{\psi_{\Sigma,0}} (E_{\alpha_{\Sigma,0}} \Xi_{\Sigma_{\nearrow 1}}(\varepsilon_{\Sigma,0}\, \iota_{\Sigma,0}\, z))^{\iota_{\Sigma,0}}.$$

**Lemma 4.10.** [10, Proposition 6.21] *Let $\Sigma$ be a positive good sequence with $\alpha_{\Sigma,0} \in \omega^{\mathbf{On}+1}$ and let $n \in \mathbb{N}$. For all $z \in \mathbf{No}$, we have*

$$\Xi_{\Sigma-n} z = L_{(\alpha_{\Sigma,0})/\omega} n\, \Xi_\Sigma z.$$

### 4.3 Hyperserial complexity

If $(I, <)$ is a well-ordered set and $(\gamma_i)_{i \in I}$ is a family of ordinals, then the ordered sum of the family is the ordinal $\dot{\sum}_{i \in I} \gamma_j$ defined inductively for $i \in I$ by

$$\dot{\sum}_{j \in [\min(I), i)} \gamma_j := \sup \left\{ \dot{\sum}_{j \in [\min(I), k]} \gamma_j : k < i \right\} \quad \text{if } [\min(I), i) \text{ has no maximum, and}$$

$$\dot{\sum}_{j \in [\min(I), i]} \gamma_j := \left( \dot{\sum}_{j \in [\min(I), i)} \gamma_j \right) + \gamma_i,$$

where $+$ is the (non-commutative) ordinal sum.

We now define a lexicographic ordering on the set $\mathcal{P}_{\max}^{<\infty}(a)$ of *finite* maximal paths in a number $a$. For $P, Q \in \mathcal{P}_{\max}^{<\infty}(a)$, we write $P \triangleleft Q$ if there is an $i < |P|$ with $P(j) = Q(j)$, for all $j < i$ and $P(i) \neq Q(i)$ and

$$P(i), Q(i) \in \operatorname{term} \psi_{P,i} \;\wedge\; P(i) \succcurlyeq Q(i), \quad \text{or}$$
$$P(i) \in \operatorname{term} \psi_{P,i} \;\wedge\; Q(i) \in \operatorname{term} u_{P,i}, \quad \text{or}$$
$$P(i), Q(i) \in \operatorname{term} u_{P,i} \;\wedge\; P(i) \succcurlyeq Q(i).$$

The fact that the paths are finite and maximal implies that only those situations can occur. Since $a$ is well-nested, we have $\mathfrak{m}_{P,i} = \min \operatorname{supp} u_{P,i}$ for large enough $i \in \mathbb{N}$, for all infinite paths $P$. It follows that $\triangleleft$ is well-founded. We define an auxiliary function pn as follows. Given $a \in \mathbf{No}$ define $\operatorname{pn}(a) \in \mathbf{On}$ as the order type of $(\mathcal{P}_{\max}^{<\infty}(a), \triangleleft)$. We have

**Lemma 4.11.** *Let $a \in \mathbf{No}$ and $\mathfrak{m} \in \operatorname{supp} a$. Then $\operatorname{pn}(a_{\succ \mathfrak{m}}) < \operatorname{pn}(a)$.*



**Lemma 4.12.** *Let* $\mathfrak{m} \in \mathbf{Mo}^{\neq}$ *and let* $\mathfrak{m} = \mathrm{e}^{\psi}(L_{\beta} E_{\alpha}^{u})^{\iota}$ *be a hyperserial expansion. Then*

$$\mathrm{pn}(a) = \mathrm{pn}(\psi) + \mathrm{pn}(u).$$

*In particular* $\mathrm{pn}(\psi) < \mathrm{pn}(a)$ *and* $\mathrm{pn}(u) \leqslant \mathrm{pn}(a)$.

**Theorem 4.13.** *Let* $\boldsymbol{\nu} \leqslant \mathbf{On}$ *with* $\boldsymbol{\nu} > 0$. *There is a unique function* $\varsigma_{\boldsymbol{\nu}} : \mathbf{No} \longrightarrow \mathbf{On}$ *such that for all* $a \in \mathbf{No}$, *we have*

a) $\varsigma_{\boldsymbol{\nu}}(a) = \omega$ *if* $a \in \mathbf{Mo}_{\boldsymbol{\lambda}}$ *and* $\varsigma_{\boldsymbol{\nu}}(1) = 1$.

b) $\varsigma_{\boldsymbol{\nu}}(a) = \dot{\sum}_{\mathfrak{m} \in \mathrm{supp}\, a} \varsigma_{\boldsymbol{\nu}}(\mathfrak{m}) \, \ell(a_{\mathfrak{m}})$.

c) *If* $\mathfrak{m} \in \mathbf{Mo}$ *has hyperserial expansion* $\mathfrak{m} = \mathrm{e}^{\psi}(L_{\beta} E_{\alpha}^{u})^{\iota}$ *and there is no* $\boldsymbol{\lambda}$*-bounded nested sequence* $\Sigma$ *such that* $\mathfrak{m}$ *is* $\Sigma$*-nested, then we have*

$$\begin{aligned} \varsigma_{\boldsymbol{\nu}}(\mathfrak{m}) &= \varsigma_{\boldsymbol{\nu}}(\psi) + \varsigma_{\boldsymbol{\nu}}(u) + \beta + 1 \quad \text{if } \alpha < \boldsymbol{\lambda}, \text{ and} \\ \varsigma_{\boldsymbol{\nu}}(\mathfrak{m}) &= \varsigma_{\boldsymbol{\nu}}(\psi) + \omega + \beta \quad \text{if } \alpha = \boldsymbol{\lambda}. \end{aligned}$$

d) *If* $\Sigma$ *is a* $\boldsymbol{\lambda}$*-bounded nested sequence and* $a$ *is* $\Sigma$*-nested, then*

$$\varsigma_{\boldsymbol{\nu}}(a) = \dot{\sum}_{i \in \mathbb{N}} (\varsigma_{\boldsymbol{\nu}}(\varphi_{\Sigma, i}) + \varsigma_{\boldsymbol{\nu}}(\psi_{\Sigma, i})).$$

**Proof.** We define $\varsigma_{\boldsymbol{\nu}}(a)$ by induction on $\mathrm{pn}(a)$. Let $a \in \mathbf{No}$ such that $\varsigma_{\boldsymbol{\nu}}(b)$ is defined for all $b$ with $\mathrm{pn}(b) < \mathrm{pn}(a)$. In view of b, we may assume that $a = \mathfrak{m}$ is a non-zero monomial. If $\mathfrak{m} = 1$, then we set $\varsigma_{\boldsymbol{\nu}}(\mathfrak{m}) = 1$. Assume that $\mathfrak{m} \neq 1$ and write

$$\mathfrak{m} = \mathrm{e}^{\psi_0}(L_{\beta_0} E_{\alpha_0}^{u_1})^{\iota_0}$$

as a hyperserial expansion. We have $\mathrm{pn}(\psi_1) < \mathrm{pn}(a)$ so, by our induction hypothesis, the number $\varsigma_{\boldsymbol{\nu}}(\psi_0)$ is defined. Consider the partial function $\mathrm{lst} : \mathfrak{M} \setminus \{1\} \longrightarrow \mathfrak{M}$ defined for non-trivial monomials $\mathfrak{n}$ that have hyperserial expansions $\mathfrak{n} = \mathrm{e}^{\varphi}(L_{\gamma} E_{\rho}^{v})^{\sigma}$ such that $\mathrm{supp}\, v$ has a minimum, and with $\mathrm{lst}(\mathfrak{n}) := \min \mathrm{supp}\, u$.

If $\mathfrak{m}$ is $\Sigma$-nested for some nested and $\boldsymbol{\lambda}$-bounded sequence $\Sigma$, then the sequence $(\mathrm{lst}^{\circ i}(\mathfrak{m}))_{i \in \mathbb{N}}$ is defined. An induction on $i \in \mathbb{N}$ using Lemmas 4.11 and 4.12 gives that

$$\mathrm{pn}(\varphi_{\Sigma, i}) = \mathrm{pn}((a_{;i})_{\succ \mathrm{lst}^{\circ i}(\mathfrak{m})}) < \mathrm{pn}((a_{;i})) \leqslant \mathrm{pn}(a) \qquad \text{and} \qquad \mathrm{pn}(\psi_{\Sigma, i}) < \mathrm{pn}((a_{;i})) \leqslant \mathrm{pn}(a)$$

for all $i \in \mathbb{N}$. So all $\varsigma_{\boldsymbol{\nu}}(\varphi_{\Sigma, i})$'s and $\varsigma_{\boldsymbol{\nu}}(\psi_{\Sigma, i})$'s are defined, and we may set

$$\varsigma_{\boldsymbol{\nu}}(\mathfrak{m}) := \dot{\sum}_{i \in \mathbb{N}} (\varsigma_{\boldsymbol{\nu}}(\varphi_{\Sigma, i}) + \varsigma_{\boldsymbol{\nu}}(\psi_{\Sigma, i})). \tag{4.3}$$

Assume next that there is an $i > 0$ such that $\mathrm{lst}^{\circ j}(\mathfrak{m})$ is defined for all $j \leqslant i$, and that $a_i := \mathrm{lst}^{\circ i}(\mathfrak{m})$ is $\Sigma_i$-nested for some nested and $\boldsymbol{\lambda}$-bounded sequence $\Sigma_i$. We then choose $i$ to be minimal to satisfy this. Writing

$$\mathfrak{m} = \mathrm{e}^{\psi_0}\left(L_{\beta_0} E_{\alpha_0}^{\varphi_1 + r_1 \mathrm{e}^{\psi_1}\left(L_{\beta_1} E_{\alpha_1}^{\cdot^{\cdot^{\varphi_i + r_i \mathfrak{m}_i}}}\right)^{\iota_1}}\right)^{\iota_0}$$

for some numbers $\varphi_1, \ldots, \varphi_i, \psi_2, \ldots, \psi_{i+1}$ and $r_1, \ldots, r_i$, we have $\mathrm{pn}(\varphi_j), \mathrm{pn}(\psi_{j+1}) < \mathrm{pn}(\mathfrak{m})$ for all $j \in \{1, \ldots, i\}$ as above. So all $\varsigma_{\boldsymbol{\nu}}(\varphi_j)$'s and $\varsigma_{\boldsymbol{\nu}}(\psi_{j+1})$'s are defined. We define $\varsigma_{\boldsymbol{\nu}}(\mathfrak{m}_i)$ as in (4.3) and for all $j < i$, we set

$$\varsigma_{\boldsymbol{\nu}}(\mathrm{lst}^{\circ j}(\mathfrak{m})) := \varsigma_{\boldsymbol{\nu}}(\psi_j) + (\varsigma_{\boldsymbol{\nu}}(\mathrm{lst}^{\circ (j+1)}(\mathfrak{m})) \, \ell(r_1)) + \beta_j + 1 \tag{4.4}$$



by induction on $i-j$.

Assume now that there is an $i \in \mathbb{N}$ such that $\mathrm{lst}^{\circ i}(\mathfrak{m})$ has the hyperserial expansion $\mathrm{lst}^{\circ i}(\mathfrak{m}) = \mathrm{e}^\psi \mathfrak{a}^\iota$ where $\mathfrak{a}$ is $L_{<\boldsymbol{\lambda}}$-atomic. Then we set $\varsigma_{\boldsymbol{\nu}}(\mathrm{lst}^{\circ i}(\mathfrak{m})) := \varsigma_{\boldsymbol{\nu}}(\psi) + \omega$ and we define $\varsigma_{\boldsymbol{\nu}}(\mathfrak{m})$ as in (4.4). The case when $\mathrm{lst}^{\circ i}(\mathfrak{m})$ has hyperserial expansion of type II is similar.

Assume lastly that there is an $i \in \mathbb{N}$, which we choose minimal, such that $b := \mathrm{lst}^{\circ i}(\mathfrak{m})$ has no minimal monomial in its support. Then by Lemma 4.11 and the induction hypothesis, we may define

$$\varsigma_{\boldsymbol{\nu}}(b) := \dot{\sum}_{\mathfrak{n} \in \mathrm{supp}\, b} \varsigma_{\boldsymbol{\nu}}(\mathfrak{n})\, \ell(b_{\mathfrak{n}})$$

We then define $\varsigma_{\boldsymbol{\nu}}(\mathfrak{m})$ as in (4.4).

Well-nestedness implies that the only remaining case is when $\mathrm{lst}^{\circ i}(\mathfrak{m}) = 1$ for some $i$, in which case again $\varsigma_{\boldsymbol{\nu}}(\mathrm{lst}^{\circ i}(\mathfrak{m})) := 1$ is directly definable, and we conclude using (4.4). □

In the sequel we write $\varsigma := \varsigma_{\mathbf{On}}$.

**Lemma 4.14.** *Let $\boldsymbol{\mu} < \boldsymbol{\nu} \leqslant \mathbf{On}$ and write $\delta := \omega^{\boldsymbol{\mu}}$. For $a \in \mathbf{No}^{>, \succ}$, we have $\varsigma_{\boldsymbol{\nu}}(\mathfrak{d}_\delta(a)) \leqslant \varsigma_{\boldsymbol{\nu}}(a)$, with equality if and only if $\mathfrak{a} \in \mathbf{Mo}_\delta$.*

**Proof.** We prove the result by induction on $\varsigma_{\boldsymbol{\nu}}(a)$. Let $a \in \mathbf{No}^{>, \succ}$ such that the result is true for all $\boldsymbol{\mu} < \boldsymbol{\nu}$ for $b \in \mathbf{No}^{>, \succ}$ with $\varsigma_{\boldsymbol{\nu}}(b) < \gamma := \varsigma_{\boldsymbol{\nu}}(a)$. The result is immediate by Theorem 4.13(b) if $\boldsymbol{\mu} = 0$. So we may assume that $\boldsymbol{\mu} > 0$ and that $a \in \mathbf{Mo}^{\succ}$. Write $a = \mathrm{e}^\psi (L_\beta E_\alpha^u)^\iota$ as a hyperserial expansion and write $\mathfrak{a} := \mathfrak{d}_\delta(a)$. If $\psi \neq 0$, then we have $\mathfrak{a} = \mathfrak{d}_\delta(\mathrm{e}^\psi)$ where $\varsigma_{\boldsymbol{\nu}}(\mathrm{e}^\psi) < \varsigma(a)$. The induction hypothesis yields in particular $\varsigma_{\boldsymbol{\nu}}(\mathfrak{a}) < \varsigma_{\boldsymbol{\nu}}(a)$. Assume that $\psi = 0$, so $\iota = 1$. Write $\beta = \beta' + \gamma$ where $\beta' \geqslant_o \delta_{/\omega}$ and $\gamma < \delta_{/\omega}$. We have $\mathfrak{a} = \mathfrak{d}_\delta(L_{\beta'} E_\alpha^u)$ by Lemma 2.3. So $\varsigma_{\boldsymbol{\nu}}(a) \geqslant \varsigma_{\boldsymbol{\nu}}(\mathfrak{a}) + \gamma$. Since $a = \mathfrak{a}$ if and only if $\gamma = 0$ and $\mathfrak{d}_\delta(L_{\beta'} E_\alpha^u) = L_{\beta'} E_\alpha^u$ this yields the result. □

## 5 Extensions by nested numbers

Throughout this section, we consider a subfield $\mathbb{T} = \mathbb{R}[[\mathfrak{M}]]$ of $\mathbf{No}$ of force $(\boldsymbol{\nu}, \boldsymbol{\nu})$. We say that a coding sequence $\Sigma$ *lies over* $\mathbb{T}$, or that $\Sigma$ is a coding sequence over $\mathbb{T}$, if $\varphi_{\Sigma, i}, \psi_{\Sigma, i} \in \mathbb{T}$ for all $i \in \mathbb{N}$. If $\Sigma$ lies over $\mathbb{T}$ and $a \in \mathbf{No}$ is $\Sigma$-nested, then we say that $a$ is *nested over* $\mathbb{T}$. We will see how to adjoin $\mathbb{T}$ with classes of nested numbers over $\mathbb{T}$.

### 5.1 Nested extensions

**Remark 5.1.** Note that if $\Sigma$ is a $\boldsymbol{\lambda}$-bounded good sequence, then the sets $L_\Sigma$ and $R_\Sigma$ of Section 4.2 lie in any subfield of $\mathbf{No}$ of force $(\boldsymbol{\nu}, \boldsymbol{\nu})$ which contains $\varphi_{\Sigma, i}, \psi_{\Sigma, i}$ for all $i \in \mathbb{N}$. So if $\Sigma$ lies over $\mathbb{T}$, then $L_\Sigma$ and $R_\Sigma$ are subsets of $\mathbb{T}$.

**Lemma 5.2.** *Let $\Sigma$ be a $\boldsymbol{\lambda}$-bounded good sequence over $\mathbb{T}$. If $\mathbf{Ad}_\Sigma \cap \mathbb{T} \neq \varnothing$, then $\mathbf{Ne}_\Sigma \cap \mathbb{T} \neq \varnothing$.*

**Proof.** Let $s \in \mathbf{Ad}_\Sigma \cap \mathbb{T}$. Since $\mathbb{T}$ has force $(\boldsymbol{\nu}, \boldsymbol{\nu})$ and $\alpha_{\Sigma, 0} < \boldsymbol{\lambda}$, we may assume that $\Sigma$ is nested. By [10, Lemma 6.16 and Theorem 6.17], there are $i \in \mathbb{N}$, $a \in \mathbf{Ne}_\Sigma$ and $\delta$, $\delta' \prec (E_{\alpha_{\Sigma, i}}^{a; i+1})^{-1}$ with

$$s_{g, i} = \varphi_{\Sigma, i} + \varepsilon_{\Sigma, i} \mathrm{e}^{\psi_{\Sigma, i}} (E_{\alpha_{\Sigma, i}}(a_{;i+1} + \delta))^{\iota_{\Sigma, i}} < s_{;i} < \varphi_{\Sigma, i} + \varepsilon_{\Sigma, i} \mathrm{e}^{\psi_{\Sigma, i}} (E_{\alpha_{\Sigma, i}}(a_{;i+1} + \delta'))^{\iota_{\Sigma, i}} = s_{h, i}.$$



We have $E_{\alpha_{\Sigma,i}} a_{;i+1} \trianglelefteq E_{\alpha_{\Sigma,i}}(a_{;i+1}+\delta), E_{\alpha_{\Sigma,i}}(a_{;i+1}+\delta')$ by [10, (3.7)]. We deduce that

$$a_{;i} = \varphi_{\Sigma,i} + \varepsilon_{\Sigma,i}\, \mathrm{e}^{\psi_{\Sigma,i}} (E_{\alpha_{\Sigma,i}} a_{;i+1})^{\iota_{\Sigma,i}} \trianglelefteq s_{g,i}, s_{h,i},$$

where $a_{;i} \trianglelefteq s_{;i}$. Thus $a_{;i} \in \mathbb{T}$. Since $\Sigma$ is $\boldsymbol{\lambda}$-bounded over $\mathbb{T}$, we deduce that $a \in \mathbb{T}_{(<\boldsymbol{\nu})}$. Since $\mathbb{T}$ has force $(\boldsymbol{\nu},\boldsymbol{\nu})$, we deduce that $a \in \mathbb{T}$. $\square$

**Lemma 5.3.** *Let $\Sigma, \Sigma'$ be distinct $\boldsymbol{\lambda}$-bounded good sequences over $\mathbb{T}$ with $\mathbb{T} \cap \mathbf{Ne}_\Sigma = \mathbb{T} \cap \mathbf{Ne}_{\Sigma'} = \varnothing$. One of the following cases occurs:*

a) $L_{\Sigma'} \not\succ R_\Sigma$ *and then* $\mathbf{Ad}_{\Sigma'} < \mathbf{Ad}_\Sigma$.

b) $L_\Sigma \not\succ R_{\Sigma'}$ *and then* $\mathbf{Ad}_\Sigma < \mathbf{Ad}_{\Sigma'}$.

**Proof.** We have $\mathbf{Ad}_\Sigma = (L_\Sigma \mid R_\Sigma)$ and $\mathbf{Ad}_{\Sigma'} = (L_{\Sigma'} \mid R_{\Sigma'})$, so by Lemma 2.1, it is enough to prove that $\mathbf{Ad}_\Sigma \cap \mathbf{Ad}_{\Sigma'} = \varnothing$. So assume for contradiction that there is a number $a$ with $a \in \mathbf{Ad}_\Sigma \cap \mathbf{Ad}_{\Sigma'}$. Let $i \in \mathbb{N}$ be minimal with $\Sigma(i) \neq \Sigma'(i)$. Considering $\Sigma_{\nearrow i}$ and $\Sigma'_{\nearrow i}$, we may assume that $i = 0$.

Assume for contradiction that $\varphi_{\Sigma,0} \neq \varphi_{\Sigma',0}$. We have $\varphi_{\Sigma,0}, \varphi_{\Sigma',0} \triangleleft a$ so we must have $\varphi_{\Sigma,0} \triangleleft \varphi_{\Sigma',0}$ or $\varphi_{\Sigma',0} \triangleleft \varphi_{\Sigma,0}$ and we may assume without loss of generality that the first case occurs. Write $\varphi_{\Sigma',0} = \varphi_{\Sigma,0} + \varepsilon_{\Sigma,0} t$ where $t \prec \operatorname{supp} \varphi_{\Sigma,0}$. So we have $\mathfrak{d}_t = \mathfrak{d}_{a-\varphi_{\Sigma,0}}$, whence $\varphi_{\Sigma,0} + \varepsilon_{\Sigma,0} \mathfrak{d}_t \in \mathbf{Ad}_\Sigma$. Since $\varphi_{\Sigma,0} + \varepsilon_{\Sigma,0} \mathfrak{d}_t \in \mathbb{T}$, we have $\mathbb{T} \cap \mathbf{Ad}_\Sigma \neq \varnothing$. We deduce with Lemma 5.2 that $\mathbb{T} \cap \mathbf{Ne}_\Sigma \neq \varnothing$: a contradiction.

So $\varphi_{\Sigma,0} = \varphi_{\Sigma',0}$. We have $a > \varphi_{\Sigma,0} \iff \varepsilon_{\Sigma,0} = \varepsilon_{\Sigma',0} = 1$ and $a \leqslant \varphi_{\Sigma,0} \iff \varepsilon_{\Sigma,0} = \varepsilon_{\Sigma',0} = -1$ so we must have $\varepsilon_{\Sigma,0} = \varepsilon_{\Sigma',0}$. The same arguments as for $\varphi_{\Sigma,0}, \varphi_{\Sigma',0}$ imply that $\psi_{\Sigma,0} = \psi_{\Sigma',0}$, and the same argument as for $\varepsilon_{\Sigma,0}, \varepsilon_{\Sigma',0}$ imply that $\iota_{\Sigma,0} = \iota_{\Sigma',0}$. So we must have $\alpha_{\Sigma,0} \neq \alpha_{\Sigma',0}$. We may assume without loss of generality that $\alpha_{\Sigma,0} > \alpha_{\Sigma',0}$. Write $b := E_{\alpha_{\Sigma,0}} a_{\Sigma,1} = E_{\alpha_{\Sigma',0}} a_{\Sigma',1}$. We have

$$c := \varphi_{\Sigma,0} + \varepsilon_{\Sigma,0}\, \mathrm{e}^{\psi_{\Sigma,0}} (\mathfrak{d}_{\alpha_{\Sigma,0}}(b))^{\iota_{\Sigma,0}} \in \mathbf{Ad}_{\Sigma'}$$

by Lemma 4.7. Let $j \in \mathbb{N}^>$ be minimal with $(\varphi_{\Sigma',j}, \psi_{\Sigma',j}) \neq (0,0)$. Let $k \leqslant j$ be maximal with $\alpha_{\Sigma',k} = \alpha_{\Sigma',0}$. We have $\alpha_{\Sigma',0} \geqslant \alpha_{\Sigma',1} \geqslant \cdots \geqslant \alpha_{\Sigma',j-1}$ and $b = E_{\alpha_{\Sigma',0}k + \alpha_{\Sigma',k+1} + \cdots + \alpha_{\Sigma',j-1}} a_{\Sigma',j}$. Note that $\mathfrak{d}_{\alpha_{\Sigma,0}}(a_{\Sigma',k+1}) = \mathfrak{d}_{\alpha_{\Sigma,0}}(\varphi_{\Sigma',j} + \varepsilon_{\Sigma',j}\,\mathrm{e}^{\psi_{\Sigma',j}}) \in \mathbb{T}$ because $\Sigma$ is $\boldsymbol{\lambda}$-bounded. If $\alpha_{\Sigma,0} = \alpha_{\Sigma',0}\omega$, then $\mathfrak{d}_{\alpha_{\Sigma,0}}(b) = E_{\alpha_{\Sigma',0}k} \mathfrak{d}_{\alpha_{\Sigma,0}}(a_{\Sigma',k+1})$, which lies in $\mathbb{T}$ because $\mathbb{T}$ has force $(\boldsymbol{\nu},\boldsymbol{\nu})$ and $\Sigma'$ is $\boldsymbol{\lambda}$-bounded. If $\alpha_{\Sigma,0} > \alpha_{\Sigma',0}\omega$, then [10, Lemma 5.5(b)] gives $\mathfrak{d}_{\alpha_{\Sigma,0}}(b) = \mathfrak{d}_{\alpha_{\Sigma,0}}(\varphi_{\Sigma',j} + \varepsilon_{\Sigma',j}\,\mathrm{e}^{\psi_{\Sigma',j}}) \in \mathbb{T}$. So in any case, we have $c \in \mathbf{Ad}_{\Sigma'} \cap \mathbb{T}$. We deduce with Lemma 5.2 that $\mathbb{T} \cap \mathbf{Ne}_{\Sigma'} \neq \varnothing$: a contradiction. $\square$

Let $\mathbf{P}$ be a class of positive, $\boldsymbol{\lambda}$-bounded good sequences $\Sigma$ over $\mathbb{T}$. Assume that for all $\Sigma \in \mathbf{P}$, we have

$$\forall i \in \mathbb{N}, (\Sigma_{\nearrow i})^+ \;\in\; \mathbf{P} \quad \text{and} \tag{5.1}$$
$$\mathbf{Ne}_\Sigma \cap \mathbb{T} \;=\; \varnothing \quad \text{and} \tag{5.2}$$
$$\Sigma + n \;\in\; \mathbf{P} \quad \text{whenever } \alpha_{\Sigma,0} \in \omega^{\mathbf{On}+1} \text{ and } \Sigma + n \text{ is good.} \tag{5.3}$$

Note that those conditions are satisfied for the class $\mathbf{P}$ of all positive good sequences over $\mathbb{T}$ with $\mathbf{Ne}_{\Sigma_{\nearrow i}} \cap \mathbb{T} = \varnothing$ for all $i \in \mathbb{N}$.

We write $\mathbf{P}_0$ for those sequences in $\mathbf{P}$ which are nested, and set

$$\mathfrak{S}_{\mathbf{P}} := \bigcup_{\Sigma \in \mathbf{P}} \mathbf{Ne}_\Sigma \subseteq \mathbf{Mo}^{\succ}$$



For $\mathfrak{g} \in \mathfrak{G}_{\mathbf{P}}$, we write $n_{\mathfrak{g}}$ for the least number $n_{\mathfrak{g}} \in \mathbb{N}$ such that $\Sigma_{\mathfrak{g}} + n_{\mathfrak{g}}$ is nested. We write $\mathfrak{m}_{\mathfrak{g}}$ for the unique $(\Sigma_{\mathfrak{g}} + n_{\mathfrak{g}})$-nested number with $\mathfrak{g} = L_{(\alpha_{\Sigma_{\mathfrak{g}}},0)/\omega} n(\mathfrak{m}_{\mathfrak{g}})$. We have $\Sigma_{\mathfrak{g}} + n_{\mathfrak{g}} \in \mathbf{P}_0$ by (5.3) so $\mathfrak{m}_{\mathfrak{g}} \in \mathfrak{G}_{\mathbf{P}}$. We write $\mathfrak{G}_{\mathbf{P},0}$ for the class of monomials $\mathfrak{m}_{\mathfrak{g}}$, $\mathfrak{g} \in \mathfrak{G}_{\mathbf{P}}$. For $\mathfrak{m} \in \mathfrak{G}_{\mathbf{P},0}$, we write $\mathfrak{G}_{\mathfrak{m}}$ for the class of monomials $\mathfrak{h}$ with $\mathfrak{m}_{\mathfrak{h}} = \mathfrak{m}$.

**Lemma 5.4.** *For $\mathfrak{g}, \mathfrak{h} \in \mathfrak{G}_{\mathbf{P}}$, we have $\mathfrak{g} \prec \mathfrak{h}$ if and only if one of the following cases occurs:*

a) $\Sigma_{\mathfrak{g}} = \Sigma_{\mathfrak{h}}$ and $z_{\mathfrak{g}} < z_{\mathfrak{h}}$.

b) $L_{\Sigma_{\mathfrak{g}}} \not\succ R_{\Sigma_{\mathfrak{h}}}$ and $\mathbf{Ne}_{\Sigma_{\mathfrak{g}}} < \mathbf{Ne}_{\Sigma_{\mathfrak{h}}}$.

*Moreover, we have $\mathcal{E}_{\alpha_{\mathfrak{g}}}[\mathfrak{g}] \prec \mathcal{E}_{\alpha_{\mathfrak{h}}}[\mathfrak{h}]$ in each case.*

**Proof.** If $\Sigma_{\mathfrak{g}} = \Sigma_{\mathfrak{h}}$, then we have $\mathfrak{g}, \mathfrak{h} \in \mathbf{Ne}_{\Sigma_{\mathfrak{g}}}$ whence $\mathfrak{g} = \Xi_{\Sigma_{\mathfrak{g}}} z_{\mathfrak{g}}$ and $\mathfrak{h} = \Xi_{\Sigma_{\mathfrak{g}}} z_{\mathfrak{h}}$. So we have $\mathfrak{g} < \mathfrak{h}$ if and only if $z_{\mathfrak{g}} < z_{\mathfrak{h}}$ in that case. We obtain $\mathcal{E}_{\alpha_{\mathfrak{g}}}[\mathfrak{g}] \prec \mathcal{E}_{\alpha_{\mathfrak{h}}}[\mathfrak{h}]$ in that case because $\mathfrak{g}$, $\mathfrak{h}$ are both $L_{<\alpha_{\mathfrak{g}}}$-atomic.

Assume now that $\Sigma_{\mathfrak{g}} \neq \Sigma_{\mathfrak{h}}$. By Lemma 5.3, we have $\mathbf{Ne}_{\Sigma_{\mathfrak{g}}} < \mathbf{Ne}_{\Sigma_{\mathfrak{h}}}$ or $\mathbf{Ne}_{\Sigma_{\mathfrak{g}}} > \mathbf{Ne}_{\Sigma_{\mathfrak{h}}}$. So we have $\mathfrak{g} \prec \mathfrak{h}$ if and only if $\mathbf{Ne}_{\Sigma_{\mathfrak{g}}} < \mathbf{Ne}_{\Sigma_{\mathfrak{h}}}$, if and only if $L_{\Sigma_{\mathfrak{g}}} \not\succ R_{\Sigma_{\mathfrak{h}}}$. We have $\mathcal{E}_{\alpha_{\mathfrak{g}}}[\mathbf{Ne}_{\Sigma_{\mathfrak{g}}}] = \mathbf{Ne}_{\Sigma_{\mathfrak{g}}} \prec \mathbf{Ne}_{\Sigma_{\mathfrak{h}}} = \mathcal{E}_{\alpha_{\mathfrak{h}}}[\mathbf{Ne}_{\Sigma_{\mathfrak{h}}}]$ by Lemma 4.7, whence $\mathcal{E}_{\alpha_{\mathfrak{g}}}[\mathfrak{g}] \prec \mathcal{E}_{\alpha_{\mathfrak{h}}}[\mathfrak{h}]$ in that case. □

**Corollary 5.5.** *For $\mathfrak{g}, \mathfrak{h} \in \mathfrak{G}_{\mathbf{P}}$ with $\mathfrak{g} \prec \mathfrak{h}$, we have $\mathfrak{L}_{<(\alpha_{\mathfrak{g}})/\omega} \circ \mathfrak{g} \prec \mathfrak{L}^{\succ}_{<(\alpha_{\mathfrak{h}})/\omega} \circ \mathfrak{h}$.*

For every family $\mathfrak{f} := (\mathfrak{f}_{\mathfrak{g}})_{\mathfrak{g} \in \mathfrak{G}_{\mathbf{P}}} \in \prod_{\mathfrak{g} \in \mathfrak{G}_{\mathbf{P}}} \mathfrak{L}_{<(\alpha_{\mathfrak{g}})/\omega}$, we write $S(\mathfrak{f}) := \{\mathfrak{g} \in \mathfrak{G}_{\mathbf{P}} : \mathfrak{f}_{\mathfrak{g}} \neq 1\}$. We let $\mathcal{F}_{\mathbf{P}}$ denote the set of families $\mathfrak{f} \in \prod_{\mathfrak{g} \in \mathfrak{G}_{\mathbf{P}}} \mathfrak{L}_{<(\alpha_{\mathfrak{g}})/\omega}$ such that the set $S_0(\mathfrak{f}) := \{\mathfrak{m}_{\mathfrak{g}} : \mathfrak{g} \in S(\mathfrak{f})\}$ is finite. For $\mathfrak{f} \in \mathcal{F}_{\mathbf{P}}$, and each $\mathfrak{m} \in S_0(\mathfrak{f})$, the sum $\sum_{\mathfrak{m}_{\mathfrak{g}} = \mathfrak{m}} (\log(\mathfrak{f}_{\mathfrak{g}}) \circ \mathfrak{g})$ is well-defined by Corollary 5.5. So $\sum_{\mathfrak{g} \in \mathfrak{G}_{\mathbf{P}}} (\log(\mathfrak{f}_{\mathfrak{g}}) \circ \mathfrak{g})$ is well-defined and lies in $\mathbf{No}_{\succ}$. We set $\tilde{\mathfrak{f}} := e^{\sum_{\mathfrak{g} \in \mathfrak{G}_{\mathbf{P}}} \log(\mathfrak{f}_{\mathfrak{g}}) \circ \mathfrak{g}} \in \mathbf{Mo}$. Note that $\mathfrak{f}^{-1} := (\mathfrak{f}_{\mathfrak{g}}^{-1})_{\mathfrak{g} \in \mathfrak{G}_{\mathbf{P}}}$ lies in $\mathcal{F}_{\mathbf{P}}$ and we have $\widetilde{\mathfrak{f}^{-1}} = \tilde{\mathfrak{f}}^{-1}$.

Therefore the class $\mathfrak{M}_{[\mathbf{P}]}$ of monomials $\tilde{\mathfrak{f}}$ for $\mathfrak{f} \in \mathcal{F}_{\mathbf{P}}$ is a subgroup of $\mathbf{Mo}$. We see that $\mathfrak{M}_{[\mathbf{P}]}$ is the subgroup of $\mathbf{Mo}$ generated by the sets $\mathfrak{L}_{<\alpha_{\Sigma,0}} \circ \mathfrak{G}_{\mathfrak{m}}$ where $\Sigma \in \mathbf{P}_0$ and $\mathfrak{m} \in \mathbf{Ne}_{\Sigma}$. By the previous argument, the ordering on $\mathfrak{M}_{[\mathbf{P}]}$ corresponds to the lexicographic ordering on $\mathcal{F}_{\mathbf{P}}$. We write $\mathfrak{M}_{\mathbf{P}}$ for the subgroup of $\mathbf{Mo}$ generated by $\mathfrak{M}$ and $\mathfrak{M}_{[\mathbf{P}]}$.

We claim that $\mathfrak{M} \cap \mathfrak{M}_{[\mathbf{P}]} = \{1\}$. Indeed assume for contradiction that there is a monomial $\mathfrak{n} \succ 1$ in $\mathfrak{M}_{[\mathbf{P}]} \cap \mathfrak{M}$. So there is $\mathfrak{f} \in \mathcal{F}_{\mathbf{P}}$ with $\mathfrak{n} = \tilde{\mathfrak{f}}$. In particular, $\mathfrak{d}_{\log \mathfrak{n}} = \mathfrak{d}_{\log \tilde{\mathfrak{f}}}$. So there are $\mathfrak{g} \in \mathfrak{G}_{\mathbf{P}}$ and $\gamma < (\alpha_{\mathfrak{g}})/\omega$ with $\ell_{\gamma+1} \circ \mathfrak{g} \in \mathfrak{M}$. But there are $\varphi \in \mathbb{T}, \psi \in \mathbb{T}_{\succ}$ and $\varepsilon, \iota \in \{-1, 1\}$ with $\varphi + \varepsilon e^{\psi} \mathfrak{g}^{\iota} \in \mathbf{Ne}_{\Sigma}$. So by Lemma 4.7, we have $\varphi + \varepsilon e^{\psi} \mathcal{E}_{\alpha_{\mathfrak{g}}}[\mathfrak{g}]^{\iota} \in \mathbf{Ne}_{\Sigma}$, whence $\varphi + \varepsilon e^{\psi} \mathcal{E}_{\alpha_{\mathfrak{g}}}[\mathfrak{g}]^{\iota} \cap \mathbb{T} = \varnothing$. We deduce that $\log \mathcal{E}_{\alpha_{\mathfrak{g}}}[\mathfrak{g}] \cap \mathbb{T} = \varnothing$, whence in particular $\ell_{\gamma+1} \circ \mathfrak{g} \notin \mathbb{T}$: a contradiction.

So any element of $\mathfrak{M}_{\mathbf{P}}$ can be uniquely written as $\mathfrak{n} = \mathfrak{m} \tilde{\mathfrak{f}}$ where $\mathfrak{m} \in \mathfrak{M}$ and $\mathfrak{f} \in \mathcal{F}_{\mathbf{P}}$. We write $\mathbb{T}_{\mathbf{P}}$ for the field $\mathbb{T}_{\mathbf{P}} = \mathbb{R}[[\mathfrak{M}_{\mathbf{P}}]]$.

**Lemma 5.6.** *Let $\mathfrak{m} \in \mathfrak{M}$ and $\mathfrak{f} \in \mathcal{F}_{\mathbf{P}}^{\succ}$. Let $\mathfrak{g} := \max S(\mathfrak{f})$. We have $\mathfrak{m} \tilde{\mathfrak{f}} \succ 1$ if and only if $\mathbf{Ne}_{\Sigma_{\mathfrak{g}}} \succ \mathfrak{m}^{-1}$.*

**Proof.** We have $|\log \tilde{\mathfrak{f}}| \in \mathcal{E}_{\alpha_{\mathfrak{g}}}[|\log \mathfrak{m}_{\mathfrak{g}}|]$ by the previous arguments. So $\tilde{\mathfrak{f}} \in \mathcal{E}_{\alpha_{\mathfrak{g}}}[\mathfrak{g}] \subseteq \mathcal{E}_{\alpha_{\mathfrak{g}}}[\mathbf{Ad}_{\Sigma}]$. We have $\mathbb{T} \cap \mathcal{E}_{\alpha_{\mathfrak{g}}}[\mathbf{Ad}_{\Sigma}] = \mathbb{T} \cap \mathbf{Ad}_{\Sigma} = \varnothing$ by (5.2) and Lemma 4.7. Since the class $\mathcal{E}_{\alpha_{\mathfrak{g}}}[\mathbf{Ad}_{\Sigma}]$ is convex in, we either have $\mathfrak{m}^{-1} \prec \mathcal{E}_{\alpha_{\mathfrak{g}}}[\mathbf{Ad}_{\Sigma}]$ or $\mathfrak{m}^{-1} \succ \mathcal{E}_{\alpha_{\mathfrak{g}}}[\mathbf{Ad}_{\Sigma}]$. We obtain the result by noticing that $\mathbf{Ne}_{\Sigma}$ is cofinal and coinitial in $\mathbf{Ad}_{\Sigma}$. □

**Lemma 5.7.** *We have $\mathbb{T}_{\mathbf{P}} = \mathcal{H}_{\nu}(\bigcup_{\Sigma \in \mathbf{P}_0} \mathbf{Ne}_{\Sigma} \cup \mathfrak{M})$. Moreover $(\mathfrak{M}_{\mathbf{P}})_{\boldsymbol{\lambda}} = \mathfrak{M}_{\boldsymbol{\lambda}}$.*



**Proof.** Let $\mathbb{U}$ be a subfield of force $\nu$ with $\mathbb{U} \supseteq \bigcup_{\Sigma \in \mathbf{P}_0} \mathbf{Ne}_\Sigma \cup \mathfrak{M}$. For $\Sigma \in \mathbf{P}_0$ and $\mathfrak{m} \in \mathbf{Ne}_\Sigma$, we have $\mathfrak{L}_{<\alpha_{\Sigma,0}} \circ \mathfrak{m} \subseteq \mathbb{U}$, so $\mathfrak{L}_{<\alpha_{\Sigma,0}} \circ \mathfrak{G}_\mathfrak{m} \subseteq \mathbb{U}$. We deduce that $\mathfrak{M}_{[\mathbf{P}]} \subseteq \mathbb{U}$, whence $\mathfrak{M}_\mathbf{P} \subseteq \mathbb{U}$, whence $\mathbb{T}_\mathbf{P} \subseteq \mathbb{U}$.

So we need only justify that $\mathbb{T}_\mathbf{P}$ is a confluent subfield of force $\nu$. Let $\mathfrak{n} \in \mathfrak{M}_\mathbf{P}^\succ$ and write $\mathfrak{n} = \mathfrak{m}\,\tilde{\mathfrak{f}}$ where $\mathfrak{f} \in \mathcal{F}_\mathbf{P}$ and $\mathfrak{m} \in \mathfrak{M}$. Let $\boldsymbol{\mu} \leqslant \boldsymbol{\nu}$ with $\boldsymbol{\mu} > 0$ and write $\mathfrak{a} := \mathfrak{d}_{\omega^\mu}(\mathfrak{n})$. By Proposition 3.3, it is enough to prove that $\mathfrak{a} \in \mathbb{T}_\mathbf{P}$, that $\mathfrak{L}_{<\omega^\mu} \circ \mathfrak{a} \subseteq \mathbb{T}_\mathbf{P}$, and that $L_{\omega^\mu}(\mathfrak{a}) \in \mathbb{T}_\mathbf{P}$ if moreover $\boldsymbol{\mu} < \boldsymbol{\nu}$.

Let $\mathfrak{g} \in \mathfrak{G}_\mathbf{P}$ be minimal with $\mathfrak{f}_\mathfrak{g} \neq 1$ and let $\gamma < (\alpha_\mathfrak{g})_{/\omega}$ be minimal with $(\mathfrak{f}_\mathfrak{g})_\gamma \neq 0$. We have $\log \mathfrak{n} \asymp \log \mathfrak{m}$ or $\log \mathfrak{n} \asymp \ell_{\gamma+1} \circ \mathfrak{g}$. Therefore we have

$$\mathfrak{a} \in \{e^{\mathfrak{d}_\omega(\log \mathfrak{m})}, \mathfrak{d}_{\omega^\mu}(\mathfrak{m})\} = \{\mathfrak{d}_{\omega^\mu}(\mathfrak{m})\} \subseteq \mathfrak{d}_{\omega^\mu}(\mathbb{T})$$

or

$$\mathfrak{a} \in \{e^{\mathfrak{d}_\omega(\ell_{\gamma+1} \circ \mathfrak{g})}, \mathfrak{d}_{\omega^\mu}(\ell_\gamma \circ \mathfrak{g})\} = \{\mathfrak{d}_{\omega^\mu}(\ell_\gamma \circ \mathfrak{g})\}.$$

In the first case, since $\mathbb{T}$ is a confluent subfield of force $\nu$, we obtain $\mathfrak{a} \in \mathbb{T}$ and $\mathfrak{L}_{<\omega^\mu} \circ \mathfrak{a} \subseteq \mathbb{T}$ and $L_{\omega^\mu}(\mathfrak{a}) \in \mathbb{T}$ if moreover $\boldsymbol{\mu} < \boldsymbol{\nu}$. So we may assume that $\mathfrak{n} = \ell_\gamma \circ \mathfrak{g}$. Write $\gamma = \gamma' + \rho$ where $\operatorname{supp} \gamma' \succcurlyeq \omega^{\mu-}$ and $\rho < \omega^{\mu-}$. So $\mathfrak{a} = \mathfrak{d}_{\omega^\mu}(\ell_{\gamma'} \circ \mathfrak{g})$, so we may assume that $\gamma = \gamma'$ and that $\operatorname{supp} \gamma \succcurlyeq \omega^{\mu-}$. Write $\Sigma := \Sigma_\mathfrak{g}$ and $\alpha_i := \alpha_{\Sigma,i}$ for all $i \in \mathbb{N}$. Consider $i > 0$ minimal with $(\varphi_{\Sigma,i}, \psi_{\Sigma,i}) \neq (0,0)$. So $\alpha_0 \geqslant \cdots \geqslant \alpha_{i-1}$ and

$$\mathfrak{g} = E_{\alpha_0 + \cdots + \alpha_{i-1}}(\varphi_{\Sigma,i} + \varepsilon_{\Sigma,i}\, e^{\psi_{\Sigma,i}} \mathfrak{m}^{\iota_{\Sigma,i}}),$$

for a certain $\mathfrak{m} \in \mathbf{Ne}_{(\Sigma_{\nearrow 1})^+}$.

Assume that $\alpha_0 \geqslant \omega^\mu$, so $\boldsymbol{\mu} < \boldsymbol{\nu}$. We have $\mathfrak{a} = \ell_\gamma \circ \mathfrak{g} \in \mathbb{T}_\mathbf{P}$. Moreover $\mathfrak{L}_{<\omega^\mu} \circ \mathfrak{a} \subseteq \mathfrak{L}_{<\alpha_0} \circ \mathfrak{G}_{\mathfrak{m}_\mathfrak{g}} \subseteq \mathbb{T}_\mathbf{P}$. If $\alpha_0 = \omega^\mu$, then we must have $\gamma = 0$ and $L_{\omega^\mu}(\mathfrak{a}) = L_{\alpha_0}(\mathfrak{g}) = \mathfrak{g}_{;1} \in \mathbb{T}_\mathbf{P}$. Otherwise write $\gamma = \gamma'' + \omega^{\mu-} n$ for $\gamma'' >_o \omega^{\mu-}$ and $n \in \mathbb{N}$. So $L_{\omega^\mu}(\mathfrak{a}) = \ell_{\gamma'' + \omega^\mu} \circ \mathfrak{g} - n$ where $\ell_{\gamma'' + \omega^\mu} \circ \mathfrak{g} \in \mathfrak{L}_{<\alpha_0} \circ \mathfrak{G}_{\mathfrak{m}_\mathfrak{g}}$. So $L_{\omega^\mu}(\mathfrak{a}) \in \mathbb{T}_\mathbf{P}$.

Assume now that $\alpha_0 < \omega^\mu$. Let $j \leqslant i-1$ be maximal with either $j = 0$ or $j > 0$ and $\alpha_j \omega = \omega^\mu$. So we have

$$\mathfrak{b} := \mathfrak{d}_{\omega^\mu}(E_{\alpha_{j+1}+\cdots+\alpha_{i-1}}(\varphi_i + \varepsilon_i\, e^{\psi_i} \mathfrak{m}^{\iota_i})) = \mathfrak{d}_{\omega^\mu}(\varphi_{\Sigma,i} + \varepsilon_{\Sigma,i}\, e^{\psi_{\Sigma,i}}) \in \mathbb{T},$$

and $\mathfrak{a} = E_{\alpha_0 j} \mathfrak{b} \in \mathbb{T}$ because $\mathbb{T}$ has force $(\boldsymbol{\nu}, \boldsymbol{\nu})$.

In this case we obtain again $\mathfrak{a} \in \mathbb{T}$, and $\mathfrak{L}_{<\omega^\mu} \circ \mathfrak{a} \subseteq \mathbb{T}$ and $L_{\omega^\mu}(\mathfrak{a}) \in \mathbb{T}$ if $\boldsymbol{\mu} < \boldsymbol{\nu}$. □

## 5.2 Nested extensions of embeddings

**Proposition 5.8.** *Let $\Phi \colon \mathbb{T} \longrightarrow \mathbf{No}$ be an embedding of force $\boldsymbol{\nu}$ and let $\Sigma$ be a $\boldsymbol{\lambda}$-bounded good sequence which lies over $\mathbb{T}$. The sequence*

$$\Phi(\Sigma) := (\Phi(\varphi_{\Sigma,i}), \varepsilon_i, \Phi(\psi_{\Sigma,i}), \iota_i, \alpha_i)_{i \in \mathbb{N}}$$

*is good. Moreover, it is nested if $\Sigma$ is nested.*

**Proof.** Note that we have $\Phi(\varphi_{\Sigma,i+1}) \in \mathbf{No}_{\succ, \alpha_{\Sigma,i}}$ and $\Phi(\psi_{\Sigma,i}) \in \mathbf{No}_\succ$ for all $i \in \mathbb{N}$. The other conditions of coding sequences are clearly preserved, so $\Phi(\Sigma)$ is a coding sequence. We have $\Phi(\Sigma + n) = \Phi(\Sigma) + n$ for all $n \in \mathbb{N}$ so we may assume that $\Sigma$ is nested and prove that $\Phi(\Sigma)$ is nested as well. By Remark 5.1, we have $L_{\Phi(\Sigma)} = \Phi(L_\Sigma) < \Phi(R_\Sigma) = R_{\Phi(\Sigma)}$. So $\Phi(\Sigma)$ is admissible. Finally to prove that $\Phi(\Sigma)$ is nested, we must justify that for $i \in \mathbb{N}$, we have

$$\mathbf{Ad}_{\Phi(\Sigma)_{\nearrow i}} \supseteq \Phi(\varphi_{\Sigma,i+1}) + \varepsilon_{\Sigma,i}\, e^{\Phi(\psi_{\Sigma,i})} (E_{\alpha_{\Sigma,i}} \mathbf{Ad}_{\Phi(\Sigma)_{\nearrow i+1}})^{\iota_{\Sigma,i}}.$$



If $\varepsilon_{\Sigma,i}\, \iota_{\Sigma,i}=1$, then (5.8) is equivalent to

$\Phi(\varphi_{\Sigma,i+1})+\varepsilon_{\Sigma,i}\mathrm{e}^{\Phi(\psi_{\Sigma,i})}(E_{\alpha_{\Sigma,i}}\Phi(L_{\Sigma_{\nearrow i+1}}))^{\iota_{\Sigma,i}}$ is cofinal with respect to $\Phi(L_{\Sigma_{\nearrow i}})$ and

$\Phi(\varphi_{\Sigma,i+1})+\varepsilon_{\Sigma,i}\mathrm{e}^{\Phi(\psi_{\Sigma,i})}(E_{\alpha_{\Sigma,i}}\Phi(R_{\Sigma_{\nearrow i+1}}))^{\iota_{\Sigma,i}}$ is coinitial with respect to $\Phi(R_{\Sigma_{\nearrow i}})$.

Those statements hold because $\Phi:\mathbb{T}\longrightarrow\Phi(\mathbb{T})$ is an order isomorphism. If $\varepsilon_{\Sigma,i}\,\iota_{\Sigma,i}=-1$, then (5.8) is equivalent to

$\Phi(\varphi_{\Sigma,i+1})+\varepsilon_{\Sigma,i}\mathrm{e}^{\Phi(\psi_{\Sigma,i})}(E_{\alpha_{\Sigma,i}}\Phi(L_{\Sigma_{\nearrow i+1}}))^{\iota_{\Sigma,i}}$ is cofinal with respect to $\Phi(R_{\Sigma_{\nearrow i}})$ and

$\Phi(\varphi_{\Sigma,i+1})+\varepsilon_{\Sigma,i}\mathrm{e}^{\Phi(\psi_{\Sigma,i})}(E_{\alpha_{\Sigma,i}}\Phi(R_{\Sigma_{\nearrow i+1}}))^{\iota_{\Sigma,i}}$ is coinitial with respect to $\Phi(L_{\Sigma_{\nearrow i}})$,

which holds as well. So $\Phi(\Sigma)$ is nested. $\square$

**Proposition 5.9.** *Let $\mathbb{T}$ and $\mathbf{P}$ be as above and let $\Phi:\mathbb{T}\longrightarrow\mathbf{No}$ be an embedding of force $\boldsymbol{\nu}$. Consider a family $(\Phi_\Sigma)_{\Sigma\in\mathbf{P}}$ of order isomorphisms $\Phi_\Sigma:\mathbf{No}\longrightarrow\mathbf{No}$ with*

$$\Phi_\Sigma(z) = \varepsilon_{\Sigma,1}\,\iota_{\Sigma,1}\,\Phi_{(\Sigma_{\nearrow 1})^+}(\varepsilon_{\Sigma,1}\,\iota_{\Sigma,1}\,z) \quad and \tag{5.4}$$
$$\Phi_{\Sigma-1}(z) = \Phi_\Sigma(z) \quad whenever\ \Sigma-1\in\mathbf{P}. \tag{5.5}$$

*for $\Sigma\in\mathbf{P}$ and $z\in\mathbf{No}$. There is a unique extension $\Phi_\mathbf{P}$ of $\Phi$ into a hyperserial embedding $\mathbb{T}_\mathbf{P}\longrightarrow\mathbf{No}$ of force $\boldsymbol{\nu}$ with*

$$\Phi_\mathbf{P}(\Xi_\Sigma z) = \Xi_{\Phi(\Sigma)}\Phi_\Sigma(z) \quad for\ all\ \Sigma\in\mathbf{P}\ and\ z\in\mathbf{No}. \tag{5.6}$$

**Proof.** By Proposition 5.8, the class $\Phi(\mathbf{P}):=\{\Phi(\Sigma):\Sigma\in\mathbf{P}\}$ satisfies (5.1,5.2,5.3) with respect to $\Phi(\mathbb{T})$. We define a function $\Psi:\mathfrak{G}_{\Phi(\mathbf{P})}\longrightarrow\mathfrak{G}_\mathbf{P}$. For $\mathfrak{g}\in\mathfrak{G}_{\Phi(\mathbf{P})}$, there is a unique $\Sigma\in\mathbf{P}$ with $\Sigma_\mathfrak{g}=\Phi(\Sigma)$. We set $\Psi(\mathfrak{g}):=\Xi_\Sigma\Phi_\Sigma^{-1}(z_\mathfrak{g})$. We claim that $\Psi$ is an order isomorphism $\mathfrak{G}_{\Phi(\mathbf{P})}\longrightarrow\mathfrak{G}_\mathbf{P}$.

Indeed let $\mathfrak{g},\mathfrak{h}\in\mathfrak{G}_{\Phi(\mathbf{P})}$ with $\mathfrak{g}\prec\mathfrak{h}$. If $\Sigma_\mathfrak{g}=\Sigma_\mathfrak{h}$, then by Lemma 5.4, we have $z_\mathfrak{g}<z_\mathfrak{h}$, so $\Phi_\Sigma^{-1}(z_\mathfrak{g})<\Phi_\Sigma^{-1}(z_\mathfrak{h})$, so $\Psi(\mathfrak{g})\prec\Psi(\mathfrak{h})$. Otherwise, by Lemma 5.4, we have $L_{\Sigma_\mathfrak{g}}\not\succ R_{\Sigma_\mathfrak{h}}$. Write $\Sigma,\Sigma'\in\mathbf{P}$ with $\Sigma_\mathfrak{g}=\Phi(\Sigma)$ and $\Sigma_\mathfrak{h}=\Phi(\Sigma')$. We deduce with Remark 5.1 that $L_\Sigma\not\succ R_{\Sigma'}$, so $\mathbf{Ne}_\Sigma<\mathbf{Ne}_{\Sigma'}$, whence in particular $\Psi(\mathfrak{g})\prec\Psi(\mathfrak{h})$. The function $\Psi$ is surjective because $\Phi:\mathbf{P}\longrightarrow\Phi(\mathbf{P})$ is surjective and each $\Phi_\Sigma:\mathbf{No}\longrightarrow\mathbf{No}$ for $\Sigma\in\mathbf{P}$ is surjective. So this proves the claim.

We set $\Phi_\mathbf{P}(\tilde{\mathfrak{f}}):=\widetilde{\mathfrak{f}\circ\Psi}$ for all $\mathfrak{f}\in\mathcal{F}_\mathbf{P}$. This defines an isomorphism of ordered groups $\mathfrak{M}_{[\mathbf{P}]}\longrightarrow(\Phi(\mathfrak{M}))_{[\Phi(\mathbf{P})]}$. We then set $\Phi_\mathbf{P}(\mathfrak{m}\,\tilde{\mathfrak{f}}):=\Phi(\mathfrak{m})\,\Phi_\mathbf{P}(\tilde{\mathfrak{f}})$ for all $\mathfrak{m}\in\mathfrak{M}$ and $\mathfrak{f}\in\mathcal{F}_\mathbf{P}$. Consider $\mathfrak{m}\in\mathfrak{M}$ and $\mathfrak{f}\in\mathcal{F}_\mathbf{P}$ with $\mathfrak{m}\,\tilde{\mathfrak{f}}\succ 1$. If $\tilde{\mathfrak{f}}=1$ then $\mathfrak{m}\succ 1$ so $\Phi_\mathbf{P}(\mathfrak{m}\,\tilde{\mathfrak{f}})=\Phi(\mathfrak{m})\succ 1$. Assume that $\tilde{\mathfrak{f}}\succ 1$. By Lemma 5.6, we have $\mathfrak{m}^{-1}\prec\mathbf{Ne}_\Sigma$ where $\mathfrak{g}=\max S(\mathfrak{f})$. So there is $b\in L_\Sigma$ with $\mathfrak{m}^{-1}<b$. We have $\Psi^{-1}(\mathfrak{g})=\max S(\mathfrak{f}\circ\Psi)$. Note that $\Phi(\mathfrak{m}^{-1})\prec\Phi(b)$ where $\Phi(b)\in\Phi(L_\Sigma)=L_{\Phi(\Sigma)}$. We deduce that $\mathfrak{m}^{-1}\prec\mathbf{Ne}_{\Phi(\Sigma)}$, so $\Phi_\mathbf{P}(\mathfrak{m}\,\tilde{\mathfrak{f}})\succ 1$. If $\tilde{\mathfrak{f}}\prec 1$, then we have $\Phi_\mathbf{P}(\mathfrak{m}\,\tilde{\mathfrak{f}})=\Phi_\mathbf{P}((\mathfrak{m}\,\tilde{\mathfrak{f}})^{-1})^{-1}$ where $\Phi_\mathbf{P}((\mathfrak{m}\,\tilde{\mathfrak{f}})^{-1})\prec 1$ by the previous arguments. We deduce that $\Phi_\mathbf{P}(\mathfrak{m}\,\tilde{\mathfrak{f}})\succ 1$. So in general $\Phi_\mathbf{P}(\mathfrak{M}_\mathbf{P}^\succ)\subseteq\mathbf{Mo}^\succ$, which implies that $\Phi_\mathbf{P}$ is an embedding of ordered groups $\mathfrak{M}_\mathbf{P}\longrightarrow\mathbf{Mo}$.

So $\Phi_\mathbf{P}$ further extends uniquely as a strongly linear embedding of ordered fields $\mathbb{T}_\mathbf{P}\longrightarrow\mathbf{No}$ which extends $\Phi$. Note that by definition, this embedding commutes with the logarithm. We claim that $\Phi_\mathbf{P}$ is a hyperserial embedding of force $\boldsymbol{\nu}$.

By Lemma 3.7, we need only prove that $L_{\omega^\eta}(\Phi_\mathbf{P}(\mathfrak{a}))=\Phi_\mathbf{P}(L_{\omega^\eta}(\mathfrak{a}))$ for all $\eta<\boldsymbol{\nu}$ with $\eta>0$ and all $\mathfrak{a}\in(\mathfrak{M}_\mathbf{P})_{<\omega^\eta}$. Consider such elements $\eta,\mathfrak{a}$. We may assume that $\mathfrak{a}\notin\mathfrak{M}$, so $\mathfrak{a}=\ell_\gamma\circ\mathfrak{g}$ for a certain $\mathfrak{g}\in\mathfrak{G}_\mathbf{P}$ and $\gamma<(\alpha_\mathfrak{g})_{/\omega}$. We must have $\alpha_\mathfrak{g}\geqslant\omega^\eta$ and $\gamma\gg\omega^{\eta-}$ because $\mathfrak{a}$ is $L_{<\omega^\eta}$-atomic. Write $\Sigma=\Sigma_\mathfrak{g}$ and rite $\gamma=\gamma'+\omega^{\eta-}n$ where $\gamma'\gg\omega^{\eta-}$ and $n\in\mathbb{N}$.



If $\alpha_{\mathfrak{g}} = \omega^\eta$, then $\gamma = 0$. So $\mathfrak{a} = \mathfrak{g}$ and $\Phi_{\mathbf{P}}(\mathfrak{a}) = \Psi^{-1}(\mathfrak{g}) = \Xi_{\Phi(\Sigma)} \Phi_\Sigma(z_\mathfrak{g})$. By Lemma 4.9, we have

$$L_{\omega^\eta}(\mathfrak{a}) = \varphi_{\Sigma,1} + \varepsilon_{\Sigma,1} \, e^{\psi_{\Sigma,1}} \, (\Xi_{(\Sigma_{\nearrow 1})^+} (\varepsilon_{\Sigma,1} \, \iota_{\Sigma,1} \, z_\mathfrak{g}))^{\iota_{\Sigma,1}}.$$

Lemma 4.9 also yields

$$\begin{aligned}
L_{\omega^\eta}(\Phi_{\mathbf{P}}(\mathfrak{a})) &= \varphi_{\Phi(\Sigma),1} + \varepsilon_{\Sigma,1} \, e^{\psi_{\Phi(\Sigma),1}} \, (\Xi_{(\Phi(\Sigma)_{\nearrow 1})^+} (\varepsilon_{\Sigma,1} \, \iota_{\Sigma,1} \, \Phi_\Sigma(z_\mathfrak{g})))^{\iota_{\Sigma,1}} \\
&= \Phi(\varphi_{\Sigma,1}) + \varepsilon_{\Sigma,1} \, \Phi(e^{\psi_{\Sigma,1}}) \, (\Xi_{\Phi((\Sigma_{\nearrow 1})^+)} (\varepsilon_{\Sigma,1} \, \iota_{\Sigma,1} \, \Phi_\Sigma(z_\mathfrak{g})))^{\iota_{\Sigma,1}} \\
&= \Phi_{\mathbf{P}}(\varphi_{\Sigma,1}) + \varepsilon_{\Sigma,1} \, \Phi_{\mathbf{P}}(e^{\psi_{\Sigma,1}}) \, (\Xi_{\Phi((\Sigma_{\nearrow 1})^+)} \Phi_{(\Sigma_{\nearrow 1})^+} (\varepsilon_{\Sigma,1} \, \iota_{\Sigma,1} \, z_\mathfrak{g}))^{\iota_{\Sigma,1}} \quad \text{(by (5.4))} \\
&= \Phi_{\mathbf{P}}(\varphi_{\Sigma,1} + \varepsilon_{\Sigma,1} \, e^{\psi_{\Sigma,1}} \, (\Xi_{(\Sigma_{\nearrow 1})^+} (\varepsilon_{\Sigma,1} \, \iota_{\Sigma,1} \, z_\mathfrak{g}))^{\iota_{\Sigma,1}}) \\
&= \Phi_{\mathbf{P}}(L_{\omega^\eta}(\mathfrak{a})).
\end{aligned}$$

If $\alpha_{\mathfrak{g}} = \omega^{\eta+1}$, then $\gamma' = 0$ so $\mathfrak{a} = \ell_{\omega^{\eta-n}} \circ \mathfrak{g}$ and $\Phi_{\mathbf{P}}(\mathfrak{a}) = \ell_{\omega^{\eta-n}} \circ \Xi_{\Phi(\Sigma)} \Phi_\Sigma(z_\mathfrak{g})$ and $L_{\omega^\eta}(\mathfrak{a}) = L_{\omega^\eta}(\mathfrak{g}) - n$. We have $\Sigma_{L_{\omega^\eta}(\mathfrak{g})} = \Sigma - 1$, and $z_{L_{\omega^\eta}(\mathfrak{g})} = z_\mathfrak{g}$ by Lemma 4.10. So

$$\begin{aligned}
\Phi_{\mathbf{P}}(L_{\omega^\eta}(\mathfrak{a})) &= \Xi_{\Phi(\Sigma-1)} \Phi_{\Sigma-1}(z_{L_{\omega^\eta}(\mathfrak{g})}) - n \\
&= \Xi_{\Phi(\Sigma-1)} \Phi_\Sigma(z_\mathfrak{g}) - n & \text{(by (5.5))} \\
&= \Xi_{\Phi(\Sigma)-1} \Phi_\Sigma(z_\mathfrak{g}) - n \\
&= L_{\omega^\eta}(\Xi_{\Phi(\Sigma)} \Phi_\Sigma(z_\mathfrak{g})) - n & \text{(by Lemma 4.10)} \\
&= L_{\omega^\eta}(\Phi_{\mathbf{P}}(\mathfrak{a})).
\end{aligned}$$

Finally, if $\alpha_{\mathfrak{g}} > \omega^{\eta+1}$, then $\Phi_{\mathbf{P}}(\mathfrak{a}) = \ell_\gamma \circ \Psi^{-1}(\mathfrak{g})$ and $L_{\omega^\eta}(\mathfrak{a}) = \ell_{\gamma'+\omega^\eta} \circ \mathfrak{g} - n$ so

$$\Phi_{\mathbf{P}}(L_{\omega^\eta}(\mathfrak{a})) = \ell_{\gamma'+\omega^\eta} \circ \Psi^{-1}(\mathfrak{g}) - n = L_{\omega^\eta}(\Phi_{\mathbf{P}}(a)).$$

This proves that $\Phi_{\mathbf{P}}$ is a hyperserial embedding of force $\boldsymbol{\nu}$. In order to see that it is unique, consider an embedding $\Phi' : \mathbb{T}_{\mathbf{P}} \longrightarrow \mathbf{Mo}$ satisfying the conditions. We have $\Phi'(\mathfrak{l} \circ \mathfrak{a}) = \mathfrak{l} \circ \Phi'(\mathfrak{a})$ for all $\mathfrak{l} \in \mathcal{L}_{<\omega^\mu}$ and $\mathfrak{a} \in (\mathfrak{M}_{\mathbf{P}})_{<\omega^\mu}$ where $\boldsymbol{\mu} \leqslant \boldsymbol{\nu}$. So it is enough to prove that $\Phi'(\mathfrak{m}_a) = \Phi_{\mathbf{P}}(\mathfrak{m}_a)$ for all $a \in \mathbf{No}$ which is $\Sigma$-nested for a certain $\Sigma \in \mathbf{P}$. For such $a = \Xi_\Sigma z$ where $z \in \mathbf{No}$, the number $\mathfrak{m}_a$ is $\Sigma_a^+$-nested. So

$$\Phi'(\mathfrak{m}_a) = \Xi_{\Sigma_a^+} \Phi_{\Sigma_a^+}(z) = \mathfrak{m}_{\Xi_{\Phi(\Sigma_a^+)} \Phi_{\Sigma_a^+}(z)} = \mathfrak{m}_{\Xi_{\Phi(\Sigma_a)} \Phi_{\Sigma_a}(z)} = \Phi_{\mathbf{P}}(\mathfrak{m}_a).$$

We deduce that $\Phi' = \Phi_{\mathbf{P}}$. □

# 6 Embedding theorems

## 6.1 Atomic embeddings

**Proposition 6.1.** *Let $\mathbb{T} = \mathbb{R}[[\mathfrak{M}]]$ be a subfield of force $\boldsymbol{\nu}$ and let $\Phi : \mathfrak{M}_{\boldsymbol{\lambda}} \longrightarrow \mathbf{Mo}_{\boldsymbol{\lambda}}$ be a strictly increasing function. If $\boldsymbol{\nu}$ is a successor, then we assume moreover that we have*

$$\Phi(L_{\boldsymbol{\lambda}/\omega}(\mathfrak{a})) = L_{\boldsymbol{\lambda}/\omega}(\Phi(\mathfrak{a})) \tag{6.1}$$

*for all $\mathfrak{a} \in \mathfrak{M}_{\boldsymbol{\lambda}}$. Then there is a unique extension $\hat{\Phi}$ of $\Phi$ into an embedding $\mathcal{H}_{\boldsymbol{\nu}}(\mathfrak{M}_{\boldsymbol{\lambda}}) \longrightarrow \mathbf{No}$ of force $\boldsymbol{\nu}$.*

**Proof.** Note that (6.1) insures that $\Phi(\mathfrak{M}_{\boldsymbol{\lambda}})$ satisfies the premises of Lemma 3.4. Consider the reciprocal $\Phi^{-1}$ of the co-restriction $\Phi : \mathfrak{M}_{\boldsymbol{\lambda}} \longrightarrow \Phi(\mathfrak{M}_{\boldsymbol{\lambda}})$. We have a group morphism

$$\Phi_1 : \mathfrak{L}_{\mathfrak{M}_{\boldsymbol{\lambda}}} \longrightarrow \mathfrak{L}_{\Phi(\mathfrak{M}_{\boldsymbol{\lambda}})} \, ; \, \tilde{\mathfrak{f}} \mapsto \widetilde{\mathfrak{f} \circ \Phi^{-1}},$$



which extends $\Phi$. Since $\Phi^{-1}$ is a strictly increasing bijection, for $\mathfrak{f} \in \mathcal{F}_{\mathfrak{M}_{\boldsymbol{\lambda}}}$, we have

$$\tilde{\mathfrak{f}} \succ 1 \Longleftrightarrow \mathfrak{f}_{\min S(\mathfrak{f})} \succ 1 \Longleftrightarrow (\mathfrak{f} \circ \Psi)_{\min S(\mathfrak{f} \circ \Psi)} \succ 1 \Longleftrightarrow \widetilde{\mathfrak{f} \circ \Psi} \succ 1.$$

So $\Phi_1$ is a strictly increasing group morphism $\mathfrak{L}_{\mathfrak{M}_{\boldsymbol{\lambda}}} \longrightarrow \mathbf{Mo}$. By , there is a unique strongly linear extension $\hat{\Phi}$ of $\Phi$, which is a field morphism. For $\tilde{\mathfrak{f}} \in \mathfrak{L}_{\mathfrak{A}}$, we have

$$\begin{aligned}
\hat{\Phi}(\log \tilde{\mathfrak{f}}) &= \hat{\Phi}(\Sigma_{\mathfrak{a} \in \mathfrak{A}} (\log \mathfrak{f}_{\mathfrak{a}}) \circ \mathfrak{a}) \\
&= \sum_{\mathfrak{a} \in \mathfrak{A}} \hat{\Phi}\left( \sum_{\gamma < \boldsymbol{\lambda}_{/\omega}} (\mathfrak{f}_{\mathfrak{a}})_{\gamma} \ell_{\gamma+1} \circ \mathfrak{a} \right) \\
&= \sum_{\mathfrak{a} \in \mathfrak{A}} \sum_{\gamma < \boldsymbol{\lambda}_{/\omega}} (\mathfrak{f}_{\mathfrak{a}})_{\gamma} \Phi_1(\ell_{\gamma+1} \circ \mathfrak{a}) \\
&= \sum_{\mathfrak{a} \in \mathfrak{A}} \sum_{\gamma < \boldsymbol{\lambda}_{/\omega}} (\mathfrak{f}_{\mathfrak{a}})_{\gamma} \ell_{\gamma+1} \circ \Phi(\mathfrak{a}) \\
&= \sum_{\mathfrak{a} \in \mathfrak{A}} (\log \mathfrak{f}_{\mathfrak{a}}) \circ \hat{\Phi}(\mathfrak{a}) \\
&= \sum_{\mathfrak{a} \in \mathfrak{A}} (\log (\mathfrak{f} \circ \Phi^{-1})_{\mathfrak{a}}) \circ \mathfrak{a} \\
&= \log \hat{\Phi}(\tilde{\mathfrak{f}}).
\end{aligned}$$

For $\eta < \boldsymbol{\nu}$ and $\mathfrak{b} \in \mathfrak{d}_{\omega^{\eta}}(\mathfrak{L}_{\mathfrak{A}})$, by (3.1), there are $\mathfrak{a} \in \mathfrak{A}$, $\gamma \gg \omega^{\eta-}$ and $n \in \mathbb{N}$ with $\mathfrak{b} = \ell_{\gamma+\omega^{\eta}-n} \circ \mathfrak{a}$. We have

$$\begin{aligned}
\tilde{\Phi}(\ell_{\omega^{\eta}} \circ \mathfrak{b}) &= \hat{\Phi}(\ell_{\gamma+\omega^{\eta}} \circ \mathfrak{a} - n) \\
&= \Phi_1(\ell_{\gamma} \circ (\ell_{\omega^{\eta}} \circ \mathfrak{a})) - n \\
&= \ell_{\gamma} \circ \Phi(\ell_{\omega^{\eta}} \circ \mathfrak{a}) - n & \text{(by definition of } \Phi_1) \\
&= \ell_{\gamma} \circ \ell_{\omega^{\eta}} \circ \Phi(\mathfrak{a}) - n & \text{(by (6.1))} \\
&= \ell_{\omega^{\eta}} \circ (\ell_{\gamma+\omega^{\eta}-n} \circ \Phi(\mathfrak{a})) \\
&= \ell_{\omega^{\eta}} \circ (\Phi(\ell_{\gamma+\omega^{\eta}-n} \circ \mathfrak{a})) & \text{(by definition of } \Phi_1) \\
&= \ell_{\omega^{\eta}} \circ \hat{\Phi}(\mathfrak{b}).
\end{aligned}$$

We deduce with Lemma 3.7 that $\hat{\Phi}$ is an embedding of force $\boldsymbol{\nu}$.

For $\mathfrak{f} \in \mathcal{F}_{\mathfrak{M}_{\boldsymbol{\lambda}}}$, there are $n \in \mathbb{N}$ and $\mathfrak{a}_0, \ldots, \mathfrak{a}_{n-1} \in S(\mathfrak{f})$ and $\mathfrak{l}_0, \ldots, \mathfrak{l}_{n-1} \in \mathfrak{L}_{<\boldsymbol{\lambda}}$ with $\tilde{\mathfrak{f}} = \prod_{i=0}^{n-1} \mathfrak{l}_i \circ \mathfrak{a}_i$. We thus have

$$\Psi(\tilde{\mathfrak{f}}) = \prod_{i=0}^{n-1} \Psi(\mathfrak{l}_i \circ \mathfrak{a}_i) = \prod_{i=0}^{n-1} \mathfrak{l}_i \circ \Psi(\mathfrak{a}_i) = \prod_{i=0}^{n-1} \mathfrak{l}_i \circ \Phi(\mathfrak{a}_i) = \hat{\Phi}(\tilde{\mathfrak{f}})$$

for any extension $\Psi$ of $\Phi$ which is a force $\boldsymbol{\nu}$ hyperserial embedding. So $\hat{\Phi}$ is unique. $\square$

## 6.2 General embedding theorem

Let $\boldsymbol{\nu} \leqslant \mathbf{On}$ with $\boldsymbol{\nu} > 0$, write $\boldsymbol{\lambda} := \omega^{\boldsymbol{\nu}}$ and let $\mathbb{T} = \mathbb{R}[[\mathfrak{M}]]$ be a confluent subfield of force $\boldsymbol{\nu}$. We define an increasing sequence $(\mathbb{T}^{[\rho]})_{\rho \in \mathbf{On}}$ of confluent subfields $\mathbb{T}^{[\rho]} = \mathbb{R}[[\mathfrak{M}^{[\rho]}]]$ of $\mathbf{No}$ of force $(\boldsymbol{\nu}, \boldsymbol{\nu})$ as follows. Set $\mathbb{T}^{[0]} := \widetilde{\mathcal{H}_{\boldsymbol{\nu}}(\mathfrak{M}_{\boldsymbol{\lambda}})}$. Let $\rho \in \mathbf{On}$ such that $\mathbb{T}^{[\iota]} = \mathbb{R}[[\mathfrak{M}^{[\iota]}]]$ is defined for all $\iota < \rho$.

If $\rho$ is a limit, then define $\mathfrak{M}^{[\rho]} := \bigcup_{\iota < \rho} \mathfrak{M}^{[\iota]}$. We set

$$\mathbb{T}^{[\rho]} := (\mathbb{R}[[\mathfrak{M}^{[\rho]}]]),$$



which is a confluent subfield of force $(\boldsymbol{\nu}, \boldsymbol{\nu})$.

If $\rho = \iota + 1$ is a successor, then consider the class $\mathbf{P}^{[\iota]}$ of all $\boldsymbol{\lambda}$-bounded positive good sequences which lie over $\mathbb{T}^{[\iota]}$, with $\mathbf{Ne}_{\Sigma \nearrow i} \cap \mathbb{T}^{[\iota]} = \varnothing$ for all $i \in \mathbb{N}$. We claim that $\mathbf{P}^{[\iota]}$ satisfies (5.1,5.2,5.3) with respect to $\mathbb{T}^{[\iota]}$. Indeed, for $\Sigma = (\varphi_i, \varepsilon_i, \psi_i, \iota_i, \alpha_i)_{i \in \mathbb{N}} \in \mathbf{P}^{[\iota]}$, we have

$$\mathbf{Ne}_{\Sigma \nearrow 1} \cap \mathbb{T}^{[\iota]} = L_{\alpha_0}\left(\frac{\mathbf{Ne}_\Sigma - \varphi_0}{\varepsilon_0 e^{\psi_0}}\right)^{\iota_0} \cap \mathbb{T}^{[\iota]} \subseteq \varnothing.$$

So $\Sigma_{\nearrow 1} \in \mathbf{P}_\iota$, so $\Sigma_{\nearrow i} \in \mathbf{P}_\iota$ for all $i \in \mathbb{N}$. Moreover, if $\alpha_0 = \beta \omega$ for a certain additively indecomposable ordinal $\beta_0$, then $\Sigma - n$ is good for all $n \in \mathbb{N}$. Let $n \in \mathbb{N}$ and assume for contradiction that there is $\mathfrak{m}$ in $\mathbb{T}^{[\iota]} \cap \mathbf{Ne}_{\Sigma - n}$. Then $E_{\beta n} \mathfrak{m} \in \mathbb{T}^{[\iota]} \cap \mathbf{Ne}_\Sigma$ because $\mathbb{T}^{[\iota]}$ has force $(\boldsymbol{\nu}, \boldsymbol{\nu})$: a contradiction. So $\Sigma - n \in \mathbf{P}^{[\iota]}$. We define

$$\mathbb{T}^{[\rho]} := ((\mathbb{T}^{[\iota]})_{\mathbf{P}^{[\iota]}})_{(<\boldsymbol{\nu})}.$$

We set $\mathbb{T}^{[\mathbf{On}]} := \bigcup_{\rho \in \mathbf{On}} \mathbb{T}^{[\rho]}$ and $\mathbb{T}' := \mathbb{T} \cap \mathbb{T}^{[\mathbf{On}]}$. We have $\mathbb{T}^{[\mathbf{On}]} = \mathbb{R}[[\mathfrak{M}^{[\mathbf{On}]}]]$ where $\mathfrak{M}^{[\mathbf{On}]} := \bigcup_{\rho \in \mathbf{On}} \mathfrak{M}^{[\rho]}$ and $\mathbb{T}^{[\mathbf{On}]}$ is a confluent subfield of force $(\boldsymbol{\nu}, \boldsymbol{\nu})$.

**Proposition 6.2.** *We have* $\mathbb{T}^{[\mathbf{On}]} = \mathbb{T}$.

**Proof.** We prove this by induction on the hyperserial complexity. Let $s \in \mathbb{T}$, write $\varsigma_{\boldsymbol{\nu}}(s) =: \eta$ and assume that we have $t \in \mathbb{T}^{[\mathbf{On}]}$ for all $t \in \mathbb{T}$ with $\varsigma_{\boldsymbol{\nu}}(t) < \eta$. If $\operatorname{supp} s$ has no least element, then we have $\operatorname{supp} s \subseteq \mathbb{T}^{[\mathbf{On}]}$ whence $s \in \mathbb{T}^{[\mathbf{On}]}$.

Otherwise, let $\mathfrak{m} = \min(\operatorname{supp} s, \succ)$ and write $\mathfrak{m} = e^\psi (L_\beta E_\alpha^u)^\iota$ as a hyperserial expansion. We have $\varsigma_{\boldsymbol{\nu}}(s_{\succ \mathfrak{m}}) < \eta$ so $s_{\succ \mathfrak{m}} \in \mathbb{T}^{[\mathbf{On}]}$, and we need only prove that $\mathfrak{m} \in \mathbb{T}^{[\mathbf{On}]}$. If $s_\mathfrak{m} \notin \{-1, 1\}$, then $\varsigma_{\boldsymbol{\nu}}(\mathfrak{m}) < \eta$ so $\mathfrak{m} \in \mathbb{T}^{[\mathbf{On}]}$. So we may assume that $s_\mathfrak{m} \in \{-1, 1\}$. We have also $\varsigma_{\boldsymbol{\nu}}(\psi) < \eta$ so $\psi \in \mathbb{T}^{[\mathbf{On}]}$ so $e^\psi \in \mathbb{T}^{[\mathbf{On}]}$. So we need only prove that $L_\beta E_\alpha^u \in \mathbb{T}^{[\mathbf{On}]}$. If $L_\beta E_\alpha^u \in \mathfrak{M}_{\boldsymbol{\lambda}}$, then we have $L_\beta E_\alpha^u \in \mathbb{T}_0 \subseteq \mathbb{T}^{[\mathbf{On}]}$. So we assume that $L_\beta E_\alpha^u$ is not $L_{<\boldsymbol{\lambda}}$-atomic. In particular, we have $\alpha < \boldsymbol{\lambda}$ and $\beta \ll \boldsymbol{\lambda}_{/\omega}$ and $(\alpha, \beta) \neq 0$. If $\beta \neq 0$, then we have $\varsigma_{\boldsymbol{\nu}}(\mathfrak{n}) < \eta$ where $\mathfrak{n} := E_\alpha^u$, so $\mathfrak{n} \in \mathbb{T}^{[\mathbf{On}]}$. But we have $\beta < \boldsymbol{\lambda}$, whence $L_\beta E_\alpha^u \in \mathbb{T}^{[\mathbf{On}]}$. So we may assume that $\beta = 0$. If $E_\alpha^u$ is not $\Sigma$-nested for a nested sequence $\Sigma$, then we have $\varsigma_{\boldsymbol{\nu}}(u) < \eta$ so $u \in \mathbb{T}^{[\mathbf{On}]}$, whence $E_\alpha^u \in \mathbb{T}^{[\mathbf{On}]}$ since $\alpha < \boldsymbol{\lambda}$.

Assume that $E_\alpha^u$ is $\Sigma$ nested for a certain nested sequence $\Sigma$. We have $\varsigma_{\boldsymbol{\nu}}(\varphi_i), \varsigma_{\boldsymbol{\nu}}(\psi_i) < \varsigma_{\boldsymbol{\nu}}(E_\alpha^u)$ for all $i \in \mathbb{N}$, so $\varphi_i, \psi_i \in \mathbb{T}^{[\mathbf{On}]}$ for all $i \in \mathbb{N}$. If there is a $i \in \mathbb{N}$, which we choose minimal, with $\alpha_{\Sigma, i} \geqslant \boldsymbol{\lambda}$, then we have $(E_\alpha^u)_{;i} \in \mathbb{T}^{[\mathbf{On}]}$ as above. Given $\rho \in \mathbf{On}$ with $(E_\alpha^u)_{;i} \in \mathbb{T}^{[\rho]}$ and $\varphi_j, \psi_j \in \mathbb{T}^{[\rho]}$ for all $j \leqslant i$, we obtain so $E_\alpha^u \in (\mathbb{T}^{[\rho]})_{(<\boldsymbol{\nu})} = \mathbb{T}^{[\rho]}$. So we may assume that $\Sigma$ is $\boldsymbol{\lambda}$-bounded.

Let $\iota \in \mathbf{On}$ be minimal such that there is a $j \in \mathbb{N}$ with with $\varphi_i, \psi_i \in \mathbb{T}^{[\iota]}$ for all $i \geqslant j$. Choosing $j$ minimal, we claim that $\Sigma_{\nearrow j} \in \mathbf{P}_\iota$. Indeed $\Sigma_{\nearrow j}$ is a good $\boldsymbol{\lambda}$-bounded sequence which lies over $\mathbb{T}_\iota$, so we need only prove that $\mathbf{Ne}_{\Sigma_{\nearrow j + i}} \cap \mathbb{T}^{[\iota]} = \varnothing$ for all $i \in \mathbb{N}$. Assume for contradiction that there is $i \in \mathbb{N}$, and $\mathfrak{m} \in \mathbf{Ne}_{\Sigma_{\nearrow j + i}} \cap \mathbb{T}^{[\iota]}$. We cannot have $\iota = 0$ by Lemma 3.5. Applying Lemma 3.11 to the path $P := (\mathfrak{m}_{;k} - \varphi_{\Sigma, i+j+k})_{k \in \mathbb{N}}$ in $\mathfrak{m}$, we see that there are a $k \in \mathbb{N}$ and a $\rho < \iota$ with $\mathfrak{m}_{;k} - \varphi_{\Sigma, i+j+k} \in \mathbb{T}^{[\rho]}$. But then for $l > k$, we have $\varphi_{\Sigma, j+i+l}, \psi_{\Sigma, j+i+l} \in \mathbb{T}^{[\rho]}$. This contradicts the minimality of $\iota$. Therefore $\Sigma \in \mathbf{P}_\iota$. We deduce that $E_\alpha^u \in \mathbb{T}^{[\iota+1]}$. This concludes the proof that $\mathbb{T} \subseteq \mathbb{T}^{[\mathbf{On}]}$. □

**Theorem 6.3.** *Let* $\mathbb{T} = \mathbb{R}[[\mathfrak{M}]]$ *be a subfield of force* $\boldsymbol{\nu}$ *and let* $\Phi : \mathfrak{M}_{\boldsymbol{\lambda}} \longrightarrow \mathbf{Mo}_{\boldsymbol{\lambda}}$ *be a function. If* $\boldsymbol{\nu}$ *is a successor, then we assume moreover that we have* $\Phi(L_{\boldsymbol{\lambda}/\omega}(\mathfrak{a})) = L_{\boldsymbol{\lambda}/\omega}(\Phi(\mathfrak{a}))$ *for all* $\mathfrak{a} \in \mathfrak{M}_{\boldsymbol{\lambda}}$. *Let* $\mathbf{P}$ *denote the class of* $\boldsymbol{\lambda}$-*bounded positive good sequences over* $\mathbb{T}$.



*Consider a family $(\Phi_\Sigma)_{\Sigma \in \mathbf{P}}$ of order isomorphisms $\Phi_\Sigma\colon \mathbf{No} \longrightarrow \mathbf{No}$ with*

$$\Phi_\Sigma(z) = \varepsilon_{\Sigma,1} \iota_{\Sigma,1} \Phi_{(\Sigma \nearrow 1)^+}(\varepsilon_{\Sigma,1} \iota_{\Sigma,1} z) \quad \text{and} \tag{6.2}$$

$$\Phi_{\Sigma-1}(z) = \Phi_\Sigma(z) \tag{6.3}$$

*for $\Sigma \in \mathbf{P}$ and $z \in \mathbf{No}$. There is a unique extension $\hat{\Phi}$ of $\Phi$ into a hyperserial embedding $\mathbb{T} \longrightarrow \mathbf{No}$ of force $\boldsymbol{\nu}$ with*

$$\hat{\Phi}(\Xi_\Sigma z) = \Xi_{\hat{\Phi}(\Sigma)} \Phi_\Sigma(z) \quad \text{for all } \Sigma \in \mathbf{P} \text{ and } z \in \mathbf{No} \text{ with } \Xi_\Sigma z \in \mathbb{T}.$$

**Proof.** By Propositions 3.9 and 3.10, we may assume that $\mathbb{T}$ has force $(\boldsymbol{\nu}, \boldsymbol{\nu})$. For each $\rho \in \mathbf{On}$, we define a hyperserial embedding $\Phi_\rho\colon \mathbb{T}^{[\rho]} \longrightarrow \mathbf{No}$ of force $\boldsymbol{\nu}$ such that $\Phi_\rho$ and $\Phi_\iota$ coincide on $\mathbb{T}^{[\iota]}$ whenever $\iota < \rho$ and that $\Phi_0$ and $\Phi$ coincide on $\mathfrak{M}_{\boldsymbol{\lambda}}$. Recall that $\mathbb{T}^{[0]} := \widetilde{\mathcal{H}_{\boldsymbol{\nu}}(\mathfrak{M}_{\boldsymbol{\lambda}})}$. By Proposition 6.1 we have a unique extension of $\Phi$ into an embedding $\Phi\colon \mathcal{H}_{\boldsymbol{\nu}}(\mathfrak{M}_{\boldsymbol{\lambda}}) \longrightarrow \mathbf{No}$ of force $\mathbf{On}$. We define $\Phi_0$ as the unique hyperserial embedding extending $\Phi$ of Proposition 3.10. Let $\rho \in \mathbf{On}$ such that $\Phi_\iota$ is defined for all $\iota < \rho$.

If $\rho$ is a limit, then $\Psi_\rho := \bigcup_{\iota < \rho} \Phi_\iota$ induces a unique hyperserial embedding $\Phi_\rho\colon \mathbb{T}^{[\rho]} \longrightarrow \mathbf{No}$ of force $\boldsymbol{\nu}$.

If $\rho = \iota + 1$ is a successor, then note that $\mathbf{P}^{[\iota]} \subseteq \mathbf{P}$. In view of (6.2, 6.3), we may apply Proposition 5.9 for $(\Phi_\iota)_\Sigma := \Phi_\Sigma$ for all $\Sigma \in \mathbf{P}^{[\iota]}$. We deduce that there is a unique extension $\Psi_\rho$ of $\Phi_\iota$ to $(\mathbb{T}^{[\iota]})_{\mathbf{P}^{[\iota]}}$ which satisfies (5.6). We then set $\Phi_\rho := (\Psi_\rho)_{(<\boldsymbol{\nu})}$ using Proposition 3.10.

Lastly, we define $\Phi_{\mathbf{On}} := \bigcup_{\rho \in \mathbf{On}} \Phi_\rho$. So $\Phi_{\mathbf{On}}$ is a hyperserial embedding of force $\boldsymbol{\nu}$. By Proposition 6.2, the map $\hat{\Phi} := \Phi_{\mathbf{On}} \upharpoonright \mathbb{T}$ is an embedding $\mathbb{T} \longrightarrow \mathbf{No}$ of force $\boldsymbol{\nu}$ extending $\Phi$. For $\Sigma \in \mathbf{P}$, as we have seen, there is $\iota \in \mathbf{On}$ with $\Sigma \in \mathbf{P}_\iota$. Given $z \in \mathbf{No}$ with $\Xi_\Sigma z \in \mathbb{T}$, we have

$$\hat{\Phi}(\Xi_\Sigma z) = (\Phi_\iota)_{\mathbf{P}^{[\iota]}}(\Xi_\Sigma z) = \Xi_{\Phi_i(\Sigma)} (\Phi_\iota)_\Sigma(z) = \Xi_{\hat{\Phi}(\Sigma)} \Phi_\Sigma(z).$$

So $\hat{\Phi}$ satisfies the condition.

Assume for contradiction that there is a distinct embedding $\Psi\colon \mathbb{T} \longrightarrow \mathbf{No}$ which extends $\Phi$ and satisfies the condition. Let $\eta \in \mathbf{On}$ be minimal such that there is $a \in \mathbb{T}$ with $\Psi(a) \neq \hat{\Phi}(a)$ and $\varsigma_{\boldsymbol{\nu}}(s) = \eta$. As usual we may assume that $a$ is a monomial with hyperserial expansion $a = L_\beta E_\alpha^u$. If $\alpha \geqslant \boldsymbol{\lambda}$, then we write $\beta = \beta' + \beta''$ where $\beta'' \ll \boldsymbol{\lambda}_{/\omega}$ and $\beta' \gg \boldsymbol{\lambda}_{/\omega}$. So $L_{\beta'} E_\alpha^u \in \mathfrak{M}_{\boldsymbol{\lambda}}$. We have $\Psi(a) = L_{\beta''} \Psi(L_{\beta'} E_\alpha^u) = L_{\beta''} \Phi(L_{\beta'} E_\alpha^u) = \hat{\Phi}(a)$: a contradiction. So $\alpha < \boldsymbol{\lambda}$. We may assume that $\varsigma_{\boldsymbol{\nu}}(u) = \eta$, so $\beta = 0$ and there is a $\boldsymbol{\lambda}$-bounded nested sequence $\Sigma$ such that $s$ is $\Sigma$-nested. Let $z \in \mathbf{No}$ with $a = \Xi_\Sigma z$. We deduce that $\Psi(s) = \Xi_{\Psi(\Sigma)} \Phi_\Sigma(z)$. But we have $\varsigma_{\boldsymbol{\nu}}(\varphi_{\Sigma,i}), \varsigma_{\boldsymbol{\nu}}(\psi_{\Sigma,i}) < \eta$ for all $i \in \mathbb{N}$, so $\Psi(\Sigma) = \hat{\Phi}(\Sigma)$. So $\Psi(a) = \Xi_{\Psi(\Sigma)} \Phi_\Sigma(z) = \Xi_{\hat{\Phi}(\Sigma)} \Phi_\Sigma(z) = \hat{\Phi}(a)$: a contradiction. $\square$

## 6.3 Automorphisms of No by scalar multiplication

We can readily define bijective embeddings $\mathbf{No} \longrightarrow \mathbf{No}$ of force $\mathbf{On}$. Let $f \in \mathbf{No}^{>}$. A simple example is that of embeddings whose effect on the hyperserial description of a number $a$ is to multiply each surreal label by $f$. Let $\mathbf{P}$ denote the class of all positive good sequences. The family $(\Phi_{f,\Sigma})_{\Sigma \in \mathbf{P}}$ of order isomorphisms $\mathrm{Id}_{f,\Sigma}\colon \mathbf{No} \longrightarrow \mathbf{No}$ given by

$$\forall z \in \mathbf{No}, \mathrm{Id}_{f,\Sigma}(z) := f z$$

satisfies (6.2) and (6.3). By Theorem 6.3, the mapping $\mathrm{Id}\colon \{\omega\} \longrightarrow \{\omega\}$ extends uniquely into an embedding $\mathrm{Id}_f\colon \mathbf{No} \longrightarrow \mathbf{No}$ of force $\mathbf{On}$ with $\mathrm{Id}_f(\Xi_\Sigma z) = \Xi_{\mathrm{Id}_f(\Sigma)} f z$ for all $\Sigma \in \mathbf{P}$ and $z \in \mathbf{No}$. Note that for $\Psi = \mathrm{Id}_f \circ \mathrm{Id}_{f^{-1}}$ or $\Psi = \mathrm{id}_{\mathbf{No}}$, the function $\Psi$ extends of $\mathrm{Id}$ into an embedding of force $\mathbf{On}$ with $\Psi(\Xi_\Sigma z) = \Xi_{\Psi(\Sigma)} z$ or all $\Sigma \in \mathbf{P}$ and $z \in \mathbf{No}$. We deduce by unicity that $\mathrm{Id}_f \circ \mathrm{Id}_{f^{-1}} = \mathrm{id}_{\mathbf{No}}$, that is, the function $\mathrm{Id}_f$ is bijective with functional inverse $\mathrm{Id}_{f^{-1}}$. Note lastly that $\mathrm{Id}_f$ fixes $\tilde{\mathbb{L}}$ pointwise.



# 7 Applications

We now apply our results by proving auxiliary facts that are crucial for defining derivations and composition laws on **No**. We will first give a characterization of the hyperexponential closure of a subclass (see below) in terms of properties of paths. We will then extend the method described in [6, Section 5] in order to adapt it to the presence of nested numbers. Lastly, we will define the fields of bounded surreal numbers and the canonical right compositions with sufficiently atomic elements on those fields.

## 7.1 Hyperexponential closure of a subclass

Let $\nu \leqslant \mathbf{On}$ with $\nu > 0$. Let $\mathbb{U} = \mathbb{R}[[\mathfrak{U}]]$ be a subfield of **No** of force $\nu$. Given a subclass $\mathbf{X} \subseteq \mathbb{U}$ with $\mathbf{X} \not\subseteq \mathbb{R}$, we write $\mathbf{X}_{(<\nu)}$ for the smallest subfield of $\mathbb{U}$ of force $(\nu, \nu)$ containing $\mathbf{X}$. We show that this can be constructed as the union of subgroups $\mathbb{R}[[\mathfrak{S}_{(\gamma)}]], \gamma \in \mathbf{On}$ where

$$\mathfrak{S}_{\langle 0 \rangle} := \bigcup_{s \in \mathbf{X}} \operatorname{supp} s$$

$$\mathfrak{S}_{\langle \gamma \rangle} := \mathfrak{S}_{\langle \gamma_- \rangle} \cup \bigcup_{\mu < \nu} \left( \{ (E_{\omega^\mu}^\varphi)^{\pm 1} : \varphi \in \mathbb{R}[[\mathfrak{S}_{\langle \gamma_- \rangle}]] \cap \mathbb{U}_{\succ, \omega^\mu} \} \cup \left( \bigcup_{\mathfrak{a} \in \mathfrak{S}_{\langle \gamma_- \rangle} \cap \mathfrak{U}_{\omega^\mu}} \operatorname{supp} L_{\omega^\mu}(\mathfrak{a}) \right) \right)$$

　　if $\gamma$ is a successor, and

$$\mathfrak{S}_{\langle \gamma \rangle} := \bigcup_{\rho < \gamma} \mathfrak{S}_{\langle \rho \rangle} \qquad \text{if } \gamma \text{ is a limit.}$$

**Proposition 7.1.** *The class* $\mathbf{X}_{(<\nu)} := \bigcup_{\gamma \in \mathbf{On}} \mathbb{R}[[\mathfrak{S}_{\langle \gamma \rangle}]]$ *is the smallest subfield of* $\mathbb{U}$ *of force* $(\nu, \nu)$ *that contains* $\mathbf{X}$.

**Proof.** It is easy to see by induction that any subfield of $\mathbb{U}$ of force $(\nu, \nu)$ that contains $\mathbf{X}$ must also contain $\mathbf{X}_{(<\nu)}$. Writing $\mathfrak{S}_{(<\nu)} := \bigcup_{\gamma \in \mathbf{On}} \mathfrak{S}_{\langle \gamma \rangle}$, we have $\mathbf{X}_{(<\nu)} := \mathbb{R}[[\mathfrak{S}_{(<\nu)}]]$ by [11, Lemma 2.1]. The class $\mathbf{X}_{(<\nu)}$ is closed under logarithms of infinite monomials, additive opposites, exponentials of 1-truncated series and reciprocals of those exponentials. Therefore $\mathfrak{S}_{(<\nu)}$ is a group, so $\mathbf{X}_{(<\nu)}$ is a field of well-based series. We also deduce that $\mathbf{X}_{(<\nu)}$ is closed under exponentials, and logarithms and real powers of strictly positive elements. Given $s \in \mathbf{X} \setminus \mathbb{R}$, we there is an $\mathfrak{u} \in \operatorname{supp} s \cap \mathfrak{U}^{\neq}$, whence $\mathfrak{S}_{(<\nu)} \supsetneq \{1\}$, i.e. $\mathbf{X}_{(<\nu)}$ is non-trivial.

　　It is enough in order to show that $\mathbf{X}_{(<\nu)}$ is a subfield of force $\nu$ to show that it is closed under all hyperlogarithms $L_{\omega^\mu}, \mu < \nu$. We prove this by induction on $\mu < \nu$. We already dealt with the case $\mu = 0$, so we may assume that $\mu \geqslant 1$ and that $\mathbf{X}_{(<\nu)}$ is closed under $L_{\omega^\eta}$ for all $\eta < \mu$. But then $\mathbf{X}_{(<\nu)}$ is also closed under the action of $\mathbb{L}_{<\omega^\mu}$. Given $s \in \mathbf{X}_{(<\nu)}$ and $\mu < \nu$, we have $L_\gamma(s) - L_\gamma(\mathfrak{a}) \prec 1$ for an $L_{<\omega^\mu}$-atomic element $\mathfrak{a} \in \mathfrak{U}_{\omega^\mu}$. Our induction hypothesis implies that $L_\gamma(\mathfrak{a})$, as the dominant monomial of $L_\gamma(s)$, lies in $\mathbf{X}_{(<\nu)}$. So $\varepsilon := L_\gamma(s) - L_\gamma(\mathfrak{a})$ lies in $\mathbf{X}_{(<\nu)}$ as well. Recall that

$$L_{\omega^\mu}(s) = L_{\omega^\mu}(\mathfrak{a}) + \sum_{k>0} \frac{(\ell_{\omega^\mu}^{\uparrow \gamma})^{(k)} \circ L_\gamma(\mathfrak{a})}{k!} \varepsilon^k,$$

where $(\ell_{\omega^\mu}^{\uparrow \gamma})^{(k)} \in \mathbb{L}_{<\omega^\mu}$ for each $k > 0$. So $(\ell_{\omega^\mu}^{\uparrow \gamma})^{(k)} \circ L_\gamma(\mathfrak{a}) \in \mathbf{X}_{(<\nu)}$. We deduce that $L_{\omega^\mu}(s) \in \mathbf{X}_{(<\nu)}$, which concludes our proof that $\mathbf{X}_{(<\nu)}$ is a subfield of force $\nu$. Now $\mathbf{X}_{(<\nu)}$ is closed under hyperexponentials of truncated elements, so $\mathbf{X}_{(<\nu)}$ has force $(\nu, \nu)$. □



## 7.2 Characterizing closure in terms of paths

**Lemma 7.2.** *Assume that* $\operatorname{supp} L_{\omega^\mu}(\mathfrak{a}) \subseteq \mathfrak{S}_{\langle 0 \rangle}$ *for all* $\mathfrak{a} \in \mathfrak{U}_{\omega^\mu} \cap \mathbf{X}$ *and all* $\mu < \boldsymbol{\nu}$. *Let* $\mathfrak{m} \in \mathfrak{S}_{(<\boldsymbol{\nu})}$ *and let $P$ be a $\boldsymbol{\lambda}$-bounded infinite path in $\mathfrak{m}$. Then we have $\mathfrak{m}_{P,i} \in \mathfrak{S}_{\langle 0 \rangle}$ for large enough $i \in \mathbb{N}$.*

**Proof.** Note that the hypothesis that $\operatorname{supp} L_{\omega^\mu}(\mathfrak{a}) \subseteq \mathfrak{S}_{\langle 0 \rangle}$ for all $\mathfrak{a} \in \mathfrak{U}_{\omega^\mu} \cap \mathbf{X}$ and all $\mu < \boldsymbol{\nu}$ implies that for any infinite $\boldsymbol{\lambda}$-bounded path $Q$ and any $i \in \mathbb{N}$, we have

$$\mathfrak{m}_{Q,i} \in \mathfrak{S}_{\langle 0 \rangle} \Longrightarrow (\forall j \geqslant i, (\mathfrak{m}_{Q,j} \in \mathfrak{S}_{\langle 0 \rangle})). \tag{7.1}$$

We will prove the statement of the lemma by induction on the least $\gamma \in \mathbf{On}$ with $\mathfrak{m} \in \mathfrak{S}_\gamma$. Consider an ordinal $\gamma$ such that the result holds for all $\mathfrak{n} \in \mathbf{Ne}_\Sigma \cap \mathfrak{S}_\gamma$ and assume that $\mathfrak{m} \in \mathfrak{S}_{\langle \gamma \rangle} \setminus \bigcup_{\rho < \gamma} \mathfrak{S}_{\langle \rho \rangle}$. Note that $\gamma$ is either $0$ or a successor. If $\gamma = 0$ then $\mathfrak{m} \in \mathbf{X}$ so we are done. So we assume that $\gamma = \rho + 1$ is a successor. By definition of $\mathfrak{S}_{\langle \rho+1 \rangle}$, there are a $\mu < \boldsymbol{\nu}$, an $\mathfrak{a} \in \mathfrak{S}_{\langle \rho \rangle} \cap \mathfrak{U}_{\omega^\mu}$ and a $\varphi \in \mathbb{R}[[\mathfrak{S}_{\langle \rho \rangle}]] \cap \mathbb{U}_{\succ, \omega^\mu}$ such that $\mathfrak{m} \in \operatorname{supp} L_{\omega^\mu}(\mathfrak{a})$ or that $\mathfrak{m} = E^\varphi_{\omega^\mu}$. In the first case, write $\sigma := L_{\omega^\mu}(\mathfrak{a})$. By [10, Corollary 5.17], there are a $\boldsymbol{\lambda}$-bounded infinite path $Q$ in $\mathfrak{d}_{E^\sigma_{\omega^\mu}} = \mathfrak{a}$ and a $k > 0$ with $P_{\nearrow 1} = Q_{\nearrow k}$. The induction hypothesis yields an $i \in \mathbb{N}$ with $\mathfrak{m}_{Q,i} \in \mathfrak{S}_{\langle 0 \rangle}$. We conclude with (7.1).

Consider the second case, when $\mathfrak{m} = E^\varphi_{\omega^\mu}$, so $\varphi = L_{\omega^\mu}(\mathfrak{m})$. By [10, Corollary 5.17], there are a path $R$ in $\varphi$ and an $l > 0$ with $P_{\nearrow 1} = R_{\nearrow l}$. As previously, the induction hypothesis gives an $i \in \mathbb{N}$ with $\mathfrak{m}_{R,i} \in \mathfrak{S}_{\langle 0 \rangle}$. We conclude with [10, Corollary 5.17] and (7.1). □

**Proposition 7.3.** *Let $\mathfrak{S} \subseteq \mathfrak{U}$ be a subclass with $\operatorname{supp} L_{\omega^\mu}(\mathfrak{a}) \subseteq \mathfrak{S}$ for all $\mathfrak{a} \in \mathfrak{U}_{\omega^\mu} \cap \mathfrak{S}$ and all $\mu < \boldsymbol{\nu}$. Set $\mathbb{G} := \mathbb{R}[[\mathfrak{S}]]$. Assume that each atomic element of $\mathbb{U}$ lies in $\mathbb{G}$. Let $\mathbf{C}$ denote the class of elements $s \in \mathbb{U}$ such that any maximal $\boldsymbol{\lambda}$-bounded path $P$ in $s$ satisfies $u_{P,i}, \psi_{P,j} \in \mathbb{G}$ for a certain $i \leqslant |P|$. Then $\mathbf{C} = \mathbb{G}_{(<\boldsymbol{\nu})}$.*

**Proof.** We first prove that each $s \in \mathbf{C}$ lies in $\mathbb{G}_{(<\boldsymbol{\nu})}$ by induction on $\varsigma_{\boldsymbol{\nu}}(s)$. Let $\gamma$ be an ordinal such that any $t \in \mathbf{C}$ with $\varsigma_{\boldsymbol{\nu}}(t) < \gamma$ lies in $\mathbb{G}_{(<\boldsymbol{\nu})}$, and let $s \in \mathbf{C}$ with $\varsigma_{\boldsymbol{\nu}}(s) = \gamma$. We may assume that $s = \mathfrak{m}$ is a monomial with $\mathfrak{m} \notin \mathbb{G}$. Since $\mathbb{G}_{(<\boldsymbol{\nu})}$ contains all atomic elements of $\mathbb{U}$, we may assume that $\mathfrak{m}$ is not atomic, so any maximal $\boldsymbol{\lambda}$-bounded path in $\mathfrak{m}$ has length $\geqslant 1$. Write $\mathfrak{m} = e^\psi (L_\beta E^u_\alpha)^\iota$ as a hyperserial expansion. If $P$ is a maximal $\boldsymbol{\lambda}$-bounded path in $\psi$, then $(\mathfrak{m}) * P$ is a $\boldsymbol{\lambda}$-bounded maximal path in $\mathfrak{m}$ with $u_{P,i}, \psi_{P,j} \in \mathbb{G}$ for a certain $i \leqslant |P|$. Since $\mathfrak{m} \notin \mathbb{G}$, the number $i$ cannot be zero, therefore $\psi \in \mathbf{C}$. But $\varsigma_{\boldsymbol{\nu}}(\psi) < \gamma$, so the induction hypothesis yields $\psi \in \mathbb{G}_{(<\boldsymbol{\nu})}$, whence $e^\psi \in \mathbb{G}_{(<\boldsymbol{\nu})}$. It remains to show that $L_\beta E^u_\alpha \in \mathbb{G}_{(<\boldsymbol{\nu})}$. Again we may assume that $\alpha < \boldsymbol{\lambda}$. If $\mathfrak{m}$ is not $\Sigma$-nested for a nested and $\boldsymbol{\lambda}$-bounded sequence $\Sigma$, then $\varsigma_{\boldsymbol{\nu}}(u) < \gamma$ and any maximal $\boldsymbol{\lambda}$-bounded path in $u$ similarly yields a maximal $\boldsymbol{\lambda}$-bounded path in $\mathfrak{m}$, so $E^u_\alpha \in \mathbb{G}_{(<\boldsymbol{\nu})}$, whence $L_\beta E^u_\alpha \in \mathbb{G}_{(<\boldsymbol{\nu})}$. We next assume that $\mathfrak{m}$ is $\Sigma$-nested for a nested and $\boldsymbol{\lambda}$-bounded sequence $\Sigma = (\varphi_i, \varepsilon_i, \psi_i, \iota_i, \alpha_i)_{i \in \mathbb{N}}$, whence in particular $\beta = 0$. Let $i \in \mathbb{N}$ and let $P_i$ be a maximal $\boldsymbol{\lambda}$-bounded path in $\varphi_i$ or in $\psi_i$. Then there is a finite path $P_{i,0}$ in $\mathfrak{m}$ such that $P_{i,0} * P_i$ is a maximal $\boldsymbol{\lambda}$-bounded path in $\mathfrak{m}$. We deduce again that $\varphi_i, \psi_i \in \mathbf{C}$. Since $\varsigma_{\boldsymbol{\nu}}(\varphi_i), \varsigma_{\boldsymbol{\nu}}(\psi_i) < \gamma$, the induction hypothesis gives $\varphi_i, \psi_i \in \mathbb{G}_{(<\boldsymbol{\nu})}$.

Consider the rightmost path $P$ in $\mathfrak{m}$. Since $P$ is infinite, we have $\tau_{P,i} \in \mathbb{G}$ for a certain $i \in \mathbb{N}$. But then $E^u_\alpha$ lies in $\mathbb{G}_{(<\boldsymbol{\nu})}$ since $\mathbb{G}_{(<\boldsymbol{\nu})}$ has force $(\boldsymbol{\nu}, \boldsymbol{\nu})$ and $\varphi_j, \psi_j \in \mathbb{G}_{(<\boldsymbol{\nu})}$ for all $j \leqslant i$. This concludes the inductive proof that $\mathbf{C} \subseteq \mathbb{G}_{(<\boldsymbol{\nu})}$.

We now prove that $s \in \mathbf{C}$ for each $s \in \mathbb{G}_{(<\boldsymbol{\nu})}$ by induction on $\varsigma_{\boldsymbol{\nu}}(s)$. Let $\gamma \in \mathbf{On}$ such that any $t \in \mathbb{G}_{(<\boldsymbol{\nu})}$ with $\varsigma_{\boldsymbol{\nu}}(t) < \gamma$ lies in $\mathbf{C}$, and let $s \in \mathbb{G}_{(<\boldsymbol{\nu})}$ with $\varsigma_{\boldsymbol{\nu}}(s) = \gamma$. Note that $\mathbf{C} = \mathbb{R}[[\mathfrak{N}]]$ where $\mathfrak{N} = \mathbf{C} \cap \mathfrak{U}$. Therefore, we may assume that $s = \mathfrak{n}$ is a monomial. Write

$$\mathfrak{n} = e^\psi (L_\beta E^u_\alpha)^\iota$$



as a hyperserial expansion. If $|P|=0$ then $\mathfrak{n}$ is atomic, so $\mathfrak{n} \in \mathbb{G} \subseteq \mathbb{G}_{(<\nu)}$. Assume that $|P|>0$ and $\tau_{P,1} \in \operatorname{supp} \psi$. Then since $\varsigma_{\nu}(\psi) < \gamma$, the path $P_{\nearrow 1}$ in $\psi$ is maximal and $\boldsymbol{\lambda}$-bounded. So there is an $i \in \mathbb{N}$ with $u_{P,i+1} = u_{P_{\nearrow 1},i} \in \mathbb{G}$ and $\psi_{P,i+1} = \psi_{P_{\nearrow 1},i} \in \mathbb{G}$.

If $|P|>0$ and $\tau_{P,1} \in \operatorname{supp} u$, then in particular $\alpha < \boldsymbol{\lambda}$. Assume that $(L_\beta E_\alpha^u)^\iota$ is $\Sigma$-nested for a $\boldsymbol{\lambda}$-bounded sequence $\Sigma = (\varphi_i, \varepsilon_i, \psi_i, \iota_i, \alpha_i)_{i \in \mathbb{N}}$, so $\beta = 0$. Note that $\varsigma_{\nu}(\varphi_i)$, $\varsigma_{\nu}(\psi_i) < \gamma$ for all $i \in \mathbb{N}$, so $\varphi_i, \psi_i \in \mathbf{C}$ for all $i \in \mathbb{N}$. Thus it suffices to show that the right-most path $R$ in $(E_\alpha^u)^\iota$ satisfies $\tau_{R,i} \in \mathbb{G}$ for a certain $i \in \mathbb{N}$. But this follows from Lemma 7.2. Assume now that $(L_\beta E_\alpha^u)^\iota$ is not $\Sigma$-nested for a $\boldsymbol{\lambda}$-bounded sequence $\Sigma$ or that $\beta \neq 0$. Then $\varsigma_{\nu}(u) < (L_\beta E_\alpha^u)^\iota$, so $u \in \mathbf{C}$. As previously, we have $u_{P,i+1} = u_{P_{\nearrow 1},i} \in \mathbb{G}$ and $\psi_{P,i+1} = \psi_{P_{\nearrow 1},i} \in \mathbb{G}$ for a certain $i < |P|$. This concludes our inductive proof that $\mathfrak{n} \in \mathbf{C}$. Therefore $\mathbf{C} = \mathbb{G}_{(<\nu)}$. □

## 7.3 A proof method, revisited

We next generalize the approach used throughout [6] (see [6, Section 5]), as an important step toward defining the composition law on **No**.

**Theorem 7.4.** *Let $\mathcal{G}$ be a subgroup of $(\mathbb{U}, +)$ with $1 \in \mathcal{G}$. Consider the class $\mathbb{E}(\mathcal{G}) = \mathbb{R}[[\mathfrak{A}(\mathcal{G})]]$ where $\mathfrak{A}(\mathcal{G})$ is the class of $\log$-atomic monomials $L_\gamma E_\alpha^\varphi \in \mathfrak{U}$ where $\alpha$ is additively indecomposable with $\omega \leqslant \alpha < \boldsymbol{\lambda}$, $\gamma \omega < \alpha$, and $\varphi \in \mathcal{G} \cap \mathbb{U}_{\succ, \alpha}$. Assume that $\mathbb{G}_{(<1)} \subseteq \mathcal{G}$ for all transserial subgroups of $\mathbb{U}$ with $\mathbb{G} \subseteq \mathcal{G}$ and that $\mathbb{E}(\mathcal{G}) \subseteq \mathcal{G}$.*

*Let $\mathfrak{S} \subseteq \mathfrak{U}$ be a subclass with $\{1\} \subsetneq \mathfrak{S}$ such that $\operatorname{supp} L_{\omega^\mu}(\mathfrak{S} \cap \mathfrak{U}_{\omega^\mu}) \subseteq \mathfrak{S}$ for all $\mu < \nu$. Write $\mathbb{G} := \mathbb{R}[[\mathfrak{S}]]$ and $\mathbb{H} := \mathbb{E}(\mathcal{G}) + \mathbb{G} = \mathbb{R}[[\mathfrak{A}(\mathcal{G}) \cup \mathfrak{S}]]$. If $\mathbb{G} \subseteq \mathcal{G}$, then*

$$\mathbb{G}_{(<\nu)} = \mathbb{H}_{(<1)} \subseteq \mathcal{G}.$$

**Proof.** Since $1 \in \mathcal{G}$, we have $\mathcal{G} - 1 = \mathcal{G}$, and consequently the class $\mathbb{E}(\mathcal{G})$ is a transserial subgroup of $\mathbb{U}$. It follows that $\mathbb{H}$ is a transserial subgroup of $\mathbb{U}$. Since $\mathbb{E}(\mathcal{G})$ and $\mathbb{G}$ are contained in $\mathcal{G}$, so is $\mathbb{H}$. Therefore $\mathbb{H}_{(<1)} \subseteq \mathcal{G}$.

We claim that $\mathbb{H}_{(<1)} = \mathbb{G}_{(<\nu)}$. We have $\mathbb{H}_{(<1)} \subseteq \mathbb{G}_{(<\nu)}$ by Proposition 7.1. We prove the inclusion $\mathbb{G}_{(<\nu)} \subseteq \mathbb{H}_{(<1)}$ by induction on the hyperserial complexity. Let $\gamma \in \mathbf{On}$ such that any $t \subseteq \mathbb{G}_{(<\nu)}$ with $\varsigma_{\nu}(t) < \gamma$ lies in $\mathbb{H}_{(<1)}$ and let $s \in \mathbb{G}_{(<\nu)}$ with $\varsigma_{\nu}(s) = \gamma$. We may assume that $s = \mathfrak{m}$ is a monomial. Write $\mathfrak{m} = e^\psi (L_\beta E_\alpha^u)^\iota$ as a hyperserial expansion. We have $\varsigma_{\nu}(\psi) < \gamma$ so $\psi, e^\psi \in \mathbb{H}_{(<1)}$, and we may assume that $\mathfrak{m} = L_\beta E_\alpha^u$.

If $\mathfrak{m}$ is $\Sigma$-nested for a certain $\boldsymbol{\lambda}$-bounded nested sequence $\Sigma = (\varphi_i, \varepsilon_i, \psi_i, \iota_i, \alpha_i)$, then applying Lemma 7.2 to the right-most path $R$ in $\mathfrak{m}$, we obtain an $i \in \mathbb{N}$ with $\mathfrak{m}_{R,i} \in \mathbb{G}$. Therefore $\mathfrak{m}_{R,i} \in \mathcal{G}$. Now for $j \leqslant i$, we have $\varsigma_{\nu}(\varphi_j), \varsigma_{\nu}(e^{\psi_j}) < \gamma$, so $\varphi_j, e^{\psi_j} \in \mathbb{H}_{(<1)}$, whence in particular $\varphi_{ij}, e^{\psi_j} \in \mathcal{G}$. We deduce by induction on $j \leqslant i$ that $\mathfrak{m}_{;i-j} \in \mathbb{E}(\mathcal{G})$, whence in particular $\mathfrak{m} \in \mathbb{E}(\mathcal{G})$. So $\mathfrak{m} \in \mathbb{H}_{(<1)}$.

Assume now that $\mathfrak{m}$ is $\Sigma$-nested for no $\boldsymbol{\lambda}$-bounded nested sequence $\Sigma$. This implies that $\varsigma_{\nu}(u) < \gamma$, so $u \in \mathcal{G}$. But then $\mathfrak{m} = L_\beta E_\alpha^u \in \mathbb{E}(\mathcal{G})$, so $\mathfrak{m} \in \mathbb{H}_{(<1)}$. This concludes our proof that $\mathbb{H}_{(<1)} = \mathbb{G}_{(<\nu)}$. □

**Proposition 7.5.** *Let $\rho$ be a limit ordinal and let $\mathfrak{U}^{[\sigma]}, \sigma < \rho$ be an increasing family of subgroups of $\mathfrak{U}$ such that each $\mathbb{U}^{[\sigma]} := \mathbb{R}[[\mathfrak{U}^{[\sigma]}]]$ is a subfield of $\mathbb{U}$ of force $(\nu, \nu)$. Let $\triangle: \bigcup_{\sigma < \rho} \mathbb{U}^{[\sigma]} \longrightarrow \mathbf{No}$ such that each restriction $\triangle \restriction \mathbb{U}^{[\sigma]} : \mathbb{U}^{[\sigma]} \longrightarrow \mathbf{No}$ for $\sigma < \lambda$ is a right composition of force $\nu$ as per [6, Definition 7.2] with a well-based relative near-support. Finally, consider a subclass $\mathfrak{S} \subseteq \mathfrak{U}$ with $\{1\} \subsetneq \mathfrak{S}$ such that $\operatorname{supp} L_{\omega^\mu}(\mathfrak{S} \cap \mathfrak{U}_{\omega^\mu}) \subseteq \mathfrak{S}$ for all $\mu < \nu$. If $\triangle \restriction \mathfrak{S}$ is well-based, then $\triangle \restriction \mathfrak{S}_{(<\nu)}$ is well-based.*



**Proof.** Set $\mathbb{G} := \mathbb{R}[[\mathfrak{S}]]$. Write $\mathcal{G}$ for the class of series $s \in \mathbb{G}_{(<\boldsymbol{\nu})}$ such that $(\triangle(\mathfrak{m}))_{\mathfrak{m} \in \operatorname{supp} s}$ is summable. So $\mathbb{G} \subseteq \mathcal{G}$. Each element of $\mathfrak{A}(\mathcal{G})$ is log-atomic, so [6, Lemma 1.16] implies that $\triangle$ is well-based on $\mathfrak{A}(\mathcal{G})$. Thus $\mathbb{E}(\mathcal{G}) \subseteq \mathcal{P}$. Given any transserial subgroup $\mathbb{T} = \mathbb{R}[[\mathfrak{X}]] \subseteq \mathcal{P}$, the function $\triangle \restriction \mathbb{H}$ extends uniquely into a transserial right composition $\triangle_1$ on $\mathbb{T}_{(<1)}$ by [6, Proposition 4.7]. But $\triangle_1$ must coincide with $\triangle$, which implies that $\triangle$ is strongly linear on $\mathbb{T}_{(<1)}$, whence $\mathbb{T}_{(<1)} \subseteq \mathcal{P}$. We can thus apply Theorem 7.4 and obtain that $\mathbb{G}_{(<\boldsymbol{\nu})} \subseteq \mathcal{P}$, hence that $\triangle \restriction \mathfrak{S}_{(<\boldsymbol{\nu})}$ is well-based. $\square$

## 7.4 Surreal numbers of bounded strength

Let $\nu$ be a non-zero ordinal and write $\lambda := \omega^\nu$. In this section, we define the force $\nu$ analogue $\mathbf{No}_\lambda$ of surreal numbers. Roughly speaking, the field $\mathbf{No}_\lambda$ contains all numbers which can be constructed using hyperexponentials and hyperlogarithms of force $< \nu$, and arbitrary surreal indexes for nested numbers.

Set $\mathbb{T}^{[0]} := \widetilde{\mathbb{L}_{<\lambda}} \circ \omega$. As in Section 6.2, we have a tower of confluent subfields $\mathbb{T}^{[\rho]}$, $\rho \in \mathbf{On}$ of force $(\nu, \nu)$. We define $\mathbf{No}_\lambda$ to be the field $\mathbb{T}^{[\mathbf{On}]} = \bigcup_{\rho \in \mathbf{On}} \mathbb{T}^{[\rho]}$, which is thus a confluent subfield of force $(\nu, \nu)$. We write $\mathbf{Mo}(\lambda)$ for the group of monomials in $\mathbf{No}_\lambda$, i.e. $\mathbf{No}_\lambda = \mathbb{R}[[\mathbf{Mo}(\lambda)]]$.

**Lemma 7.6.** *The $L_{<\lambda}$-atomic elements in $\mathbf{No}_\lambda$ are $\omega$ if $\nu$ is a limit, and among $\{L_{(\lambda/\omega)n}\omega, E_{(\lambda/\omega)n}\omega : n \in \mathbb{N}\}$ if it is a successor.*

**Proof.** Assume that $\nu$ is a limit. By [6, Lemma 7.25], the only atomic element of $\mathbb{L}_{<\lambda} \circ \omega$ is $\omega$. It follows by induction using [6, Lemma 7.25] and Lemma 5.7 that $\omega$ is the only atomic element of $\mathbf{No}_\lambda$.

Assume now that $\nu = \mu + 1$ is a successor. Then by [6, Lemma 7.25], the set $\{L_{(\lambda/\omega)n}\omega, E_{(\lambda/\omega)n}\omega : n \in \mathbb{N}\}$ is the set of atomic elements of $\mathbb{L}_{<\lambda} \circ \omega$. It follows by induction using [11, p. 66] and Lemma 5.7 that $\{L_{(\lambda/\omega)n}\omega, E_{(\lambda/\omega)n}\omega : n \in \mathbb{N}\}$ is the set of atomic elements of $\mathbf{No}_\lambda$. $\square$

**Lemma 7.7.** *Consider a property $\mathcal{X}$ of paths such that if every path in a number $a \in \mathbf{No}$ satisfies $\mathcal{X}$, then every subpath in $a$ satisfies $\mathcal{X}$. Assume that every $L_{<\lambda}$-atomic number whose paths all satisfy $\mathcal{X}$ lie in $\mathbf{No}_\lambda$. Let $a \in \mathbf{No}$ such that every path in $a$ satisfies $\mathcal{X}$. Then $a \in \mathbf{No}_\lambda$.*

**Proof.** We prove this by induction on $\varsigma_\nu(a)$. Write $\mathbf{X}$ for the class of numbers whose paths all satisfy $\mathcal{X}$. Let $\gamma \in \mathbf{On}$ such that $a \in \mathbf{X}$ with $\varsigma_\nu(a) < \gamma$ lie in $\mathbf{No}_\lambda$. Let $a \in \mathbf{No}$ with $\varsigma_\nu(a) = \gamma$. We may assume that $a = \mathfrak{m}$ is a non-trivial monomial. If $\mathfrak{m}$ is $L_{<\lambda}$-atomic, then $\mathfrak{m} \in \mathbf{No}_\lambda$ by the hypothesis on $\mathcal{X}$. We next assume that $\mathfrak{m}$ is not $L_{<\lambda}$-atomic. Write $\mathfrak{m} = e^\psi (L_\beta E_\alpha^u)^\iota$ as a hyperserial expansion. Every path $P$ in $\psi$ or $u$, is a subpath in $\mathfrak{m}$. Since $\varsigma_\nu(\psi) < \gamma$, we deduce that $\psi \in \mathbf{No}_\lambda$. If $\beta \neq 0$ or $\mathfrak{m}$ is not $\Sigma$-nested for a $\lambda$-bounded nested sequence $\Sigma$, then we also have $\varsigma_\nu(u) < \gamma$ whence $u \in \mathbf{No}_\lambda$, whence $\mathfrak{m} \in \mathbf{No}_\lambda$. So we may assume that $\mathfrak{m}$ is $\Sigma$-nested for a $\lambda$-bounded nested sequence $\Sigma = (\varphi_i, \varepsilon_i, \psi_i, \iota_i, \alpha_i)_{i \in \mathbb{N}}$. Given $i \in \mathbb{N}$, we have $\varsigma_\nu(\varphi_i), \varsigma_\nu(\psi_i) < \gamma$. Furthermore any path in $\varphi_i$ or $\psi_i$ is a subpath in $\mathfrak{m}$. We deduce by the induction hypothesis that $\varphi_i, \psi_i \in \mathbf{No}_\lambda$. Let $\iota \in \mathbf{On}$ be minimal such that $\varphi_i, \psi_i \in \mathbb{T}^{[\iota]}(\nu)$ for large enough $i \in \mathbb{N}$, and fix the least corresponding index $i_0 \in \mathbb{N}$. So $\Sigma_{\nearrow i_0} \in \mathbf{P}^{[\iota]}$ where $\mathbb{T}^{[\iota+1]}(\nu) = ((\mathbb{T}^{[\iota]}(\nu))_{\mathbf{P}^{[\iota]}})_{(<\nu)}$. But $\mathfrak{m} \in (\mathbb{T}^{[\iota]}(\nu))_{\mathbf{P}^{[\iota]}}$, so $\mathfrak{m} \in \mathbf{No}_\lambda$. $\square$

**Lemma 7.8.** *Let $\alpha \in \omega^{\mathbf{On}}$ and $\mathfrak{a} \in \mathbf{Mo}_{\alpha\omega}$. Then one of the following cases occurs:*

  a) *The hyperserial expansion of $E_\alpha^\mathfrak{a}$ is $E_{\alpha\omega}^{L_{\alpha\omega}\mathfrak{a}+1}$.*



    b) *There are a $k \in \mathbb{N}^>$ and a $\mathfrak{b} \in \mathbf{Mo}_{\alpha\omega^2}$ with $\mathfrak{a} = L_{\alpha k} \mathfrak{b}$.*

**Proof.** Note that $E_{\alpha\omega}^{L_{\alpha\omega}\mathfrak{a}+1} = E_\alpha^{\mathfrak{a}}$. If this is a hyperserial expansion then we are done. Otherwise, we have $E_\alpha^{\mathfrak{a}} \in L_{<\alpha\omega} \mathbf{Mo}_{\alpha\omega^2}$ so there is $n \in \mathbb{N}$ with $\mathfrak{b} := E_{\alpha n} E_\alpha^{\mathfrak{a}} \in \mathbf{Mo}_{\alpha\omega^2}$. Setting $k := n + 1 \in \mathbb{N}$, we have $\mathfrak{a} = L_{\alpha k} \mathfrak{b}$ as in b. $\square$

**Proposition 7.9.** *We have the following characterization of $\mathbf{No}_\lambda$.*

    a) *If $\nu$ is a non-zero limit, then a number a lies in $\mathbf{No}_\lambda$ if and only if every path in a is $\lambda$-bounded.*

    b) *If $\nu = \mu + 1$ is a successor, then a number a lies in $\mathbf{No}_\lambda$ if and only if path in a is either $\lambda$-bounded or satisfies*

$$\exists n \in \mathbb{Z}, (u_{P,|P|-1}, \alpha_{P,|P|-1}) = (L_\lambda \omega + n, \lambda).$$

    *for the first index $i < |P|$ with $\sigma_{P,i+1} = 1$ and ($\alpha_{P,i} \geqslant \lambda$ or $\beta_{P,i} \geqslant \lambda$).*

**Proof.** We write $\mathbf{No}_\lambda = \bigcup_{\rho \in \mathbf{On}} \mathbb{T}^{[\rho]}(\nu)$ as above, irrespective of the nature of $\nu$. Suppose that $\nu$ is a limit. Every path in elements of $\mathbb{L}_{<\lambda} \circ \omega$ is $\lambda$-bounded. We deduce with Proposition 7.3 that every path in $\mathbb{T}^{[0]}(\nu)$ is $\lambda$-bounded. It follows by induction using Proposition 7.3 that every path in elements of $\mathbf{No}_\lambda$ is $\lambda$-bounded.

Consider the property $\mathcal{X}_0$ of paths of being $\lambda$-bounded. Consider a $L_{<\lambda}$-atomic number $\mathfrak{a}$ all of whose paths satisfy $\mathcal{X}_0$. Write $\mathfrak{a} = L_\gamma E_\alpha^u$ as a hyperserial expansion, where $\gamma\omega < \alpha$. By [10, Lemma 5.5(a)], we have $\gamma \ggg \lambda_{/\omega} = \lambda$, and we either have $\alpha = 0$ or $\alpha \geqslant \lambda$. The path $(\mathfrak{a})$ is $\lambda$-bounded so we must have $\alpha = 0$ and $\gamma < \lambda$, whence $\gamma = 0$ and $\mathfrak{a} = \omega \in \mathbf{No}_\lambda$. Therefore Lemma 7.7 applies for $\mathcal{X}_0$ and entails that numbers with only $\lambda$-bounded paths lie in $\mathbf{No}_\lambda$. This concludes the proof of a).

Assume now that $\nu = \mu + 1$ is a successor and write $\beta := \omega^\mu$. Consider the property $\mathcal{X}_1$ of a path $P$ of being $\lambda$-bounded or satisfying

$$\exists n \in \mathbb{Z}, (u_{P,|P|-1}, \alpha_{P,|P|-1}) = (L_\lambda \omega + n, \lambda).$$

for the first index $i < |P|$ with $\sigma_{P,i+1} = 1$ and $\alpha_{P,i} \geqslant \lambda$ or $\beta_{P,i} \geqslant \lambda$. Let $P$ be a path in $\widetilde{\mathbb{L}_{<\lambda}} \circ \omega$ and assume that $P$ is not $\lambda$-bounded. Let $i < |P|$ be minimal to witness this. So $\sigma_{P,i+1} = 1$ and $\alpha_{P,i} \geqslant \lambda$ or $\beta_{P,i} \geqslant \lambda$. Then in particular $\mathfrak{a} := L_{\beta_{P,i}} E_{\alpha_{P,i}}^{u_{P,i+1}}$ is $L_{<\lambda}$-atomic, so there is an $n \in \mathbb{Z}$ with $\mathfrak{a} = L_{\beta n} \omega$. Therefore $\beta_{P,i} = 0$ and the hyperserial expansion of $\mathfrak{a}$ is $\mathfrak{a} = L_{\beta_{P,i}} E_{\alpha_{P,i}}^{u_{P,i+1}} = E_\lambda^{L_\lambda \omega + n}$ by Lemma 7.8. So $P$ satisfies $\mathcal{X}_1$. It follows by induction on $\rho$ using Proposition 7.3 that every path in elements of $\mathbf{No}_\lambda$ satisfies $\mathcal{X}_1$. As we've seen, each path in $L_{\beta n} \omega$ for $n \in \mathbb{Z}$ satisfies $\mathcal{X}_1$. We deduce with Lemma 7.6 that we can apply Lemma 7.7 for $\mathcal{X}_1$. So any path in an element of $\mathbf{No}_\lambda$ satisfies $\mathcal{X}_1$. $\square$

## 7.5 Right compositions with atomic elements

Let $\nu \in \mathbf{On}^>$, write $\lambda := \omega^\nu$ and fix an $\mathfrak{a} \in \mathbf{Mo}_\lambda$. We now define an embedding $\triangle_\mathfrak{a} : \mathbf{No}_\lambda \longrightarrow \mathbf{No}$ of force $\nu$ which is to be thought as the canonical composition on the right with $\mathfrak{a}$. Recall by Lemma 7.6 that $(\mathbf{Mo}(\lambda))_\lambda = \{\omega\}$ if $\nu$ is a limit and that $(\mathbf{Mo}(\lambda))_\lambda = \{E_{(\lambda_{/\omega})n} \omega : n \in \mathbb{Z}\}$ if $\nu$ is a successor. Consider the function $\triangle_\mathfrak{a} : (\mathbf{Mo}(\lambda))_\lambda \longrightarrow \mathbf{Mo}_\lambda$ defined by $\triangle_\mathfrak{a}(\omega) := \mathfrak{a}$ if $\nu$ is a limit, and $\triangle_\mathfrak{a}(E_{(\lambda_{/\omega})n} \omega) := E_{(\lambda_{/\omega})n} \mathfrak{a}$ for all $n \in \mathbb{Z}$ if $\nu$ is a successor. If $\Sigma$ is a $\lambda$-bounded nested sequence over $\mathbf{No}_\lambda$, then we define $(\triangle_\mathfrak{a})_\Sigma$ to be the identity map on $\mathbf{No}$. The family $(\triangle_\mathfrak{a}, ((\triangle_\mathfrak{a})_\Sigma))$ satisfies the premises of Theorem 6.3. So $\triangle_\mathfrak{a}$ extends uniquely into a hyperserial embedding

$$\triangle_\mathfrak{a} : \mathbf{No}_\lambda \longrightarrow \mathbf{No}$$

of force $\nu$. We call $\triangle_{\mathfrak{a}}$ the *canonical right composition with* the atomic element $\mathfrak{a}$.

**Example 7.10.** Suppose that $\nu=1$ and $\mathfrak{a}=\log \omega$. Then $\triangle_{\mathfrak{a}}$ should be thought as a formal substitution of $\log \omega$ for $\omega$ in hyperserial representations (see [10, Section 7]) of numbers in $\mathbf{No}_\omega$. The unicity property of Theorem 6.3 entails that for $n \in \mathbb{N}$, the $n$-fold iterate of $\triangle_{\mathfrak{a}}$ (taking inverses if $n<0$) is $\triangle_{L_n(\omega)}$.

# Bibliography

[1] N. L. Alling. *Foundations of analysis over surreal number fields*, volume 141. North-Holland Mathematical Studies, North-Holland Amsterdan, 1987.

[2] M. Aschenbrenner, L. van den Dries, and J. van der Hoeven. *Asymptotic Differential Algebra and Model Theory of Transseries*. Number 195 in Annals of Mathematics studies. Princeton University Press, 2017.

[3] M. Aschenbrenner, L. van den Dries, and J. van der Hoeven. On numbers, germs, and transseries. In *Proc. Int. Cong. of Math. 2018*, volume 1, pages 1–24. Rio de Janeiro, 2018.

[4] M. Aschenbrenner, L. van den Dries, and J. van der Hoeven. The surreal numbers as a universal H-field. *Journal of the European Mathematical Society*, 21(4), 2019.

[5] M. Aschenbrenner, L. van den Dries, and J. van der Hoeven. Maximal Hardy fields. https://arxiv.org/abs/2304.10846, 2023.

[6] V. Bagayoko. Hyperexponentially closed fields. https://www.hal.inserm.fr/X-LIX/hal-03686767v1, 2022.

[7] V. Bagayoko. *Hyperseries and surreal numbers*. PhD thesis, UMons, École Polytechnique, 2022. https://theses.hal.science/tel-04105359.

[8] V. Bagayoko. On ordered groups of regular growth rates. https://arxiv.org/abs/2402.00549, 2024.

[9] V. Bagayoko and J. van der Hoeven. The hyperserial field of surreal numbers. https://hal.science/hal-03232836, 2021.

[10] V. Bagayoko and J. van der Hoeven. Surreal numbers as hyperseries. https://hal.science/hal-03681007, 2022.

[11] V. Bagayoko, J. van der Hoeven, and E. Kaplan. Hyperserial fields. https://hal.science/hal-03196388, 2021.

[12] A. Berarducci and V. L. Mantova. Surreal numbers, derivations and transseries. *JEMS*, 20(2):339–390, 2018.

[13] A. Berarducci and V. L. Mantova. Transseries as germs of surreal functions. *Transactions of the American Mathematical Society*, 371:3549–3592, 2019.

[14] J. H. Conway. *On numbers and games*. Academic Press, 1976.

[15] L. van den Dries, J. van der Hoeven, and E. Kaplan. Logarithmic hyperseries. *Transactions of the American Mathematical Society*, 372, 2019.

[16] J. Écalle. *Introduction aux fonctions analysables et preuve constructive de la conjecture de Dulac*. Actualités Mathématiques. Hermann, 1992.

[17] J. Écalle. The natural growth scale. *Journal of the European Mathematical Society*, CARMA, vol 1:93–223, 2020.

[18] P. Ehrlich. Absolutely Saturated Models. *Fundamenta Matematicae*, 133(1):39–46, 1989.

[19] H. Gonshor. *An Introduction to the Theory of Surreal Numbers*. Cambridge Univ. Press, 1986.

[20] H. Hahn. Über die nichtarchimedischen Grö$\beta$ensysteme. *Sitz. Akad. Wiss. Wien*, 116:601–655, 1907.

[21] J. van der Hoeven. *Automatic asymptotics*. PhD thesis, École polytechnique, Palaiseau, France, 1997.

[22] J. van der Hoeven. Operators on generalized power series. *Illinois Journal of Math*, 45(4):1161–1190, 2001.

[23] J. van der Hoeven. *Transseries and real differential algebra*, volume 1888 of *Lecture Notes in Mathematics*. Springer-Verlag, 2006.

[24] S. MacLane. The universality of formal power series fields. *Bulletin of the American Mathematical Society*, 45(12.P1):888–890, 1939.

[25] M. Rosenlicht. Hardy fields. *Journal of Mathematical Analysis and Applications*, 93(2):297–311, 1983.

[26] M. C. Schmeling. *Corps de transséries*. PhD thesis, Université Paris-VII, 2001.